\documentclass[11pt,a4paper]{article}

\usepackage{auxhook}
\usepackage{latexsym}
\usepackage{blindtext}

\PassOptionsToPackage{force}{filehook}

\usepackage[T1]{fontenc}
\usepackage[utf8]{inputenc}
\usepackage[main=english,french]{babel}
\usepackage{graphicx}
\usepackage{graphics,subfigure}
\usepackage{fancyhdr}
\usepackage{lmodern}
\usepackage{color}
\usepackage{transparent}
\usepackage{float}
\usepackage{pythontex}
\usepackage{tikz-cd}
\usepackage{listings}
\usepackage[linktocpage=true,citebordercolor=green, urlcolor=black]{hyperref}
\usepackage{tabularx}
\usepackage[title,titletoc]{appendix}
\usepackage{xcolor}
\usepackage[mathscr]{eucal}
\usepackage{calrsfs}
\usepackage[skins]{tcolorbox}

\usetikzlibrary{shapes.misc}
\usepackage{pgf,tikz}
\usepackage{caption}
\usetikzlibrary{decorations.markings}
\usetikzlibrary{arrows}
\usetikzlibrary{patterns}
\usepackage{tikz-3dplot}
\usetikzlibrary{calc}
\usepackage{amsfonts}
\usepackage{amsmath,amssymb,amsthm}
\usepackage{mathtools}
\usepackage{relsize}
\usepackage{titlesec}
\usepackage{tcolorbox}

\newcolumntype{P}[1]{>{\centering\arraybackslash}p{#1}}

\usepackage{stmaryrd}
\usepackage{trimclip}

\usepackage[shortlabels]{enumitem}

\usepackage{wasysym}

\makeatletter
\DeclareRobustCommand{\shortto}{%
	\mathrel{\mathpalette\short@to\relax}%
}

\newcommand{\short@to}[2]{%
	\mkern2mu
	\clipbox{{.5\width} 0 0 0}{$\m@th#1\vphantom{+}{\shortrightarrow}$}%
}
\makeatother

\theoremstyle{definition}
\newtheorem{definition}{Definition}[section]

\newtheorem{remark}[definition]{Remark}
\newtheorem{corollary1}[definition]{Corollary}
\newtheorem{lemma1}[definition]{Lemma}

\newtheorem*{proposition*}{Proposition}
\newtheorem*{conjecture*}{Conjecture}

\titleformat{\paragraph}
{\normalfont\normalsize\bfseries}{\theparagraph}{1em}{}
\titlespacing*{\paragraph}
{0pt}{3.25ex plus 1ex minus .2ex}{1.5ex plus .2ex}

\DeclareMathOperator{\Hess}{Hess}
\DeclareMathOperator{\grad}{grad}
\DeclareMathOperator{\Tr}{Tr}

\usepackage[margin=2.5cm]{geometry}

\usepackage{mathptmx}

\usepackage{authblk}

\usepackage{comment}
\usepackage[strict]{csquotes}
\usepackage[
style=numeric,
autolang=other,
backend=biber,
sorting=nyt,
citereset=none,
maxnames=10,
natbib=true,
hyperref=true,
url=true,
giveninits=true,
]{biblatex}

\title{Studying network of symmetric periodic orbit families of the Hill problem via symplectic invariants}
\author[1]{Cengiz Aydin}
\author[2]{Alexander Batkhin}
\affil[1]{Institut für Mathematik, Universität Heidelberg, Germany}
\affil[ ]{\textit {E-mail address: \href{mailto:cengiz.aydin@hotmail.de}{cengiz.aydin@hotmail.de}}}
\affil[2]{Technion -- Israel Institute of Technology, Haifa, Israel}
\affil[ ]{\textit {E-mail address: \href{mailto:batkhin@technion.ac.il}{batkhin@technion.ac.il}}}

\begin{document}
	
\setcounter{page}{1}
\pagenumbering{arabic}
	
\maketitle

\begin{abstract}
In the framework of the spatial circular Hill three-body problem we illustrate the application of symplectic invariants to analyze the network structure of symmetric periodic orbits families.\ The extensive collection of families within this problem constitutes a complex network, fundamentally comprising the so-called basic families of periodic solutions, including the orbits of the satellite $g$, $f$, the libration (Lyapunov) $a,c$, and collision $\mathcal B_0$ families.\ Since the Conley--Zehnder index leads to a grading on the local Floer homology and its Euler characteristics, a bifurcation invariant, the computation of those indices facilitates the construction of well-organized bifurcation graphs depicting the interconnectedness among families of periodic solutions.\ The critical importance of the symmetries of periodic solutions in comprehending the interaction among these families is demonstrated.
\end{abstract}

\begin{center}
\begin{tabular}{ll}
	\textbf{Keywords} & $\quad$ Hill three-body problem $\cdot$ periodic solution $\cdot$ Conley--Zehnder index $\cdot$\\
    & $\quad$ symmetry $\cdot$ bifurcation diagram \\
	\textbf{MSC 2020} & $\quad$ 70G45 $\cdot$ 70F07 $\cdot$ 70H12
\end{tabular}
\end{center}
	
\makeatletter{\renewcommand*{\@makefnmark}{}

\setcounter{tocdepth}{2}
\tableofcontents

\section{Introduction}

The Hill three-body problem, a limiting case of the circular restricted three-body problem, is a well-known scenario in celestial mechanics which provides an approximation of the dynamics of the infinitesimal body in the vicinity of the smaller primary.\ The remaining primary is pushed infinitely far away in a way that it acts as a velocity independent gravitational perturbation of the rotating Kepler problem formed by the smaller primary and the infinitesimal body.\ In its original application, Hill \cite{hill} reformulated the lunar theory in which one of his main contribution was his discovery of one periodic solution with period equal to the synodic month of the Moon.\ Further applications of Hill's approach were captured in the dynamics of natural or artificial satellites \cite{henonVI, lidov}, distant moons of asteroids \cite{hamilton_krivov}, low-energy escaping trajectories \cite{villac_scheeres} or frozen orbits around planetary satellites \cite{lara_palacian}, and dynamics of star clusters~\cite{Gurfil2005}.\ Extensions of Hill's concept to four-body systems have also been developed and studied \cite{scheeres, burgos_gidea, aydin_contact}.\

Since Poincaré's pioneering work \cite{Poincare1893} periodic solutions are used as the central structure for the study of the global dynamics.\ Therefore, it becomes of vital importance to catalog as many periodic orbits as possible.\ The total energy of the system is a first integral, thus by the implicit function theorem, non-degenerate periodic orbits always come in a smooth one-parameter family and hence form a smooth orbit cylinder, parameterized by the energy \cite[Section~6.4]{meyer}.\ Natural families of periodic orbits are those that can be approached via solvable systems and analytical considerations, either from the Kepler problem for very low energies (planar direct, or pro-grade, and retrograde infinitesimal circular orbits, and rectilinear vertical consecutive infinitesimal collision orbits in regularized system) or from the linear behavior of the flow around the collinear Lagrange points (planar and vertical Lyapunov orbits).\ All families of periodic solutions can be conditionally divided into two groups: \textit{open families}, continuing up to the limit values, and \textit{closed families}, existing only on a finite interval of the energy values.

An overwhelming amount of significant results on the study of basic periodic orbits can be found in the literature, summarized in chronological order:
\begin{itemize}[noitemsep]
    \item In 1969 Hénon~\cite{henonV} gave a first numerical exploration of the basic families of planar Lyapunov (family $a$ and $c$), direct (family $g$), and retrograde orbits (family $f$), and studied their horizontal stability properties.\ While the orbits of the family $g$ and $f$ are doubly symmetric, the orbits of the family $a$ are simple symmetric.\ Especially, the family $g$ undergoes a planar symmetry-breaking pitchfork bifurcation, where two families of simple symmetric direct orbits appear (family $g'$).\  Also in this paper was provided a first attempt to describe the known families of periodic solutions in terms of generating solutions obtained from the limiting problem. This limiting problem is the planar case of the Clohessy--Wiltshire equations~\cite{Clohessy1960}.\ The complete description of this approach was provided in his latter papers in 2003 and 2005 (see below).\ A year later the paper~\cite{henonVI} with description of some planar Hill's problem applications to natural and artificial satellites movement, stars dynamics in clusters were published.\
    \item The vertical stability of the above basic planar orbits were studied by Hénon~\cite{henon74} in 1974, in which he characterized their vertical critical orbits, up to period-doubling.\    
    \item Based on the latter work, in 1980 Michalodimitrakis~\cite{michalodimitrakis} examined some out-of-plane bifurcations and described three branches from family $a$:\ family $a_{1v}$ (its orbits are nowadays known as halo orbits), $a_{2v}$ and $a_{3v}$ (period-doubling).\ He found out that halo orbits terminate at a rectilinear vertical collision orbit (that is the family $\mathcal{B}_0^{\pm}$, see below), and that the family $a_{2v}$ terminates at a spatial orbit, which is doubly symmetric and whose projection onto the $yz$-plane is of ``figure-eight'' shape (that is a vertical Lyapunov orbit, see below).\ In particular, the latter closed family forms a ``bridge'' (i.e.,\ an orbit cylinder, as defined in \cite{llibre_meyer_soler}).\ He also computed period-doubling bifurcations from $g$ and $g'$:\ families $g_{1v}$, $g'_{1v}$ and $g'_{2v}$.\ The first two families are identical, closed and form a connection between $g$ and $g'$; the latter is open which is also described in \cite{batkhin2009}.\
    \item In 1981 Perko applying technique of generating solutions proved the existences of the countable set of families of periodic orbits with one symmetry~\cite{PerkoI,PerkoII}. He also provided some asymptotic forms of the initial conditions of families $f$, $g$, $a$ ($c$), $g'$ and predict some new families, but the proposed formulas turned out to be unsuitable for computations.\ 
    \item The first attempt to obtain generating solution of families of spatial periodic orbits interacting with the family of rectilinear vertical consecutive collision orbits (not starting out of plane ones), denoted by family $\mathcal{B}_0^{\pm}$ ($\mathcal{B}_0^+$ indicates collision from above and $\mathcal{B}_0^-$ indicates collision from below), was provided by Lidov and Rabinovich 1979~\cite{Lidov_Rabinovich_1979KosIs}.\ This family was also discussed in a recent work\cite{belbruno_frauenfelder_vankoert} and can be found in the book~\cite[Subsect.~8.4.2]{marchal}, without any description of it.\ For a comprehensive study of such families in 1982 Lidov worked out an original numerical method based on KS regularization of the Hill3BP problem~\cite{lidov82a,lidov82}.\ The results of the application of this method was described in~\cite{lidov_lyakhova83}, in which eight families of symmetric spatial periodic solutions were numerically investigated with various level of details.\ One of the discovered family was subsequently used to design the orbit of a radio interferometer space station~\cite{Lidov_Lyakhova_1988PAZh,Lidov_Lyakhova_1988SvAL}.\ 
    Since the mantioned above papers by Lidov and his coauthors have not been translated into English, they were not well known, hence we will give a short description about $\mathcal{B}_0^{\pm}$ family in our paper.\
    \item In 1985 Zagouras--Markellos \cite{zagouras_markellos} studied the vertical Lyapunov orbits (family $_H L^e_{2v}$).\ They discovered that the vertical Lyapunov orbits undergo a symmetry-breaking pitchfork bifurcation, which corresponds exactly to the bridge mentioned above (family $a_{2v}$).
    \item In 1998 Sim\'o and Stuchi~\cite{SIMO20001} provided a qualitative description of the global dynamics of the planar Hill problem and also mentioned a lot of families of periodic solutions playing the role bridges between families of direct and retrograde orbits.
    \item Hénon studied in one of his last published works~\cite{Henon2003} from 2003 the horizontal multiple cover bifurcations from basic planar orbits.\ In this paper the theory of generating solutions in the context of the planar Hill problem, worked out for the restricted three-body problem in H\'enon's book~\cite[Ch.~5]{Henon97}, was presented.\ In particular, he described doubly-symmetric triple-revolution orbits branching out from $f$ and denoted this family by $g_3$.\ This family was found by T.\,Matukuma in 1957 and was firstly described by Hénon \cite{henonVI} in 1970.\ In 1994 A.D.\,Bruno mentioned this family in one of his preprints where it was denoted as $f_3$ due to that fact that this family intersects twice with the~$f$. More detailed description of this family was provided in~\cite[\S 4.3.5]{Hillbook}.\ In this paper we also use the last notation of this family as $f_3$.\ The results on families of asymmetric periodic solutions were given in~\cite{Henon2005}.\ Up to now this is the only paper on the Hill's problem asymmetric solutions.
    \item Paper~\cite{gomez_marcote_mondelo} by Gomez--Marcote--Mondelo presented results on the invariant manifolds of the spatial Hill’s problem associated to the two liberation points.\ A systematic exploration of the homoclinic and heteroclinic connections between the center manifolds of the liberation points was also given by applying a specific normalization of the problem's Hamiltonian.
    \item In 2009 Batkhin--Batkhina~\cite{batkhin2009} studied out-of-plane bifurcations from basic families and formed a network relating those basic orbits with the determined new branches.\ Especially, it was worked out that $g'$ is related to double cover of halo orbits, and $g$ is connected to double cover of $\mathcal{B}_0^{\pm}$.\
    \item There is a paper~\cite{Tsiroganis2009}, where by the brute force method 31 families of planar periodic solutions making different covering of the family $f$ were obtained.
    \item In 2014 the second author~\cite{BatkhinDAN14} proposed a certain generalization of the Hill problem.\ This generalization demonstrated that all the known families of planar periodic solutions form a common network interacting to each other either by $k$ covering or by sharing the same generating solutions.
    \item Similar to the paper~\cite{batkhin2009}, in 2020 Kalantonis~\cite{kalantonis} studied vertical period-tripling and -quadrupling bifurcation from families $g$, $g'$ and $a$.\ In particular, he could derive connections (closed families) between $g^3$ and $f^5$, $g^4$ and $g'^4$, and $g'^4$ and $f^6$, where the upper index corresponds to number of covering where vertical bifurcation occurs.\
\end{itemize}

Such numerical explorations show that various families of periodic orbits can interact with each other at bifurcation points.\ Connections of this kind and their data base are indeed hard to organize by pure computations.\ For the natural question to ponder, \textit{``Which families can interact with each other at bifurcation points?''},\ classical methods (such as continuation method, stability indices) are not enough, hence additional tools and structures are necessary.\ This problem was addressed by Aydin in his recent work \cite{aydin_cz}, by applying symplectic techniques.\ One of the key features is to consider the Conley--Zehnder index of periodic orbits and study its interaction at bifurcation points.\ The Conley--Zehnder index \cite{conley_zehnder} assigns a mean winding number for the linearized flow of non-degenerate periodic orbits.\ On an orbit cylinder, all periodic orbits have the same index, and at bifurcation points the index jumps.\ When working locally near a family of non-degenerate periodic orbits, then there is a fascinating bifurcation invariant \cite{ginzburg}:\ the local Floer homology and thus its Euler characteristic, the alternating sum of the ranks of the local Floer homology groups.\ Significantly, the index leads to a grading on local Floer homology and thus, the index provides important information how different families are related to each other at bifurcation points.\ By using such symplectic invariants, Aydin \cite{aydin_cz} constructed bifurcation graphs demonstrating a well-organized picture of the interconnectedness associated to $g$, $g'$ and $f$, and their multiple cover bifurcations.\ By a ``bifurcation graph'' we understand a labeled graph, whose vertices correspond to bifurcation points and whose edges correspond to periodic orbits families, labeled with Conley--Zehnder indices.\ This symplectic approach provides additional structure to periodic orbit families and their relations at bifurcation points from a topological point of view.\ It is also worth mentioning that Aydin showed in \cite{aydin_babylonian} that Conley--Zehnder indices are astronomically significant.\ He used those indices and the Floquet multipliers of Hill's lunar orbit with a period of one synodic month to express the anomalistic and draconitic months of the Moon.\

Our purpose in this paper is to provide bifurcation graphs which illustrate a common network in association to the natural periodic orbit families and their bifurcations.\ Table \ref{table_1} shows our main bifurcation results in each row, where in each row the integer $n$ indicates each $n$-th cover bifurcation of the underlying family in the first row.\ In particular, our results highlight that spatial bifurcation connections related to $n$-th cover of $g$, $n+1$-th cover of $\mathcal{B}_0^{\pm}$ and $n+2$-th cover of $f$ are provided for $n = 3,4,5$.\ One of the aims of this work is to demonstrate that by Conley--Zehnder indices one can expect such pattern.\

\begin{table}[bh]\centering
\begin{tcolorbox}[width=6.8cm]
\begin{tabular}{P{0.5cm}P{0.5cm}P{0.5cm}P{0.5cm}P{0.5cm}P{0.5cm}}
     $g'$ & $g$ & $\mathcal{B}_0^{\pm}$ & $f$ & $f_3$ & halo
\end{tabular}
\tcblower
\begin{tabular}{P{0.5cm}P{0.5cm}P{0.5cm}P{0.5cm}P{0.5cm}P{0.5cm}}
     1 & & & & & 2 \\
     & 1 & 2 & & & \\
     2 & 2 & 3 & & & \\
     3 & 3 & 4 & 5 & 1 & \\
     4 & 4 & 5 & 6 & 2 & \\
     & 5 & 6 & 7 & &
\end{tabular}
\end{tcolorbox}
\caption{Bifurcation results.}
\label{table_1}
\end{table}

This paper is organized as follows.\ In Section \ref{sec2} we briefly discuss the equations for the spatial circular restricted three-body problem and the Hill three-body problem (from now on, SCR3BP and Hill 3BP), their linear symmetries and Lagrange points.\ Moreover, we give the description of KS regularization of the Hill 3BP together with monodromy matrix computation in regularized variables.\ Here also a brief description of Lidov's method of generating solutions computation is provided.\ We also discuss some peculiarities of branching of families of spatial periodic solutions from planar doubly symmetric critical orbits and vertical collision orbits as well.\ In Section~\ref{sec3} we give an overview of how the Conley--Zehnder index interacts with bifurcation points and describe the relevant computational tools.\ Section \ref{sec4} is devoted to the Conley--Zehnder indices of natural periodic orbit families at their origin.\ In Section \ref{sec5} we provide detailed descriptions of our results.\ Finally Appendix contains detail information in the form tables of data.

\section{SCR3BP and Hill 3BP}
\label{sec2}

\subsection{Equations of motion}

In the SCR3BP we consider two primary masses $B_1$ and $B_2$, orbiting on circles around their common barycenter following the Newton’s law, and a third body $B_3$, that is significantly smaller than the other two and thus has a negligible effect on their motion.\ It is convenient to choose the units of mass, length and time such that the gravitational constant is 1, and the period of the motion of the primaries is $2\pi$.\ In these units the distance between the two primaries is normalized to be one.\ Furthermore, it is convenient to use a rotating frame that rotates with an angular velocity of the orbital angular rate of the primaries.\ Then, if we denote by $\mu \in [0,\frac{1}{2}]$ the mass of $B_2$, the primary $B_1$ is fixed at $(-\mu,0,0)$ and $B_2$ is fixed at $(1-\mu,0,0)$.\ The motion of $B_3$ is described by the second order differential equations \cite{szebehely}:\
\begin{align}
	\ddot{x} &= 2 \dot{y} + x - (1 - \mu) \dfrac{x + \mu}{r_1^3} - \mu \dfrac{x - 1 + \mu}{r_2^3}, \nonumber \\
	\ddot{y} &= - 2 \dot{x} + y - \left( \dfrac{1 - \mu}{r_1^3} + \dfrac{\mu}{r_2^3} \right) y, \label{ham_equation1} \\
	\ddot{z} &= - \left( \dfrac{1-\mu}{r_1^3} + \dfrac{\mu}{r_2^3} \right) z, \nonumber
\end{align}
where $r_1 = \left( (x+\mu)^2 + y^2 + z^2 \right)^{\frac{1}{2}}$ and $r_2 = \left( (x - 1 + \mu)^2 + y^2 + z^2 \right)^{\frac{1}{2}}$ indicate the distances from $B_3$ to the primaries.\ In most astronomical applications one of the primaries has a small mass $\mu$ compared to the second one \cite[p.\ 5]{restrepo_russell}, e.g., $\mu = 3.0404317 \times 10^{-6}$ corresponds to the Sun--Earth system.\

Defining kinetic moments by $p_x = \dot{x} - y$, $p_y = \dot{y} + x$ and $p_z = \dot{z}$, the system can be written in Hamiltonian form with the corresponding Hamiltonian function
\begin{align} \label{hamiltonian}
	H(x,y,z,p_x,p_y,p_z) = \frac{1}{2} \left( p_x^2 + p_y^2 + p_z^2 \right) - \frac{1-\mu}{r_1} - \frac{\mu}{r_2} + p_x y - p_y x,
\end{align}
which is a first integral of the SCR3BP.\ An equivalent first integral is the Jacobi integral defined by $\Gamma = -2H$, where $\Gamma$ is called the Jacobi constant.\ Notice that the phase space  $\left( \mathbb{R}^3 \setminus \{ B_1,B_2 \} \right) \times \mathbb{R}^3$ is endowed with the standard symplectic form $\omega = \sum k \wedge p_k$ ($k=x,y,z$) and the flow of its Hamiltonian vector field $X_H$, defined by $dH(\cdot) = \omega( X_H , \cdot )$, is equivalent to the equations~\eqref{ham_equation1}.\

In the Hill 3BP one of the primaries is significantly much heavier than the other one and the infinitesimal body moves in the vicinity of the smaller primary.\ For very small $\mu$, one shifts the smaller primary to the origin, scales the coordinates by a factor $\mu^{\frac{1}{3}}$ and after a Taylor expansion \cite{meyer_schmidt} the Hamiltonian becomes
\begin{equation*}\label{eq:HillHam_mu} 
H^{\mu}(x,y,z,p_x,p_y,p_z) = \frac{1}{2} \left( p_x^2 + p_y^2 + p_z^2 \right) - \frac{1}{r} + p_x y - p_y x - x^2 + \frac{1}{2}y^2 + \frac{1}{2}z^2 + \mathscr{O} ( \mu^{\frac{1}{3}} ), 
\end{equation*}
where $r = \left( x^2 + y^2 + z^2 \right)^{\frac{1}2}$.\ Since $\mu$ is very small, all the $\mathscr{O} ( \mu^{\frac{1}{3}} )$ terms can be neglected and the limiting Hamiltonian for $\mu \to 0$ corresponds to the Hill 3BP.\ Dynamically, Hill's Hamiltonian consists of the rotating Kepler problem (that is the case of the SCR3BP for $\mu = 0$) with a velocity independent gravitational perturbation produced by the massive primary (the degree 2 term $- x^2 + \frac{1}{2}y^2 + \frac{1}{2}z^2$).\ This difference between the rotating Kepler problem and Hill's system gives a dramatic dynamical change.\ While the rotating Kepler problem is an integrable system, the Hill 3BP is non-integrable \cite{meletlidou}, \cite{morales-ruiz}, whose equations are given by
\begin{equation}\label{ham_equation2}
\begin{aligned}
\ddot{x} &= 2 \dot{y} + 3x - \frac{x}{r^3}, \\	\ddot{y} &= - 2 \dot{x} - \frac{y}{r^3},  \\
\ddot{z} &= - \left( \frac{1}{r^3} + 1 \right) z.
\end{aligned}
\end{equation}

\subsection{Linear symmetries and Lagrange points}

A ``symmetry'' $\sigma$ is an involution on the phase space which leaves the system invariant and is symplectic or anti-symplectic ($\sigma ^* \omega = \pm \omega$).\ Anti-symplectic symmetries are used to denote time-reversal symmetries in the Hamiltonian context, see e.g., \cite{lamb_roberts}.\ A periodic solution $\mathbf q \equiv (x,y,z,p_x,p_y,p_z)$ is symmetric with respect to an anti-symplectic symmetry $\rho$ if $\mathbf q (t) = \rho \left( \mathbf q (-t) \right)$ for all $t$, and symmetric with respect to a symplectic symmetry $\sigma$ if $\mathbf q (t) = \sigma \left( \mathbf q (t) \right)$ for all $t$.\

\begin{figure}[t!]
	\centering
	\includegraphics[scale=1]{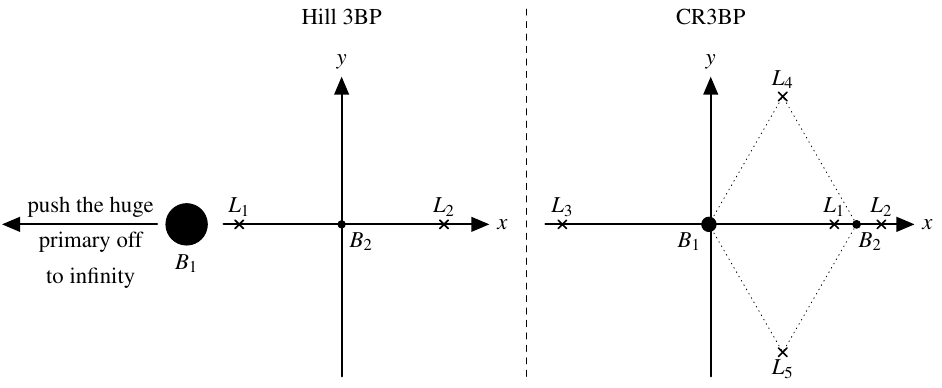}
	\caption{Right:\ Lagrange points in the CR3BP; left: Lagrange points in the Hill 3BP.}
	\label{figure_lagrange_points}
\end{figure}

In the SCR3BP, the reflection at the ecliptic $\{z=0\}$ gives rise to the linear symplectic symmetry
\begin{align}\label{sigma}
			\sigma (x,y,z,p_x,p_y,p_z) = (x,y,-z,p_x,p_y,-p_z),
\end{align}
whose fixed point set Fix$(\sigma) = \{ (x,y,0,p_x,p_y,0) \}$ corresponds to the planar system.\ By the third equation in \eqref{ham_equation1}, it is obvious that all Lagrange points are located at the ecliptic.\ There are five Lagrange points \cite{szebehely} (see right part in Figure \ref{figure_lagrange_points}):\ Three collinear points $L_1$ (located within the primaries), $L_2$ and $L_3$ (outside the interval joining the primaries) on the $x$-axis determined by the equations
\begin{align}\label{x_coord_langrange}
	x_{L_j} + A \frac{1-\mu}{(x_{L_j} + \mu)^2} + B \frac{\mu}{(x_{L_j} - 1 + \mu)^2} = 0,\quad j=1,2,3,
\end{align}
where $(A,B) = (-1,+1)$ if $j=1$, $(A,B) = (-1,-1)$ if $j=2$ and $(A,B) = (+1,+1)$ if $j=3$.\ In addition, there are two equilateral triangular points $L_4$ and $L_5$, determined by $x_{L_4} = x_{L_5} = \frac{1}{2} - \mu$ and $y_{L_j} = C \frac{\sqrt{3}}{2}$, where $C=+1$ if $j=4$ and $C=-1$ if $j=5$.\ The two equilateral triangular points $L_4$ and $L_5$ are related by the planar system's anti-symplectic symmetry given by
\begin{align*}
\rho^p_{OX} (x,y,0,p_x,p_y,0) = (x,-y,0,-p_x,p_y,0) \text{ (reflection at the $x$-axis)},
\end{align*}
which can be extended to two anti-symplectic symmetries,
\begin{align}
XOZ := \rho_1 (x,y,z,p_x,p_y,p_z) = &(x,-y,z,-p_x,p_y,-p_z) \text{ (reflection at the $xz$-plane)},\label{eq:rho1}\\
OX := \rho_2 (x,y,z,p_x,p_y,p_z) = &(x,-y,-z,-p_x,p_y,p_z) \text{ ($\pi$-rotation around the $x$-axis)}.\label{eq:rho2}
\end{align}
Their composition is $\rho_1 \circ \rho_2 = \rho_2 \circ \rho_1 = \sigma$.\

The planar Hill 3BP is additionally invariant under the anti-symplectic symmetry
\begin{align}\label{rho_p_2}
\rho^p_{OY} (x,y,0,p_x,p_y,0) = (-x,y,0,p_x,-p_y,0) \text{ (reflection at the $y$-axis)}.
\end{align}
From the equations~\eqref{ham_equation2} one readily calculates that the Hill 3BP only has the first two Lagrange points $L_1$ and $L_2$ with equal distance to the origin (see left part in Figure \ref{figure_lagrange_points}).\ They are located at $(\pm 3^{-1/3},0,0)$ with Jacobi constant $3 \sqrt[3]{3}$.\ The remaining Lagrange points $L_3$, $L_4$ and $L_5$ have been removed to infinity by the limiting procedure.\ Moreover, one can extend~\eqref{rho_p_2} to two further anti-symplectic symmetries in the spatial problem,
\begin{align}
YOZ := \rho_3 (x,y,z,p_x,p_y,p_z) = &(-x,y,z,p_x,-p_y,-p_z) \text{ (reflection at the $yz$-plane)},\label{eq:rho3}\\
OY := \rho_4 (x,y,z,p_x,p_y,p_z) = &(-x,y,-z,p_x,-p_y,p_z) \text{ ($\pi$-rotation around the $y$-axis)},\nonumber
\end{align}
Their composition is $\rho_3 \circ \rho_4 = \rho_4 \circ \rho_3 = \sigma$ as well.\ Other linear symplectic symmetries are $- \sigma$ and $\pm$id, where $- \sigma$ corresponds to a $\pi$-rotation around the $z$-axis, hence the $z$-axis is invariant under $-\sigma$.\ These eight form a group isomorphic to $\mathbb{Z}_2 \times \mathbb{Z}_2 \times \mathbb{Z}_2$, and in the planar problem the set $\{ \rho^p_{OX/OY} , \pm \text{id} \}$ form a Klein four-group, $\mathbb{Z}_2 \times \mathbb{Z}_2$.\ It is proved in \cite{aydin_sym} that there no other linear symmetries.\

\subsection{Kustaanheimo--Stiefel regularization of Hill 3BP}\label{sec:KS}

Connections between families of planar and spatial periodic solutions cannot be understood without the family $\mathcal B_0$ of vertical rectilinear collision orbits.\ Such orbits form two mutually symmetric families denoted $\mathcal B_0^+$ for upper half-space and $\mathcal B_0^-$ for lower half-space correspondingly.

The authors failed to prioritize the study of the orbits of the $\mathcal B_0$ family.\ A brief description of known publications, especially those that were published in the USSR in Russian and never been translated into English, is given in the historical review in Introduction.

\subsubsection{Regularization of equations of motion}
Considering the fact that all orbits of the families $\mathcal B_0^\pm$ are collisional, their dynamical properties can be studied only after regularization of the equations of motion.\ In the opinion of the second author the most elegant way to obtain Kustaanheimo-Stiefel (KS) regularization is a quaternion approach proposed by Waldvogel~\cite{Waldvogel2008}. 

Let $\mathbf u=u_1+\mathfrak i u_2+\mathfrak j u_2+\mathfrak k u_3$ be an arbitrary quaternion, where $u_k\in\mathbb{R}$, $k=1,2,3,4$, and $\mathfrak i, \mathfrak j, \mathfrak k$ are imaginary units satisfying the well-known Hamilton's identity 
\[
\mathfrak i^2 = \mathfrak j^2 = \mathfrak k^2 =\mathfrak i \mathfrak j \mathfrak k = -1.
\]
Physical coordinates are represented as a ``truncated'' quaternion $\tilde{\mathbf{r}}=x+\mathfrak i y+\mathfrak j z$ associated with the positional vector $\mathbf r=(x,y,z)$.\ To define mapping $\mathbf u\to\tilde{\mathbf{r}}$ a \textit{star conjugate} transformation $\mathbf u^*=u_1+\mathfrak i u_2+\mathfrak j u_2-\mathfrak k u_3$ is introduced. Star conjugation agrees with the basic operations on quaternions and norm:
\[
\left(\mathbf u^*\right)^*=\mathbf u, \quad (\mathbf u \mathbf v)^*=\mathbf v^* \mathbf u^*,\quad |\mathbf u^*| = |\mathbf u|,
\]
and provides the mapping 
\begin{equation}\label{eq:mapping}
    \mathbf u \to \tilde{\mathbf{r}}= \mathbf u \mathbf u^*.
\end{equation}
The norms of a positional vector $\mathbf r$ and a corresponding quaternion $\mathbf u$ satisfy relation $r\equiv|\mathbf r|=|\mathbf u|^2$.

To obtain the inverse mapping it is necessary to make the following trick.\ Firstly we find a particular solution $\tilde{\mathbf u}\tilde{\mathbf u}^*=\tilde{\mathbf{r}}$ in the form of a quaternion with vanishing $\mathfrak k$-component: $\tilde{\mathbf u}=\tilde u_1+\mathfrak i \tilde u_2+\mathfrak j \tilde u_2$.\ This particular solution can be written in the well-known form of the square root of a quaternion:
\begin{equation}\label{eq:reg_part_sol}
\tilde{\mathbf u}=\frac{\tilde{\mathbf{r}}+|\tilde{\mathbf{r}}|}{\sqrt{2\left(|\tilde{\mathbf{r}}|+x\right)}}.
\end{equation}
The entire family of solutions is obtained from particular~\eqref{eq:reg_part_sol} by multiplying by $e^{\mathfrak k\varphi}$:
\begin{equation}\label{eq:reg2phy}
\mathbf u = \tilde{\mathbf u}\cdot e^{\mathfrak k\varphi}=\tilde{\mathbf u}\left(\cos\varphi+\mathfrak k\sin\varphi\right).   
\end{equation}
For $\tilde{\mathbf{r}}$ with negative real part to avoid singularity at the point $\tilde{\mathbf{r}}=x+0\cdot\mathfrak i+0\cdot\mathfrak j$ another branch of the square root should be chosen:
\[
\tilde{\mathbf u}=\frac{\tilde{\mathbf{r}}-|\tilde{\mathbf{r}}|}{\mathfrak i\sqrt{2\left(|\tilde{\mathbf{r}}|-x\right)}}.
\]

Usually mapping~\eqref{eq:mapping} is written in coordinate form
\begin{equation}\label{eq:r2u}
\tilde{\mathbf r}\equiv
\begin{pmatrix}
    x \\ y \\ z \\ 0
\end{pmatrix} = 
\begin{pmatrix}
    u_1^2-u_2^2-u_3^2+u_4^2 \\
    2(u_1u_2-u_3u_4) \\
    2(u_1u_3+u_2u_4) \\
    0
\end{pmatrix}=L(\mathbf u)\mathbf u,
\end{equation}
where $\tilde{\mathbf r}$ is positional vector $\mathbf r$  augmented with forth zero coordinate. This mapping is defined by so-called $KS$-matrix:
\begin{equation*}\label{eq:Lmatrix}
    L(\mathbf u)=
    \begin{pmatrix}
    u_1 & -u_2 & -u_3 &  u_4\\
    u_2 &  u_1 & -u_4 & -u_3\\
    u_3 &  u_4 &  u_1 &  u_2\\
    u_4 & -u_3 &  u_2 & -u_1
    \end{pmatrix}.
\end{equation*}
For any vector $\mathbf u^\top=(u_1,u_2,u_3,u_4)$ matrix $L(\mathbf u)$ can be represented in the form
\[
L(\mathbf u)=\begin{pmatrix}
    \mathbf u^\top K_1 & \mathbf u^\top K_2 & \mathbf u^\top K_3 &
    \mathbf u^\top K_4
\end{pmatrix}^\top,
\]
where square matrices $K_j$, $j=1,2,3$, are square roots of identity matrix $\mathbb I_4$: $K_j^2=\mathbb I_4$,
\begin{equation}\label{eq:K123}
K_1=\begin{pmatrix}
    1 & 0 & 0 & 0 \\ 0 & -1 & 0 & 0 \\ 0 & 0 & -1 & 0 \\ 0 & 0 & 0 & 1\\
\end{pmatrix},\quad
K_2=\begin{pmatrix}
    0 & 1 & 0 & 0 \\ 1 & 0 & 0 & 0 \\ 0 & 0 & 0 &-1 \\ 0 & 0 & -1 & 0 
\end{pmatrix},\quad
K_3=\begin{pmatrix}
    0 & 0 & 1 & 0 \\ 0 & 0 & 0 & 1 \\ 1 & 0 & 0 & 0 \\ 0 & 1 & 0 & 1 \\
\end{pmatrix},
\end{equation}
and square matrix $K_4$ has another property: $K_4^2=-\mathbb I_4$ and
$K_4=\begin{pmatrix}
    0 & 0 & 0 & 1 \\ 0 & 0 & -1 & 0 \\ 0 & 1 & 0 & 0 \\ -1 & 0 & 0 & 0
\end{pmatrix}$.

As shown in~\cite{Poleshchikov2003} there exists a family of $KS$-matrices depending on six parameters, providing regularization of equations of motion and possessing the following properties:
\begin{gather}
L^\top(\mathbf u)L(\mathbf u) = L(\mathbf u)L^\top(\mathbf u)=\left|\mathbf u\right|^2E,\label{eq:L_orthogonal}\\
\left(L(\mathbf u)\mathbf v\right)_j=\left(L(\mathbf v)\mathbf u\right)_j,\quad j=1,2,3,\notag\\
\left(L(\mathbf u)\mathbf v\right)_4=-\left(L(\mathbf v)\mathbf u\right)_4.\notag
\end{gather}

Physical momenta transformation provided by 
\begin{equation}\label{eq:p2v}
\tilde{\mathbf p}=\frac{\boldsymbol\xi}{2\alpha|\mathbf r|},\text{ where }\boldsymbol{\xi}=L(\mathbf u)\mathbf v=L(\mathbf v)\mathbf u.
\end{equation}
Here $\alpha\in\mathbb R\backslash\{0\}$ is an arbitrary constant, and $\tilde{\mathbf p}$ vector of momenta $\mathbf p$ augmented with forth zero coordinate. If $\mathbf r=0$, then vector $\mathbf p$ is undefined, but vector $\boldsymbol{\xi}=0$.\ 
So, KS-coordinates $\mathbf u$ and $\mathbf v$ are connected to each other with additional particular integral $\mathbf{u}^\top K_4\mathbf{v}=0$.\
The inverse transformation according to~\eqref{eq:L_orthogonal} is the following
\begin{equation*}\label{eq:v2p}
\mathbf v=2\alpha L^\top(\mathbf u)\tilde{\mathbf p}.
\end{equation*}

Let the original Hamiltonian be represented as a perturbation of the Kepler's problem Hamiltonian 
\begin{equation}\label{eq:Kepler_perturbed}
H(t,\mathbf x,\mathbf p)=\frac12|\mathbf p|^2 - \frac1r+W(t,\mathbf x,\mathbf p).
\end{equation}
For a general non-autonomous Hamiltonian system with function $H(t,\mathbf x,\mathbf p)$ the regularization of the equations of motion consists of the following consecutive steps.
\begin{description}
\item[Step 1.] Transition to the extended autonomous system with Hamiltonian $\widetilde H=H(t,\mathbf x,\mathbf p) + h$, where $h$ is the momentum canonically conjugate to the variable $t$.\ The new momentum $h$ should be chosen in such a way that $\widetilde H\equiv0$.
\item[Step 2.] Introduction of a new independent variable (fictitious time) $\tau$:
\[
dt = \alpha  r d\tau,\quad \alpha\in\mathbb R\backslash\{0\},
\]
and transition to the Hamiltonian $H^*=  r \widetilde H$.
\item[Step 3.] Substitution of regular variables $\mathbf u$, $\mathbf v$ according to Formulas~\eqref{eq:r2u} and~\eqref{eq:p2v}.
\end{description}

As a result of these steps we obtain the regularized Hamiltonian in the form of
\begin{equation*}\label{eq:Ham3Dreg}
\mathcal H(\mathbf u,\mathbf v)=\frac1{8\alpha}|\mathbf v|^2+\alpha h |\mathbf u|^2+V(\tau,\mathbf u,\mathbf v)-\alpha,
\end{equation*}
where $V(\tau,\mathbf u,\mathbf v)=\alpha r W(t,\mathbf x,\mathbf p)$ after substitutions~\eqref{eq:r2u} and~\eqref{eq:p2v}.\ 
The constant coefficient $\alpha$ is usually chosen as $1$, but some authors prefer to put it equal to $1/4$ as applied in \cite{SIMO20001,gomez_marcote_mondelo}.\ In the last case the coefficient at $v^2$ term becomes $1/2$ as in the initial Hamiltonian~\eqref{eq:Kepler_perturbed} and the coefficient at $u^2$ term becomes $(h/2)/2$. Moreover, in the mentioned above papers was proposed a special scaling of the phase variables $\mathbf u,\mathbf v$ and independent variable $\tau$
\begin{equation}\label{eq:h_scaling}
\mathbf u=2^{3/4}|h|^{1/4}\mathbf Q,\quad \mathbf v=2^{1/4}|h|^{3/4}\mathbf P,\quad \mathcal H=\frac{1}{\sqrt 2|h|^{3/2}}\mathcal H,
\end{equation}
that reduces the regularized Hamiltonian into the form without parameter $h$. It is obvious that at the point $h=0$ the transformation~\eqref{eq:h_scaling} is not applicable.

\subsubsection{Lidov's method of computing critical solutions of family $\mathcal B_0$}\label{subsub:Lidov_method}

Here we provide a brief description of Lidov's method for computation of the so-called \textit{ critical solutions} (see Definition~\ref{def:critical_solution}) at which the spatial periodic orbits intersecting the $\mathcal B_0$ family. This method was described in detail in papers~\cite{lidov82,lidov82a}, but, unfortunately, the scientific journal ``Kosmicheskie issledovania'' (``Space Research'') in the USSR period was never translated into English. This description is just an attempt to familiarize English readers with it and discuss its some advantages and drawbacks. 

Before describing the method, here are the basic notations that will be needed later.\ 
$\mathbb I_n$ denotes the identity $n\times n$ matrix and $\mathbb J$ denotes the symplectic unit 
$\mathbb J=\begin{pmatrix}
    0_4 & \mathbb I_4 \\ -\mathbb I_4 & 0_4
\end{pmatrix}$. 
Phase vectors of original (physical) coordinates and KS-coordinates are denoted by $\mathbf q=(\mathbf r,\mathbf{p})$ and $\mathbf{w}=(\mathbf{u},\mathbf{v})$ correspondingly. 

In brief, Lidov's method consists of the following steps.
\begin{description}
\item[Step 1.] For the system of canonical equations of the Hill problem in KS-variables $\mathbf w=(\mathbf u,\mathbf v)$ 
\begin{equation}\label{eq:eqreg}
\frac{d\mathbf w}{d\tau} = \mathbb J\left(\frac{\partial\mathcal H}{\partial\mathbf w}\right)^\top,\quad \frac{dt}{d\tau}=\frac{\partial\mathcal{H}}{\partial h}=|\mathbf w
|^2,
\end{equation}
we write down a system of equations in variations. Its form differs from the analogous system in variations in physical coordinates, since the variation of regular variables occurs under the additional condition $\mathcal H(\mathbf u,\mathbf v,h)=0$:
\begin{equation}\label{eq:Vareqreg}
    \frac{dY}{d\tau}=\mathbb J \Hess\mathcal H(\mathbf u,\mathbf v)Y+\mathbb J\frac{\partial\grad\mathcal H}{\partial h}\circ\boldsymbol{\eta},\quad \frac{d\boldsymbol{\zeta}}{d\tau}=\frac{\partial}{\partial\mathbf w}\left(\frac{\partial\mathcal{H}}{\partial h}\right)Y,
\end{equation}
where $\circ$ is a dyadic product, vector $\boldsymbol{\zeta}=dt/d\mathbf w_0$ and vector $\boldsymbol{\eta}$ is the constant vector
\[
\boldsymbol{\eta}=-\left(\frac1{|\mathbf u|^2}\frac{\partial\mathcal H}{\partial\mathbf w}\right)_{\mathbf w=\mathbf w_0}.
\]
Fundamental matrix $Y$ of the system~\eqref{eq:Vareqreg} is defined by the initial conditions
\[
Y(0)=\mathbb I_8,\quad \boldsymbol{\zeta}=0.
\]

\item[Step 2.] For a given initial condition with a single nonzero coordinate $z_0=a$ at $t=0$, the simultaneous integration of the equations of motion~\eqref{eq:eqreg} and the system of variational equations~\eqref{eq:Vareqreg} along the periodic solution is performed.\ 
For symmetric solutions, it is sufficient to integrate these equations either on half a period (for solutions with one symmetry) or on a quarter of a period (for solutions with two symmetries).\ The monodromy matrix $Y$ and the auxiliary vector $\boldsymbol\zeta$ are then obtained at the $\tau=\tau^*$.

\item[Step 3.] Since the obtained values are used for calculations in the original coordinates, a special transformation of the fundamental matrix $Y$ in KS-coordinates to the monodromy matrix $M$ in the original coordinates is applied.\ The general scheme of such a transformation is given in the paper~\cite{lidov82a}, but here we will only give the formulas for the special case of the Hill problem and the KS regularization.\ Using proposed in~\cite{lidov82a} notation connection fundamental matrix $X(\tau^*)={\partial\mathbf q}/{\partial\mathbf{q}_0}$ of variational system in original coordinates with matrices $Y,\boldsymbol{\zeta}$ computed in KS-variables is the following:
\begin{equation}\label{eq:Y2X}
    X=\boldsymbol{\varepsilon}\begin{pmatrix}
        A_2 & \mathbf{0}_4 \\ A_1 & C_1
    \end{pmatrix},\quad 
    \boldsymbol{\varepsilon}=\frac{\partial\mathbf q}{\partial\mathbf{w}}Y-\mathbb J\left(\frac{\partial H}{\partial\mathbf q}\right)^\top\circ\boldsymbol{\zeta}.
\end{equation}
Matrix $\partial\mathbf q/\partial\mathbf{w}$ is the Jacobi matrix of KS transformation and its blocks can be obtained by differentiating formulas~\eqref{eq:r2u} and~\eqref{eq:p2v} by vectors $\mathbf{u}$ and $\mathbf{v}$ correspondingly.\ Together with matrices $Y,\boldsymbol{\zeta}$ all these matrices computed at the moment $\tau=\tau^*$. 

Structure of the matrices $A_1,A_2,C_1$ defined by KS-transformation.\ In the mentioned Lidov's paper all of them are expressed by Jacobi generating function $S(\mathbb p,\mathbf{u})=\mathbf{p}^\top L(\mathbf u)\mathbf u$:
\[
A_1=S_{\mathbf{u}\mathbf{u}}A_2,\quad A_2=S_{\mathbf{p}\mathbf{u}}B^{-1},\quad B=S_{\mathbf{u}\mathbf{p}}S_{\mathbf{p}\mathbf{u}},\quad C_1=S_{\mathbf{p}\mathbf{u}}.
\]
Their explicit expressions are the following:
\begin{equation}\label{eq:A1A2C1}
    A_1=\frac{1}{|\mathbf{u}|^2}\left(\sum_{j=1}^3 p_jK_j\right)L^\top(\mathbf{u}),\quad 
    A_2=\frac{1}{2|\mathbf{u}|^2}L^\top(\mathbf{u}),\quad
    C_1=2L^\top(\mathbf{u}),
\end{equation}
where matrices $K_j$ are taken from Formula~\eqref{eq:K123}. All these matrices are computed at the initial point $\mathbf{u}=\mathbf{u}_0$, $\mathbf{v}=\mathbf{v}_0$.

\item[Step 4.] At this step we have the monodromy matrix $M$ which is the fundamental matrix $Y$ computed on the whole period of solution.\ Two stability indices $S_1$ and $S_2$ are computed by formulas 
\begin{equation*}\label{eq:S12}
S_{1,2}=-\frac14\left(a\pm\sqrt D\right),\quad a=2-\Tr M,\quad D=2\Tr\left(M^2\right)-a^2+4.
\end{equation*}
Here we have 4 situations with stability of spatial periodic solution in non-degenerate case.\ If $D<0$ both indices are complex and mutually conjugate and so called \textit{complex instability} takes place. If $D>0$ then both indices are real.\ In the case $S_{1,2}\in(-1;+1)$ we have stable periodic solution, in the case $S_{1,2}\not\in(-1;+1)$ we have unstable, and, finally, if only one of the indices belongs to the open interval $(-1;+1)$ we have semi-stable situation. 

\item[Step 5.] Critical solution with initial condition $z_0=a^*$ and period $T^*$ takes place if one of the indices $S_j=\cos(2\pi p/q)$, $j=1,2$, where $p\in\mathbb N_0$, $q\in\mathbb N$.\ As was shown in~\cite{lidov82} the spatial periodic solutions in the vicinity of critical one can only be singly or doubly symmetric with respect to planes $XOZ$ and/or $YOZ$, i.e. they are invariant under transformations   $\rho_1$~\eqref{eq:rho1} or $\rho_3$~\eqref{eq:rho3}. Period $T$ of spatial solution is $q$-multiple to the period of the corresponding critical collision solution of the family $\mathcal{B}_0$.


\item[Step 6.] In the neighborhood of the generating solution of the family $\mathcal B_0$, one looks for corrections $\boldsymbol{\delta}=(\delta_1,\delta_2,\delta_3)$ to the initial conditions and $\delta T$ to the period, which ensure the periodicity condition with a certain accuracy.\ For initial conditions with $\rho_j$, $j=1,3$, symmetry these corrections are the following:
\[
\begin{aligned}
    \rho_1\colon & \boldsymbol{\delta}=(\delta x,\delta z,\delta p_y),\\
    \rho_3\colon & \boldsymbol{\delta}=(\delta y,\delta z,\delta p_x).
\end{aligned}
\]
\end{description} 

The method described above has the following drawback. If the initial condition $\mathbf{r}_0$ or the final point $\mathbf{u}(\tau^*)$ are close to the origin, i.e. to the singular point, then formulas~\eqref{eq:Y2X} and~\eqref{eq:A1A2C1} cannot be applicable.\ In this case all computations should be provided in regular coordinates only.

\subsection{Bifurcations of periodic  orbit families with double symmetry}

\begin{definition}\label{def:critical_solution}
Critical orbit is a periodic solution for which one of the stability indices $S_{1,2}$ is equal to the value $\cos(2\pi p/q)$, where $p\in\mathbb N_0$ and $q\in\mathbb N$. 
\end{definition}

Critical orbit can be considered from different points of view.\ According to H.\,Poincar\'e~\cite[Ch.~XXX]{Poincare1893} it is a bifurcation solution in the vicinity of which there exist one or more families of \textit{periodic solutions of second generation} ({\selectlanguage{french}`deuxi\`eme genre'}).\  M.\,H\'enon considered such orbit as a common orbit for two or more families of periodic solutions~\cite[Ch.~2]{Henon97}.\ Here we use both approaches, but prefer the second one.\ Since in this paper the main object of study is the network of families of symmetric spatial periodic solutions and their interaction with each other by means of critical orbits, we will point out here the peculiarities of bifurcations of families of doubly symmetric orbits.\ 
For families of planar singly symmetric solutions the possible types of bifurcations (interactions) described in~\cite[Ch.~VIII]{BrunoRTBPeng} based on the normal formal analysis, or in~\cite{Kreisman2005} based on the structure of the monodromy matrix.\ The case of families of doubly symmetric periodic solutions was considered in~\cite{BatkhinPCS2020}. In the generic case, the following proposition takes place.
\begin{proposition*}[{\cite{BatkhinPCS2020}}]
Let a doubly symmetric planar periodic solution with period $T$ have the stability index~$S$ equal to $\cos(2\pi p/q)$, where $p, q$ are coprime integrals.
\begin{itemize}
    \item If both numbers $p$ and $q$ are odd, then there exists only one family of doubly symmetric solutions with period $T' = qT$ in the vicinity of the original solution.
    \item If one of the numbers $p$ or $q$ is even, then there exist two pairs of families of singly symmetric solutions with different types of symmetry for each pair with period $T' = qT$ in the vicinity of the original solution.
\end{itemize}
\end{proposition*}

Out of plane branching is considered in~\cite{robin_numerical_1980} (see also~\cite{kalantonis} in more compact form), where only two types of symmetries: $\rho_1$~\eqref{eq:rho1} and $\rho_2$~\eqref{eq:rho2}.\ Here we have two scenarios of branching depending on the parity of denominator $q$.

In our research, we investigate the connection between families of spatial periodic orbits which branch mainly from doubly symmetric planar critical orbits (families $g$, $f$ or $\mathbf{f}_3$) or from vertical collision critical orbits of family $\mathcal{B}_0$.\ Branching from vertical collision critical orbits has its own peculiarities and would apparently be described.\ In this analysis, rigorous proofs for the assertions presented are not provided.

\subsubsection{Out-of-plane branching from planar doubly symmetric critical orbits}

Let $\gamma$ be a doubly symmetric planar critical orbit with vertical stability index $S_v=\cos(2\pi p/q)$.
\begin{enumerate}
    \item The denominator $q$ is odd, and $q\neq1$.\ There are two families of doubly symmetric spatial periodic orbits with complementary symmetries $\rho_j$ -- $\rho_k$ in the vicinity of $\gamma$, where $j=1,2$ and $k=3,4$.\ The complementarity of symmetries means that if, for example, the first family has symmetry $\rho_{1}$ -- $\rho_{3}$, then the second family has symmetry $\rho_{2}$ -- $\rho_{4}$.
    \item The denominator $q$ is even. There exist two pairs of families of doubly symmetric spatial periodic orbits with complementary symmetries $\rho_j$ -- $\rho_k$ in the vicinity of $\gamma$, where the pairs of indices $(j,k)$ should be $(1,2)$ ($(2,1)$) or $(3,4)$ ($(4,3)$).
\end{enumerate}

If denominator $q=1$ then there exists only one family of spatial doubly symmetric periodic orbits, as described in item 1.\ There are no limitations on any combination of symmetries that is evident from the description of critical orbits of the family $\mathbf{f}_3$ given in Subsection~\ref{sec:f_3}. Here the parity of the numerator does not affect the structure of symmetries of solutions in the neighborhood of critical solutions.

\subsubsection{Branching from vertical collision symmetric critical orbits}
The collision orbits of the family $\mathcal B_0$ lie at the intersection of two planes of symmetries: $XOZ$ and $YOZ$.\ In this case, the symmetries of spatial periodic solutions appearing in the vicinity of the critical orbits demonstrate completely different behavior.\ After providing the numerical experiment, it turns out that if the denominator $q$ in Definition~\ref{def:critical_solution} is odd, then there are two families of spatial periodic orbits with only one symmetry $\rho_1$ or $\rho_3$.\ If $q$ is even, then there are two families of spatial periodic orbits with two symmetries $\rho_1$ -- $\rho_3$ and $\rho_3$ -- $\rho_1$.\ Below is a conjecture that summarizes the obtained results.

\begin{conjecture*}[On symmetries of spatial solutions interacting with $\mathcal B_0$ family]\label{hyp:B0}
Let $\gamma$ be a critical solution of family $\mathcal{B}_0^+$ with the initial condition $x_0=y_0=0$, $z_0=a^*$, period $T^*$ and one of the stability indices $S_j=\cos(2\pi p/q)$, $j=1,2$.
\begin{enumerate}
\item If the denominator $q$ is odd, i.e. $q=2k+1$, $k\in\mathbb N$, then there exist two families of spatial periodic solutions with a single symmetry.\  Orbits of the first family have symmetry $\rho_1$ -- $\rho_1$, and orbits of the second one have symmetry $\rho_2$ -- $\rho_2$.\ For spatial orbits close to the critical solution $\gamma$ one point of the orthogonal intersection of the symmetry plane is located near $z\approx a^*$ and another point is located near the origin. 
\item If the denominator $q$ is even, i.e. $q=2k$, $k>1$, then there exist two families of doubly symmetric spatial periodic solutions.\ Orbits of the first family have symmetry $\rho_1$ -- $\rho_2$, and orbits of the second one have symmetry $\rho_2$ -- $\rho_1$.\ Location of the points of orthogonal intersection depends on the parity of $k$.\ If $k$ is odd ($q$ has the form $4n+2$) then these points are located as described above but in two different planes.\ If $k$ is even ($q$ has the form $4n$) then the orbits of one family have both points near $z\approx a^*$ and the orbits of the other family have both points near the origin.\ As mentioned at the end of Subsection~\ref{subsub:Lidov_method} such family could not be computed by the Lidov's method.
\end{enumerate}
\end{conjecture*}

If denominator $q=1$ then there exists only one family of spatial singly symmetric periodic orbits with $\rho_1$ or $\rho_3$ symmetry.\ If denominator $q=2$ then there exists only one family of spatial doubly symmetric periodic orbits, as described in item 2 above.\

\section{Interaction of Conley--Zehnder index with bifurcation points}
\label{sec3}

Let $\phi^t(\xi)$ be the Hamiltonian flow and let us consider the first return time map of the linearized flow of a $T$-periodic orbit through $\xi(0)$, given by the monodromy matrix $M = D \phi^T \left(\xi(0)\right)$.\ Since the Hamiltonian is time-independent, 1 appears trivially twice as an eigenvalue of $M$.\ Fixing the energy and modding out the direction of the Hamiltonian vector field yields the reduced monodromy~$M^{red}$ whose eigenvalues are called ``Floquet multipliers''.\ A periodic orbit is non-degenerate if 1 is not among its Floquet multipliers.\

The (transversal) Conley--Zehnder index $\mu_{CZ}$ assigns a mean winding number to non-degenerate periodic orbits, which is defined in terms of a path of symplectic matrices $\Psi(t)$ generated by the linearized flow along the whole periodic orbit (see \cite{hofer_w_z}, \cite{hofer_w_z_1} for the planar case and \cite{salamon_zehnder} for the general setting).\ This path is characterized by fixing a symplectic frame transverse to the direction of the Hamiltonian vector field in order to measure the rotation of the linearized flow.\ In particular, the symplectic path $\Psi(t)$ is a map from the interval $[0,T]$ to the group of symplectic matrices, where $\Psi(0) = \text{id}$ and $\Psi(T) = M^{red}$ is the reduced monodromy whose Floquet multipliers are different from 1 due to the non-degeneracy of the periodic orbit.\ The set of symplectic matrices with eigenvalue 1 is called ``Maslov cycle'', which corresponds to the space of reduced monodromies that are degenerate.\ The Conley--Zehnder index $\mu_{CZ}$ measures the twisting of the symplectic path $\Psi(t)$ by counting the number of crossing the Maslov cycle that lies between the starting point of a periodic orbit and its end point.\ An algorithmic use of the formal definition in numerical computations of periodic orbit families was performed in the recent work~\cite{moreno_aydin_koert_frauenfelder_koh}.\ In this paper, we determine the index by studying analytically the origin of the natural families of periodic orbits (see Section \ref{sec4}), then we continue these families, follow their Floquet multipliers and study the index jump at bifurcation points as discussed in the following way.\

We first consider planar periodic orbits that obey special geometric properties.\ Let $\xi^n$ be the $n$-times iteration of a planar orbit $\xi$.\ Assume that $\xi^n$ is non-degenerate for all $n \geqslant 1$.\ By virtue of \eqref{sigma} the reduced monodromy of planar orbits splits symplectically into planar and spatial blocks,
$$ M^{red} = \begin{pmatrix}
	M_p^{red} & 0\\
	0 & M_s
\end{pmatrix},\quad M_p^{red}, M_s \in Sp(2) = SL(2,\mathbb{R}). $$
Consequently, the Floquet multipliers are real (hyperbolic case) or lie on the unit circle (elliptic case), and the Conley--Zehnder index splits additively into a planar and spatial index, $ \mu_{CZ} = \mu_{CZ}^p + \mu_{CZ}^s$\vspace{0.5em}\\
\textbf{In the elliptic case.}\ The Floquet multipliers are of the form $e^{ \pm i n \theta_p}$ and $e^{\pm i n \theta_s}$, which means that each of the planar and spatial block of the linearized flow is conjugate to a rotation in $\mathbb{R}^2$,
\begin{align}\label{elliptic_case}
	M^{red} = \begin{pmatrix}
		M_p^{red} & 0\\
		0 & M_s
	\end{pmatrix} \sim \begin{pmatrix}
		\cos n \theta_p & - \sin n \theta_p & 0 & 0\\
		\sin n \theta_p & \cos n \theta_p & 0 & 0\\
		0 & 0 & \cos n \theta_s & - \sin n \theta_s \\
		0 & 0 & \sin n \theta_s & \cos n \theta_s
	\end{pmatrix}.
\end{align}
Notice that according to the formal definition of the index, each rotation angle is computed as a real number, and not modulo $2 \pi$.\ The index measures the number of times that Floquet multipliers crosses 1 during the first return time.\ In particular, each index is determined by
\begin{align}\label{index_rot}
    \mu_{CZ}^p (\xi^n) = 1 + 2 \cdot \lfloor n \cdot \theta_p /(2\pi) \rfloor, \quad \mu_{CZ}^s (\xi^n) = 1 + 2 \cdot \lfloor n \cdot \theta_s /(2\pi) \rfloor,
\end{align}
i.e., for every complete rotation the index jumps by 2 and is odd.\ The integer $ \text{rot}_k(\xi^n) := \lfloor n \cdot \theta_k /(2\pi) \rfloor$ is called planar or spatial ``rotation number'', correspondingly for $k = p,s$.\vspace{0.5em}\\
\textbf{In the positive/negative hyperbolic case.}\ Assume that for the simple closed orbit $\xi$ one pair of the Floquet multipliers is positive/negative real and of the form $\lambda$, $1 / \lambda$.\ Then the corresponding eigenvectors are rotated by $m\pi$ for an integer~$m$, and the corresponding index of $\xi^n$ equals
$$ \mu_{CZ}^{p/s}(\xi^n) = n m,\quad m \in \begin{cases}
	2 \mathbb{Z} & \text{ for the pos. hyperbolic case}\\
	2 \mathbb{Z} + 1 & \text{ for the neg. hyperbolic case.}
\end{cases} $$
\noindent
\textbf{Transition and bifurcation.}\ If the orbit changes its stability from elliptic to positive hyperbolic, or vice versa, then the eigenvalue 1 is crossed (corresponding to a degenerate orbit) and the index jumps by $\pm1$, according to direction of crossing the eigenvalue 1, see Figure \ref{scenario_1}.\
\begin{figure}[t]
	\centering
	\includegraphics[scale=1]{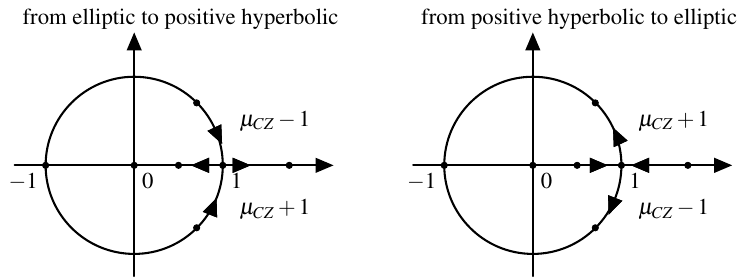}
	\caption{The index jump.}
	\label{scenario_1}
\end{figure}
If it stays elliptic, the index jumps by~$\pm 2$, e.g., if the rotation angle is a $k$-th root of unity, for $k \geq 3$ the $k$-th cover is still elliptic and goes through the eigenvalue 1.\ The crossing of the eigenvalue 1 generates bifurcation of new families \cite[Ch.~8]{abraham_marsden}; a planar-to-planar bifurcation appears when 1 is among its planar Floquet multipliers and a planar-to-spatial bifurcation appears when 1 is among its spatial Floquet multipliers.\ For the cases from elliptic to negative hyperbolic, or vice versa, the index does not change, but the index of its even cover jumps.\vspace{0.5em}\\
\textbf{Good and bad orbits.}\ In the language of SFT (Symplectic Field Theory) \cite{eliashberg_givental_hofer}, if 
\begin{align} \label{index_parity}
	\mu_{CZ}^p(\xi^n) \equiv \mu_{CZ}^p(\xi) \text{ mod } 2, \quad \mu_{CZ}^s(\xi^n) \equiv \mu_{CZ}^s(\xi) \text{ mod } 2,
\end{align}
or both equations in \eqref{index_parity} are not simultaneously satisfied, then $\xi^n$ is a ``good orbit''.\ Otherwise, $\xi^n$ is a ``bad orbit'' contributing nothing to the local Floer homology and its Euler characteristic.\ Therefore, if one checks the Euler characteristic at a bifurcation point and there appears bad orbits, they will not be counted.\ Note that all simple closed periodic orbits are good, and bad orbits occur as $n$-th cover of orbits whose exactly one pair of Floquet multipliers is of negative hyperbolic type, where~$n$ is even.\vspace{0.5em}\\
\textbf{Krein signature for symmetric periodic orbits and specification of index jump at bifurcation points.}\ The classical Krein signature \cite[Appendix~29]{arnold_avez} associates a $\pm$ sign to each elliptic Floquet multiplier.\ For symmetric periodic orbits, \cite{frauenfelder_moreno} constructed a $\pm$ sign to each elliptic or hyperbolic Floquet multiplier, which coincides with the classical Krein signature in the elliptic case.\ We now briefly recall the basic constructions in the general case and refer the curious reader to the latter reference for details.\ Let us consider a periodic orbit $\xi$ which satisfies $\xi(t) = \rho(\xi(-t))$, for all $t$.\ The tangent space at $\xi(0)$ splits into two Lagrangian submanifolds, $L_+ \oplus L_{-}$, where $L_{\pm}$ are the eigenspaces of the differential of $\rho$ at $\xi(0)$ to the eigenvalue~$\pm 1$.\ One can choose a symplectic basis with respect to this Lagrangian splitting such that the (reduced) monodromy matrix is of the special form
\begin{align}\label{monodromy_special_form}
    M = \begin{pmatrix}
	A & B\\
	C & A^T
\end{pmatrix},\quad B = B^T,\quad C = C^T,\quad AB = BA^T,\quad A^T C = CA, \quad A^2 - BC = \text{id}.
\end{align}
In particular, the spectrum of $M$ is determined by the spectrum of $A$:\ If $\lambda$ is an eigenvalue of $M$ then its stability index $\frac{1}{2}\left( \lambda + 1/\lambda \right)$ is an eigenvalue of $A$; if $a$ is an eigenvalue of $A$ then $\lambda(a) = a + \sqrt{a^2 - 1}$ is an eigenvalue of~$M$.\

Assume that $a$ is a real, simple and non-trivial eigenvalue of $A$, which means that $\lambda(a) \in S^1 \setminus \{ \pm 1 \}$ or $|\lambda(a)| \neq 1$, then the $B$-sign of the eigenvalue $\lambda(a)$ of $M$, defined by
$$ B\text{-sign} \left( \lambda(a) \right) := \text{sign} (v^T B v) = \pm, $$
is invariant under the choice of the symplectic basis to write down the matrix, where $v$ is an eigenvector of $A^T$ with eigenvalue $a$.\ The $C$-sign of $\lambda(a)$ is defined similarly by replacing $B$ with $C$, and $A^T$ by~$A$.\ For a pair of elliptic Floquet multipliers the $(C/B)$-signatures are different, while for a pair of hyperbolic Floquet multipliers the $(C/B)$-signatures are the same.\

In the 2 dimensional case, \eqref{monodromy_special_form} has the form $\begin{pmatrix}
	a & b\\
	c & a
\end{pmatrix}$, $a^2 - bc = 1$, which can be already found in the work by Darwin \cite{darwin}.\ Consequently for symmetric planar orbits, the signature of $b$ specifies the Floquet multipliers.\ We choose the orientation such that in the elliptic case~\eqref{elliptic_case}, if $b < 0$ (or $c > 0$) then the rotation is determined by $n \theta_p$, and if $b > 0$ (or $c < 0$) then the rotation is determined by $- n \theta_p$ (similar for the spatial rotation).\ This determines the direction of crossing the eigenvalue 1, and especially the index jump at bifurcation points.\

For spatial orbits we treat each pair of Floquet multipliers in the same way as described before.\ The case of complex instability can occur if two pairs of elliptic (or hyperbolic) Floquet multipliers collide at $S^1 \setminus \{ \pm 1 \}$ (or $\mathbb{R} \setminus \{ \pm 1 \}$) and move off it.\ In both such cases each pair consists of Floquet multipliers with different $B$-sign.\

\section{Conley--Zehnder indices of natural periodic orbits families}
\label{sec4}


\textbf{Families $g$, $f$ and $\mathcal{B}_0^{\pm}$.}\ For very low energies, the regularized Kepler problem is the source of the families $g$, $f$ and $\mathcal{B}_0^{\pm}$.\ In \cite[Chapter 8]{frauenfelder_koert} it is shown how the families $g$ and $f$ branch out from the regularized planar Kepler problem.\ In \cite{aydin_babylonian} this technique was extended to the spatial case and the Conley--Zehnder indices of the families $g$, $f$ and $\mathcal{B}_0^{\pm}$ were explored, which are given by
\begin{align}\label{indices}
		\mu_{CZ} = \begin{cases}
			6 = \mu_{CZ}^p + \mu_{CZ}^s = 3 + 3 & \text{for family $g$} \\
			4 & \text{for family $\mathcal{B}_0^{\pm}$} \\
            2 = \mu_{CZ}^p + \mu_{CZ}^s = 1 + 1 & \text{for family $f$}.
		\end{cases}
	\end{align}

\noindent
\textbf{Planar and vertical Lyapunov orbits.}\ It is well-known that the analysis of the linear behavior of the flow around collinear Lagrange points of the CR3BP is of the type $\text{saddle} \times \text{center} \times \text{center}$ \cite[Ch.~5]{szebehely}.\ Hence, by the Lyapunov center theorem \cite[Section~9.2]{meyer}, each collinear Lagrange point generates a pair of one-parameter families of periodic solutions which in the limit have planar and vertical frequencies, $\omega_p$ and $\omega_v$, related to both centers.\ These are the families of planar and vertical Lyapunov orbits.\ We now examine the relation of the frequencies and briefly recall their computations.\ We first write the Hamiltonian~\eqref{hamiltonian} as
$$ H(x,y,z,p_x,p_y,p_z) = \frac{1}{2} \left( (x-p_y)^2 + (y+p_x)^2 + p_z^2 \right) - U(x,y,z), $$
where $U(x,y,z) = \frac{1}{2} \left( x^2 + y^2 \right) + \frac{1 - \mu}{r_1} + \frac{\mu}{r_2}$.\ We expand $(x + \xi, y + \eta, z + \zeta)$ near $(x_{eq},0,0)$ by linearizing the equations~\eqref{ham_equation1}, where $x_{eq}$ denotes the x-coordinate of any of $L_{1/2/3}$, determined by~\eqref{x_coord_langrange}.\ This yields
\begin{align} \label{linearized_equations}
	\ddot{\xi} - 2 \dot{\eta} = U_{xx}^{eq} \xi,\quad \ddot{\eta} + 2 \dot{\xi} = U_{yy}^{eq} \eta,\quad \ddot{\zeta} = U_{zz}^{eq} \zeta,
\end{align}
where $U_{ii}^{eq}$ $(i=x,y,z)$ denotes the evaluation of the second order partial derivative of $U(x,y,z)$ at the collinear Lagrange point, given by
\begin{align}\label{equation_1}
	U_{xx}^{eq} = 2c + 1,\quad U_{yy}^{eq} = 1-c,\quad U_{zz}^{eq} = -c,\quad c := \frac{1 - \mu}{|x_{eq} + \mu|^3} + \frac{\mu}{|x_{eq} - 1 + \mu|^3}.
\end{align}
Since $c > 1$ for all $\mu \in (0,\frac{1}{2}]$, we have that at the three collinear Lagrange points $U_{xx}^{eq} > 0$, $U_{yy}^{eq} < 0$ and $U_{zz}^{eq} < 0$, for all $\mu \in (0,\frac{1}{2}]$.\ Evidently, the out-of-plane motion is uncoupled from the planar direction and it describes a harmonic oscillator with frequency $\omega_v = \sqrt{c}$ (vertical frequency).\ The characteristic polynomial associated to the two planar equations in~\eqref{linearized_equations} is
\begin{align}\label{equation_2}
	p(\lambda) = \lambda^4 + 2 \beta_1 \lambda^2 - \beta_2,\quad \beta_1 = 2 - \frac{1}{2} ( U_{xx}^{eq} + U_{yy}^{eq} ),\quad \beta_2 = - U_{xx}^{eq} U_{yy}^{eq} > 0.
\end{align}
In the variable $\alpha = \lambda^2$ the discriminant $\beta_1^2 + \beta_2 > 0$, therefore we have two non-zero real roots,
$$ \alpha_1 = - \beta_1 + \sqrt{\beta_1^2 + \beta_2},\quad \alpha_2 = - \beta_1 - \sqrt{\beta_1^2 + \beta_2}. $$
By Vieta's formula, $\alpha_1 \alpha_2 = - \beta_2 < 0$, which means that the two real roots have opposite sign.\ Hence, the eigenvalues are
$$ \lambda_{1/2} = \pm \sqrt{\alpha_1},\quad \lambda_{3/4} = \pm i \sqrt{-\alpha_2} =: \pm i \omega_p, \quad \lambda_{5/6} = \pm i \omega_v,  $$
where $\omega_p = \left( \beta_1 + \sqrt{\beta_1^2 + \beta_2} \right)^{\frac{1}{2}}$ indicates the planar frequency.\ The planar instability implies that the planar index $\mu_{CZ}^p$ equals 2 (see \cite[Ch.~8]{frauenfelder_koert}), and for the remaining indices we need the following lemma.\
\begin{lemma1}
	In the vicinity of the collinear Lagrange points the frequencies satisfy $\omega_v < \omega_p < 2 \omega_v$.\
\end{lemma1}
\begin{proof}
Recall that $c = \omega_v^2$.\ In view of \eqref{equation_1} and \eqref{equation_2} we have $\omega_p^2 = 1 - \frac{1}{2} c + \sqrt{ {9}c^2/{4} - 2 c }$.\
Via the equivalences
	$$ 1 < c \quad \Leftrightarrow \quad \left( \frac{3}{2}c - 1 \right)^2 < \frac{9}{4}c^2 - 2c \quad \Leftrightarrow \quad c < 1 - \frac{1}{2}c + \sqrt{ \frac{9}{4}c^2 - 2c } $$
we obtain $\omega_v^2 < \omega_p^2$ which shows the first inequality.\ In order to show the second inequality we consider the polynomial function $f(x) = 18x^2 - 7x + 1$.\ Since the discriminant is negative, it has no real roots, and because $f(0)=1$, $f(x) > 0$ for all $x \in \mathbb{R}$.\ Then the following equivalences
	$$ 0 < f(c) \quad \Leftrightarrow \quad \frac{9}{4}c^2 - 2c < \left( \frac{9}{2}c - 1 \right)^2 \quad \Leftrightarrow \quad 1 - \frac{1}{2} c + \sqrt{ \frac{9}{4}c^2 - 2 c } < 4 c $$
    yield $\omega_p^2 < 4 \omega_v^2$ which proves the second inequality.\
\end{proof}
\begin{corollary1}
	In the vicinity of the collinear Lagrange points it holds that
	\begin{align*}
		\mu_{CZ} = \begin{cases}
			3 = \mu_{CZ}^p + \mu_{CZ}^s = 2 + 1 & \text{for the family of planar Lyapunov orbits} \\
			5 & \text{for the family of vertical Lyapunov orbits.}
		\end{cases}
	\end{align*}
\end{corollary1}
\begin{remark}
	 In the Hill 3BP one readily obtains the same result by computing the eigenvalues,
	$$ \lambda_{1/2} = \pm \sqrt{2\sqrt{7} + 1},\quad \lambda_{3/4} = \pm i \sqrt{2\sqrt{7}-1} =: \pm i \omega_p,\quad \lambda_{5/6} = \pm 2 i =: \pm i \omega_v. $$
\end{remark}

\section{Results}
\label{sec5}

In this section, we construct bifurcation graphs as proposed in \cite{aydin_cz}, namely as labeled graphs, whose vertices are bifurcation points and whose edges correspond to families of periodic orbits, labeled with Conley--Zehnder indices.\ For its construction, we draw from bottom to top in the direction of decreasing the Jacobi constant $\Gamma = -2H$.\ Thick families indicate doubly symmetric orbits and thin families indicate simple symmetric orbits.\ Edges at critical orbits associated to planar periodic orbit families, we draw in black, vertically and shortly before and after the bifurcations.\ Edges of spatial periodic orbit families we draw colored.\ Symmetric families have same color, if one is dashed then the corresponding symmetry is the symplectic one $\sigma$ (reflection at the ecliptic), and if it is not dashed then the corresponding symmetry is one of the anti-symplectic symmetries.\ Families labeled with overlined indices indicate bad orbits.\ By ``b-d'' we denote a periodic of birth-death type.\ In general, a periodic orbit of ``birth-death'' type is a degenerate orbit from which two families bifurcate with an index difference of 1 and into the same energy direction.\ Its local Floer homology and its Euler characteristic are therefore zero.\

Plots of orbits have the same color as the corresponding bifurcation graph, and are labeled with~$\Gamma$.\ We plot spatial orbits with projections onto the three coordinate planes.\ From light to dark of orbits in one plot indicates the decreasing of the Jacobi constant.\ In case of more plots of the same family in a row, the beginning is on the right and its continuation towards the left, from row to row.\ In case of more plots of the same family in a column, the beginning is on the top and continues downwards.\

In stability diagrams we denote relevant underlying critical periodic orbits at bifurcation points in the corresponding family with upper index $(n,k)$, where $n$ denotes the $n$-times iteration where bifurcation takes place, and $k$ denotes the order associated to $n$-th cover bifurcation.\ The covering number~$n$ is also labeled on the right $y$-axis.\

Furthermore, considered families are closed by default and in the case when they are open we point it out explicitly.

\subsection{Families $g$, $g'$ and $f$}

\begin{figure}[t!]
	\centering
	\includegraphics[width=1\linewidth]{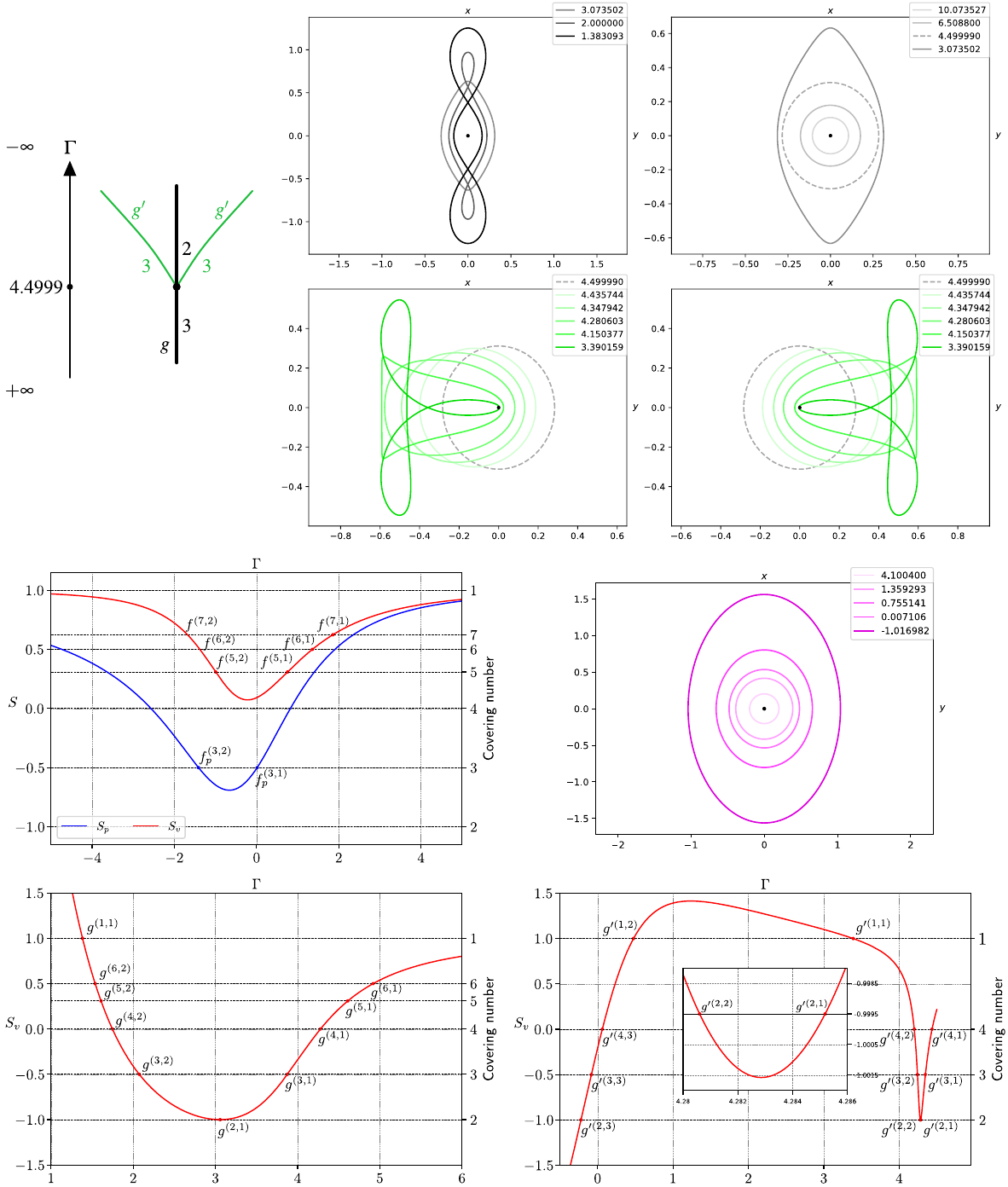}
	\caption{Left top:\ Bifurcation graph related to $g$ and $g'$.\ Plots on the top shows $g$-orbits.\ The gray dashed $g$-orbit at $\Gamma = 4.49999$ is where symmetry-breaking bifurcation of $g'$ happens; some $g'$-orbits are plotted in the second row in green with symmetric ones.\ Middle right shows $f$-orbits in purple and middle left their planar and vertical stability indices.\ Diagrams at the bottom show vertical stability indices associated to families $g$ (left) and $g'$ (right).}
	\label{bifurcation_graph_0}
\end{figure}

The Conley--Zehnder indices of families $g$, $g'$ and $f$ are discussed in detail in \cite{aydin_cz}.\ Their data base are collected in Table \ref{data_g_g'_f}, some orbits and stability diagrams are shown in Figure \ref{bifurcation_graph_0}.\ The orbits of the families~$g$ and $f$ start with infinitesimal circular direct and retrograde periodic orbits around the origin, which are doubly symmetric with respect to the axis $OX$ and $OY$.\ The $g$- and $f$-orbits start being planar and spatial elliptic, thus by \eqref{indices} and \eqref{index_rot} we have
\begin{align*}
		\mu_{CZ}^p = \mu_{CZ}^s = 1 + 2 \cdot \text{rot}_{p/s} = \begin{cases}
			  3 = 1 + 2 \cdot 1 & \text{for family }g \text{ (rotation angles have negative $B$-sign)} \\
            1 = 1 + 2 \cdot 0 & \text{for family }f \text{ (rotation angles have positive $B$-sign)}. \\
		\end{cases}
	\end{align*}

This property of family $f$ holds for all Jacobi constants.\ In view of the corresponding stability diagram in Figure~\ref{bifurcation_graph_0} and since their $B$-signs are positive, the planar and spatial rotation angles decrease, reach a minimum, and then they increase.\ The minimum rotation angle generating multiple cover planar bifurcations is 3 (at $f_p^{(3,1)}$ and $f_p^{(3,2)}$).\ The minimum one where multiple cover out-of-plane bifurcations occur is 5 (at $f^{(5,1)}$ and $f^{(5,2)}$).\ In the limit for $\Gamma \to - \infty$ the $f$-orbits converge to a retrograde ellipse around the origin with period $2 \pi$, which is degenerate and whose semi-minor axis is double as its semi-major axis (both rotation angles converge to $2 \pi$ as well).\

The properties for family $g$ are summarized as follows (signs in brackets correspond to $B$-signs):\
\begin{center}
    \begin{tabular}{cccc}
         $\Gamma$ & planar & spatial & $\mu_{CZ} = \mu_{CZ}^p + \mu_{CZ}^s$ \\
         $(+ \infty , 4.499990)$ & elliptic ($-$) & elliptic ($-$) & $6 = 3 + 3$ \\
         $(4.499990 , 3.057470)$ & pos.\ hyperbolic ($-$) & elliptic ($-$) & $5 = 2 + 3$ \\
         $(3.057470 , 1.383093)$ & pos.\ hyperbolic ($-$) & elliptic ($+$) & $5 = 2 + 3$ \\
         $(1.383093 , - \infty)$ & pos.\ hyperbolic ($-$) & pos.\ hyperbolic ($-$) & $6 = 2 + 4$
    \end{tabular}
\end{center}

Due to spatial $B$-signs and stability diagram $g$ in Figure~\ref{bifurcation_graph_0}, we see that the spatial rotation angles increase without becoming neg.\ hyperbolic.\ At the point of tangency, $g^{(2,1)}$ at $\Gamma = 3.057470$ the spatial rotation angle equals $\pi$, where the $B$-sign changes.\ At $g^{(1,1)}$ ($\Gamma = 1.383093$) the eigenvalue 1 is crossed from below (positive $B$-sign), thus the spatial index jumps from 3 to 4.\ In between in the Jacobi constant range $(+ \infty , 1.383093)$ are critical orbits associated to multiple cover out-of-plane bifurcations that are marked in the stability diagram for family $g$ in Figure~\ref{bifurcation_graph_0}.\

At $\Gamma = 4.49999$ the planar rotation angle crosses the eigenvalue 1 from above (negative $B$-sign), hence the planar index jumps from 3 to 2.\ This bifurcation point corresponds to the symmetry-breaking pitchfork bifurcation of $g'$, as shown in Figure~\ref{bifurcation_graph_0}.\ The $g'$-orbits begin with egg-shaped orbits and are simple symmetric with respect to the $OX$-axis, whose symmetric orbits are obtained by the $OY$-axis symmetry.\

The $g'$-orbits start being elliptic with planar index 3.\ Let us verify that this is in accordance with the Euler characteristics before and after bifurcation of $g'$, as shown in Figure \ref{bifurcation_graph_0}.\ By counting periodic orbits with indices, before and after bifurcation the Euler characteristics are $(-1)^3 = -1$ and $(-1)^2 + 2 \cdot (-1)^3 = -1$.\ Relevant critical orbits associated to $g'$ are highlighted in the corresponding stability diagram in Figure~\ref{bifurcation_graph_0}.\ Their indices and properties can be found in their data base in Table \ref{data_g_g'_f}.\

\subsection{Planar and vertical Lyapunov, halo orbits and $\mathcal{B}_0^{\pm}$}
\label{sec:5.2}

\subsubsection{Planar and vertical Lyapunov orbits}

\begin{figure}[t!]
	\centering
	\includegraphics[width=1\linewidth]{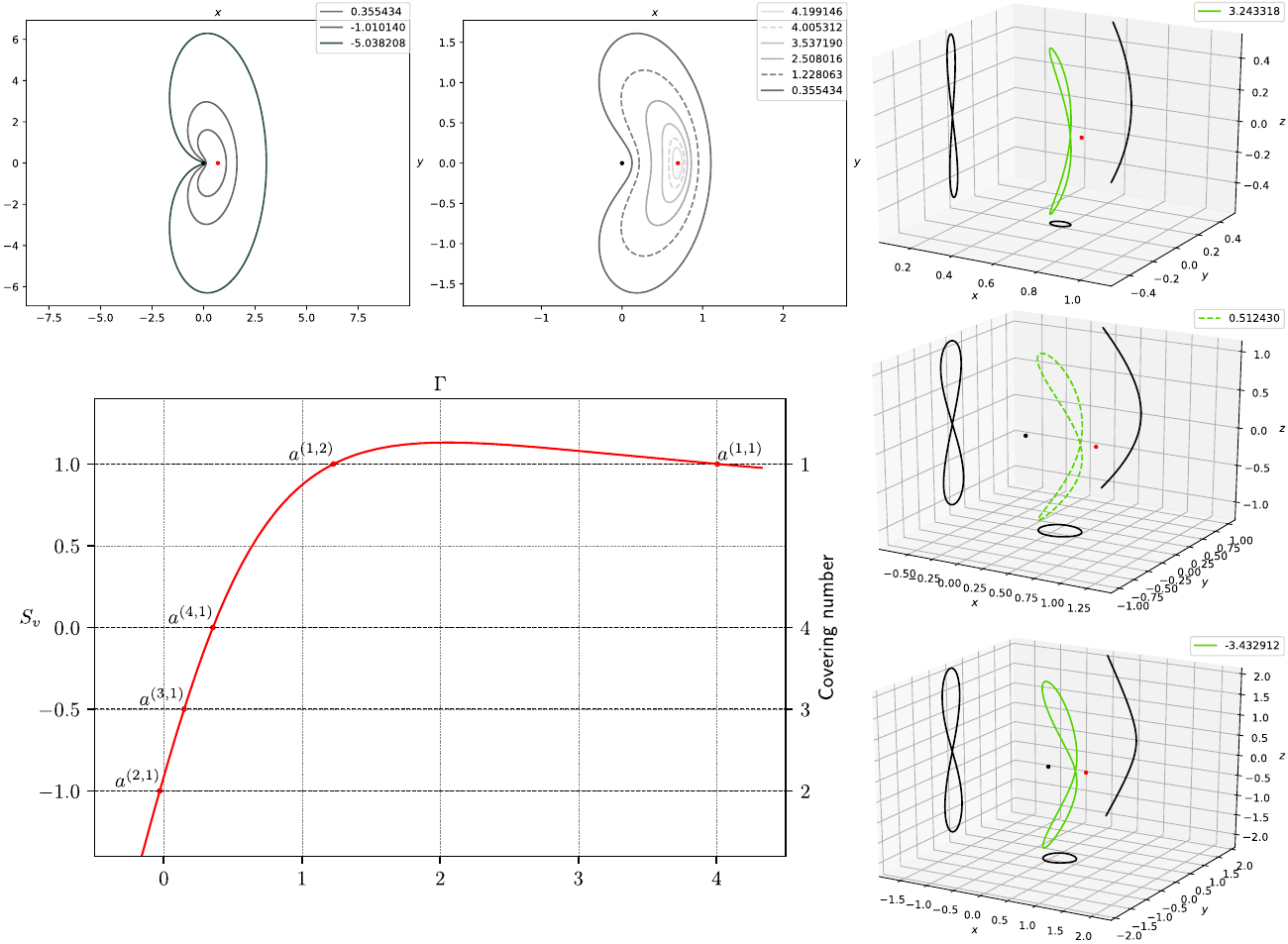}
	\caption{Top middle and left:\ Plots of planar Lyapunov orbits in gray (red dot corresponds to $L_2$); below is the corresponding vertical stability diagram.\ First gray dashed orbit at $\Gamma = 4.005312$ corresponds to $a^{(1,1)}$ and second gray dashed orbit at $\Gamma = 1.228063$ corresponds to $a^{(1,2)}$.\ Right:\ Plots of vertical Lyapunov orbits in green, from top to bottom; second dashed orbit is degenerate.}
	\label{figure_p_v_lyapunov}
\end{figure}

Since $L_1$ and $L_2$ are symmetric to each other, we restrict our study to $L_2$ located on the positive $x$-axis.\ The family $a$ of planar Lyapunov orbits start with infinitesimal retrograde elliptic orbits around~$L_2$, and are simple symmetric with respect to the $OX$-axis.\ Some orbits and the vertical stability diagram are shown in Figure~\ref{figure_p_v_lyapunov}.\ Their data are collected in Table \ref{data_a}, from which we deduce their properties:\

\begin{center}
    \begin{tabular}{cccc}
         $\Gamma$ & planar & spatial & $\mu_{CZ} = \mu_{CZ}^p + \mu_{CZ}^s$ \\
         $(4.326749 , 4.005312)$ & pos.\ hyperbolic ($-$) & elliptic ($+$) & $3 = 2 + 1$ \\
         $(4.005312 , 1.228063)$ & pos.\ hyperbolic ($-$) & pos.\ hyperbolic ($+$) & $4 = 2 + 2$ \\
         $(1.228063 , -0.029389)$ & pos.\ hyperbolic ($-$) & elliptic ($-$) & $5 = 2 + 3$ \\
         $(-0.029389 , - \infty)$ & pos.\ hyperbolic ($-$) & neg.\ hyperbolic ($+$) & $5 = 2 + 3$
    \end{tabular}
\end{center}

While the planar behavior is positive hyperbolic with negative $B$-sign and constant planar index 2 for all Jacobi constants, the vertical stability property allows us to generate out-of-plane bifurcations, as shown in the stability diagram in Figure \ref{figure_p_v_lyapunov}.\ The vertical stability starts being elliptic with positive $B$-sign in the Jacobi constant range $(4.326749 , 4.005312)$, then at the critical orbit $a^{(1,1)}$ the eigenvalue 1 is crossed from below, from which on the vertical stability is positive hyperbolic with positive $B$-sign and the spatial index jumps from 1 to 2.\ Then the spatial Floquet multipliers at $\Gamma = 1.228063$ become elliptic again with negative $B$-sign, which means that the spatial index jumps from 2 to 3 at $a^{(1,2)}$.\ At the Jacobi constant $\Gamma = - 0.029389$ the spatial Floquet multipliers become negative hyperbolic with positive $B$-sign, and stays so for all Jacobi constants after $\Gamma = - 0.029389$.\ In the further discussions in this paper we consider bifurcations from $a^{(1,1)}$ (gives rise to halo orbits, see Section \ref{sec:halo}) and $a^{(1,2)}$ (see Section~\ref{sec:pl_vert_lyapunov}).\ We note that we have also studied bifurcation from the critical points $a^{(j,1)}$, $j=2,3,4$, but since we could not find any connections to other families, we forgo their description in this paper.\

The vertical Lyapunov orbits begin as an infinitesimal oscillation in the $z$ direction and are of ``figure-eight'' shape, with the eight clamped at $L_2$, one lobe is above and one lobe is below the ecliptic $\{ z = 0 \}$.\ These orbits are doubly symmetric with respect to the $OX$-axis and to the $XOZ$-plane, hence they are invariant under $\sigma$ (reflection at the ecliptic).\ Some orbits are plotted in Figure \ref{figure_p_v_lyapunov} and their data can be found in Table \ref{data_vertical_lyapunov}.\ One pair of Floquet multipliers stays positive hyperbolic with negative $B$-sign for all the time.\ The other pair of Floquet multipliers starts being elliptic with negative $B$-sign and becomes at $\Gamma = 0.512430$ positive hyperbolic with negative $B$-sign, therefore the eigenvalue 1 is crossed from above, which means that the index jumps from 5 to 4.\ For their further continuation these orbits stay positive hyperbolic and positive hyperbolic, both with negative $B$-sign.\ This is an example of an open family.

\subsubsection{Bridge between planar and vertical Lyapunov orbits}
\label{sec:pl_vert_lyapunov}

The second critical point $a^{(1,2)}$ of family $a$ gives rise to a branch of spatial orbits that terminates at the critical orbit corresponding to vertical Lypaunov orbits at $\Gamma = 0.512430$.\ The data for the orbits of this branch are collected in Table \ref{data_a_2v} and some orbits are plotted in Figure~\ref{figure_second_vertical}.\ The orbits of this branch are positive hyperbolic and  positive hyperbolic, both with negative $B$-sign, for all $\Gamma$-values in the range $(1.228063 , 0.512430)$ with constant $\mu_{CZ} = 4$.\ Therefore, these orbits form a bridge between planar and vertical Lyapunov orbits.\
		
\begin{figure}[t!]
	\centering
	\includegraphics[width=1\linewidth]{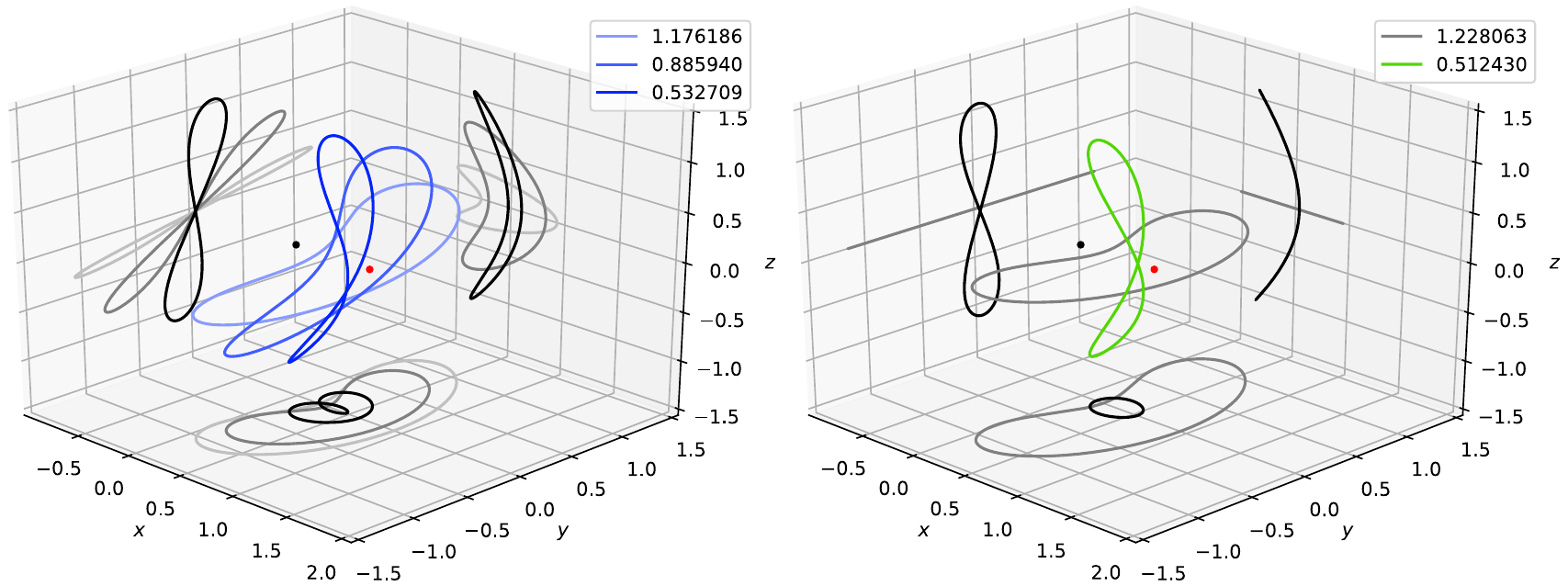}
	\caption{Bridge between planar and vertical Lyapunov orbits.\ Right:\ Critical planar Lyapunov orbit $a^{(1,2)}$ (gray) at $\Gamma = 1.228063$ and critical vertical Lyapunov orbit (green) at $\Gamma = 0.512430$.\ Left:\ Some orbits in blue associated to the bridge.}
	\label{figure_second_vertical}
\end{figure}
		
\subsubsection{Family $\mathcal B_0^{\pm}$}

Here we describe the family $\mathcal B_0$ of vertical collision orbits which in regular coordinates become the solutions to nonlinear oscillator.\ Let initial condition for any periodic solution of the family $\mathcal B_0$ be defined by the vector $\mathbf{q}_0=(\mathbf{r}_0,\mathbf{p}_0)^\top$, where $\mathbf r_0=(0,0,a)^\top$ and $\mathbf{p}_0=\mathbf{0}_3^\top$.\ Assuming in~\eqref{eq:reg2phy} $\varphi=0$ one gets that in KS-variables initial condition $\mathbf{w}_0$ as follows
\[
\mathbf u_0=\left(\sqrt{a/2},0,\sqrt{a/2},0\right)^\top,\quad  \mathbf v_0=\mathbf{0}_4^\top. 
\]

We analyze the stability of the family $\mathcal{B}_0$ in the linear approximation using KS regularization as described in Subsection~\ref{sec:KS}.\ Figure \ref{fig:B0stab} shows the behavior of the stability indices of this family according to the initial condition, the initial $z$-coordinate $a$.\ Only the part of the graph corresponding to the case where both the indices $S_1$ and $S_2$ are real is shown here.\ When $a>1.31808046$, the family enters the region of complex instability when two pairs of Floquet multipliers come off the unit circle, with one pair shifting inside the circle and the other shifting outside.

\begin{figure}[t]
\centering
\includegraphics[width=.7\textwidth]{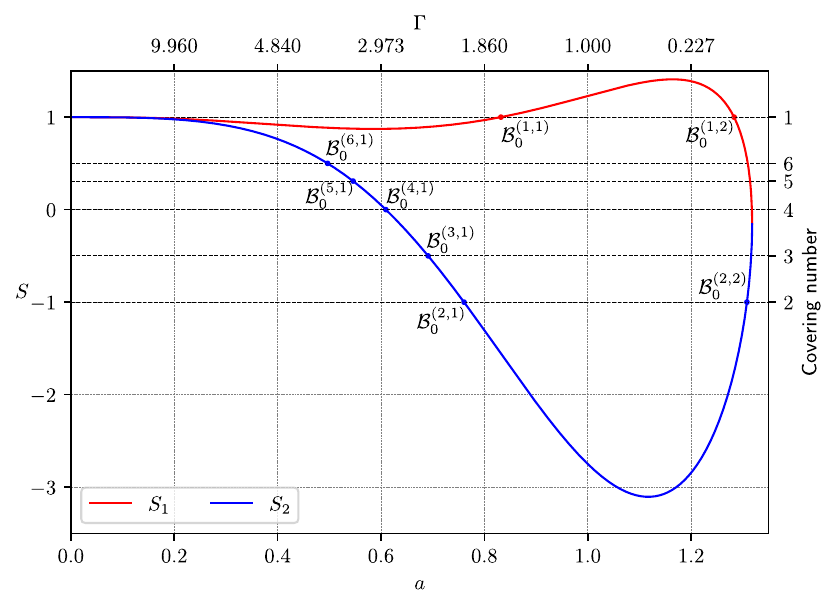}
\caption{The behavior of the stability indices of the $\mathcal{B}_0$ family.}\label{fig:B0stab}
\end{figure}

Two stability regions of the $\mathcal{B}_0$-family can be distinguished.\ The first region corresponds to values of the initial parameter $a$ from the interval $[0;0.83211]$, the second region corresponds to values from the interval $[1.28312;1.31808046]$.\ For small values of $a$ ($a<0.2$, $\Gamma\geqslant10$), the stability indices are almost equal to $1$, indicating that the Hill problem is close to an integrable one for large values of $\Gamma$.\ At $a\leqslant0.76070$, both indices take values from the interval $[-1;+1]$, with the index $S_1$ decreasing up to the value $0.8722$, then increasing, passing the critical value $+1$.\ The index $S_2$ decreases monotonically, taking successively all values from the stability interval.\ The first index $S_1$ reaches a maximum at $a\approx1.16$, then decreases passing the critical value $+1$, and the second index $S_2$ reaches a minimum at $a\approx1.12$ and increases passing the critical value $-1$.\ 

In Figure~\ref{fig:B0stab}, the horizontal lines indicate the stability indices corresponding to the critical values (from top to bottom) $+1$; $+0.5$; $+0.309\approx\cos(2\pi/5)$; $0$; $-0.5$; $-1$.\ At corresponding values of the parameter $a$, families of spatial periodic solutions providing $q$-covering, where $q\in\{1, 6, 5, 4, 3, 2\}$, appear in the neighborhood of the orbit $\mathcal{B}_0$, respectively.\ The numbers in Figure~\ref{fig:B0stab} mark those values of $a$ for which we were able to find and study the generated periodic solutions.\ Some of these solutions are partially described in~\cite{lidov_lyakhova83}.
We summarize the properties of $\mathcal{B}_0$ orbits as follows, based on data from Table \ref{data_B0}:\
\begin{center}
    \begin{tabular}{cccc}
         $\Gamma$ & $a$ & stability behavior & $\mu_{CZ}$ \\
         $(+ \infty , 2.051489)$ & $(0, 0.7607)$ & elliptic ($-$) \& elliptic ($+$) & $4$ \\
         $(2.051489 , 1.709346)$ & $(0.7607 , 0.8321)$ & elliptic ($-$) \& neg. hyperbolic ($-$) & $4$ \\
         $(1.709346, -0.087809)$ & $(0.8321 , 1.2831)$ & pos. hyperbolic ($+$) \& neg. hyperbolic ($-$) & $3$ \\
         $(-0.087809 , -0.182001)$ & $(1.2831 , 1.3080)$ & elliptic ($+$) \& neg. hyperbolic ($-$) & $2$ \\
         $(-0.182001 , -0.219862)$ & $(1.3080 , 1.3180)$ & elliptic ($+$) \& elliptic ($-$) & $2$
    \end{tabular}
\end{center}

We did not provide a precise computation of $\mathcal B_0$ stability indices on the long interval of $a$.\ So we cannot state that there are no other intervals of stability.

\subsubsection{Halo orbits:\ From 1st out-of-plane bifurcation from planar Lyapunov to $\mathcal B_0^{\pm}$-orbits}
\label{sec:halo}

\begin{figure}[t!]
	\centering
	\includegraphics[width=1\linewidth]{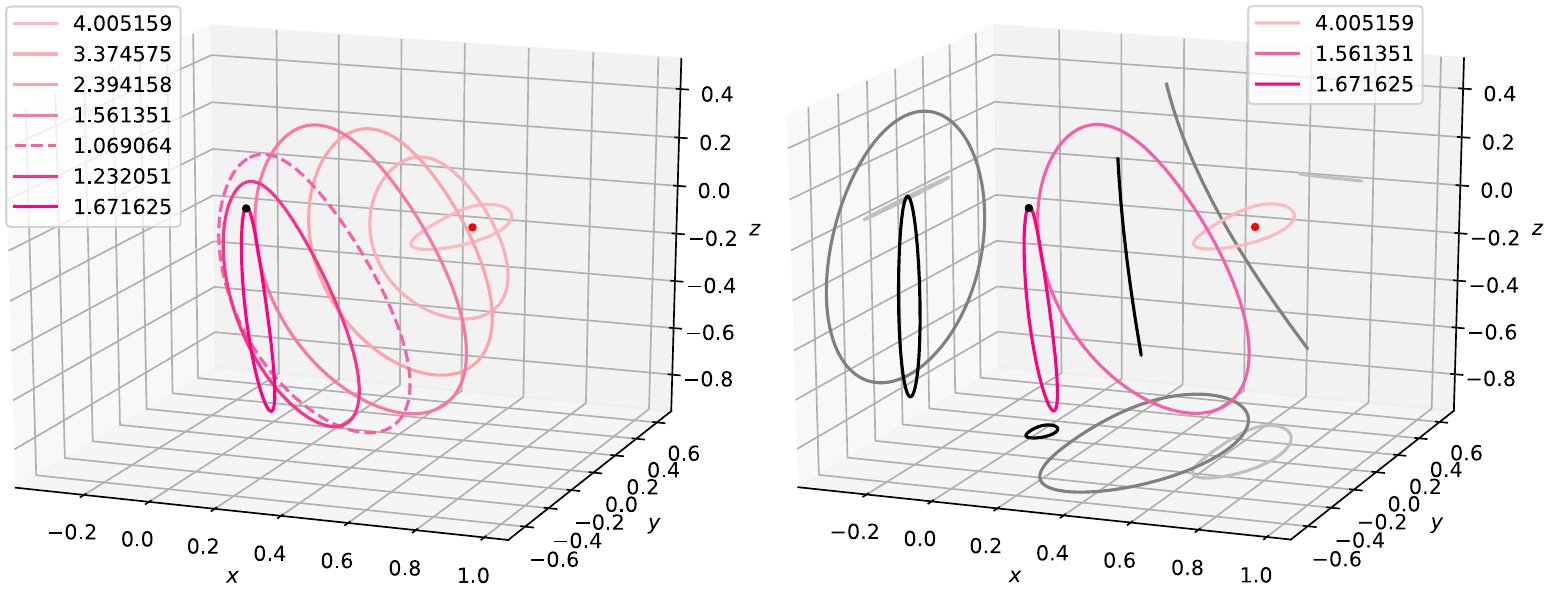}
	\caption{Plot of halo orbits in purple, starting from $a^{(1,1)}$ ($\Gamma = 4.005312$) and after a birth-death (dashed orbit on the left at $\Gamma = 1.069064$), they terminate at $B_0^-$ ($\Gamma = 1.709346$).}
	\label{figure_halo}
\end{figure}

\begin{figure}[t!]
	\centering
	\includegraphics[width=.7\textwidth]{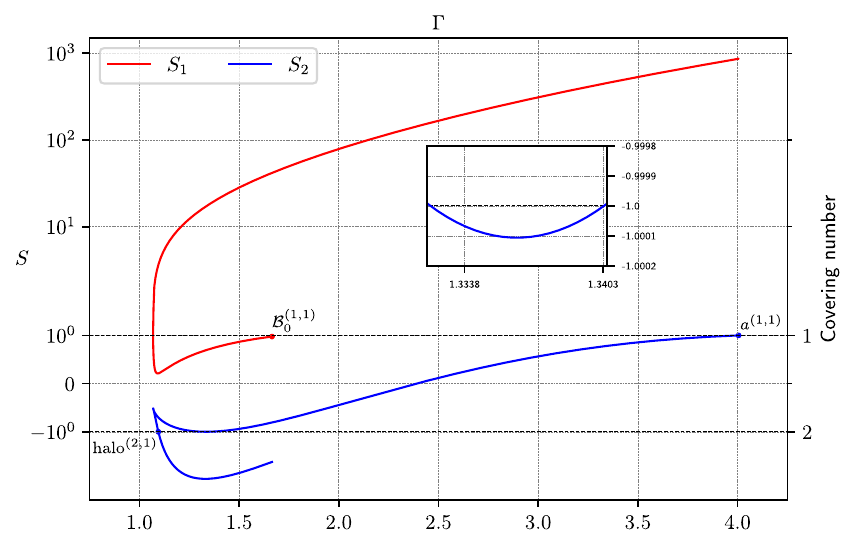}
	\caption{Stability diagram associated to halo orbits; $\text{halo}^{(2,1)}$ at $\Gamma = 1.095146$ plays a significant role in our investigations (see Section \ref{sec:g'_halo}).}
	\label{figure_halo_stab}
\end{figure}

Halo orbits are spatial orbits branching out from $a^{(1,1)}$ at $\Gamma = 4.005312$.\ These spatial orbits are simple symmetric with respect to the $XOZ$-plane, and by using $\sigma$ (reflection at the ecliptic) one obtains the symmetric family.\ Some orbits are plotted in Figure \ref{figure_halo}, the data are collected in Table \ref{data_halo} and the stability diagram is shown in Figure \ref{figure_halo_stab}.\ After a birth-death bifurcation at around $\Gamma = 1.069064$, they terminate at $\mathcal B_0^{(1,1)}$ ($\Gamma = 1.709346$).\ Therefore, halo orbits can be generated from $a^{(1,1)}$ or from $\mathcal B_0^{(1,1)}$.\ In addition to the stability diagram in Figure \ref{figure_halo_stab} our symplectic approach gives the following results:\

\begin{center}
    \begin{tabular}{ccc}
         $\Gamma$ & stability behavior & $\mu_{CZ}$ \\
         $(4.005312 , 1.339486)$ & pos.\ hyperbolic ($-$) \&  elliptic ($+$) & $3$ \\
         $(1.339486,1.317297)$ & pos.\ hyperbolic ($-$) \&  neg. hyperbolic ($+$) & $3$ \\
         $(1.317297, 1.069064)$ & pos.\ hyperbolic ($-$) \&  elliptic ($-$) & $3$ \\
         $(1.069064, 1.095146)$ & elliptic ($-$) \&  elliptic ($-$) & $4$ \\
         $(1.095146, 1.709346)$ & elliptic ($-$) \&  neg. hyperbolic ($-$) & $4$
    \end{tabular}
\end{center}

\subsubsection{Bifurcation graph}

The core results about connections at bifurcation points related to the previous discussions in this Section~\ref{sec:5.2} are illustrated in the bifurcation graph in Figure \ref{bifurcation_graph_1}.\ In particular, this bifurcation result demonstrates how essential it is to work in the spatial system in regularized coordinates.\ If one considers only the planar problem, the planar Lyapunov orbits are isolated and have no connection to other families.\
		
\begin{figure}[t!]
	\centering
	\includegraphics[width=1\linewidth]{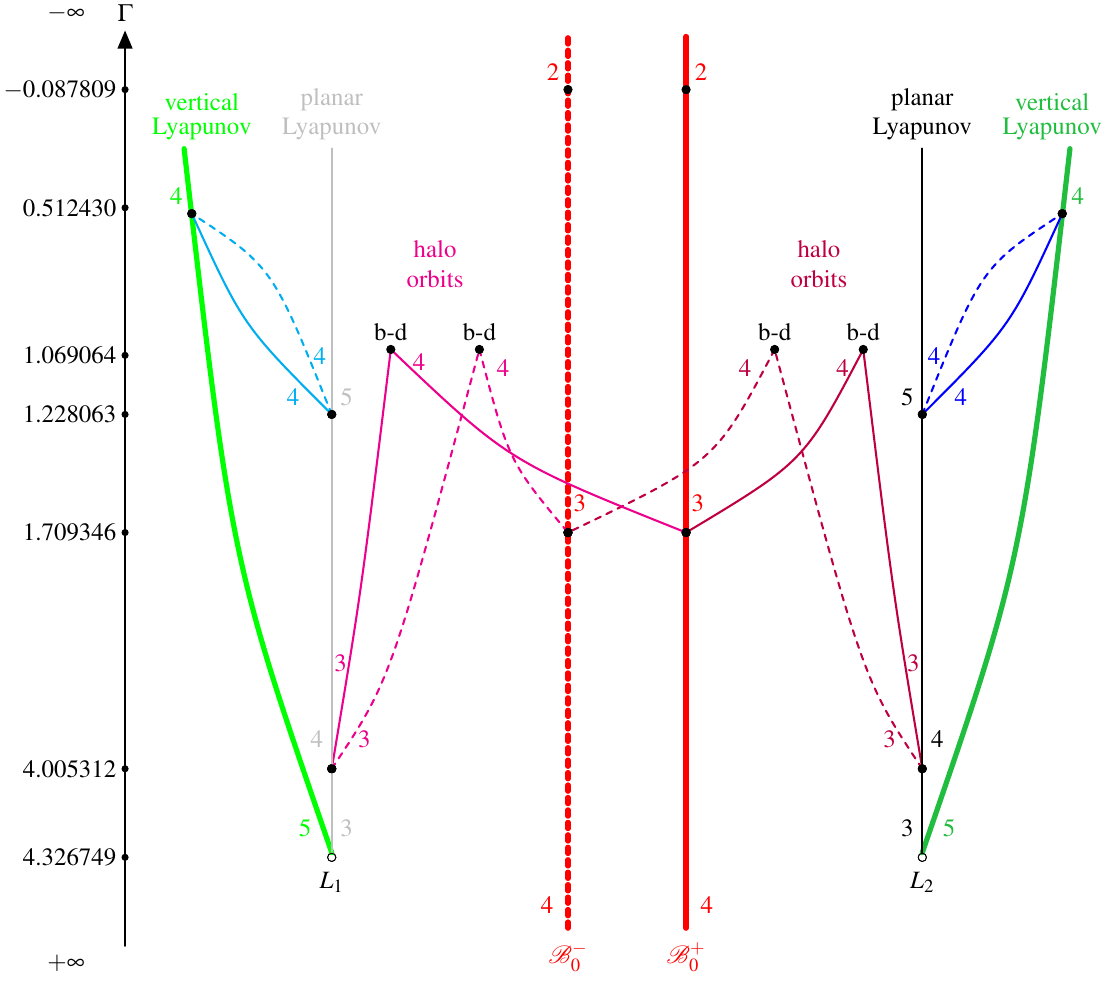}
	\caption{Bifurcation graph between families of planar and vertical Lyapunov, halo orbits and $\mathcal B_0^{\pm}$.}
	\label{bifurcation_graph_1}
\end{figure}


\subsection{Bifurcation graph between $g$ and double cover of $\mathcal{B}_0^{\pm}$}
\label{sec5.3}

At $g^{(1,1)}$ ($\Gamma = 1.383093$) the spatial index jumps from 3 to 4 (planar index is 2).\ The continuation of this new branch of spatial orbits terminates at $\mathcal B_0^{(2,2)}$ ($\Gamma = -0.182001$), where corresponding pair of Floquet multipliers become from negative hyperbolic to elliptic, hence its double cover before the period-doubling bifurcation are bad orbits (index jumps from $\overline{5}$ to 4).\ The bifurcation graph and some orbits are shown in Figure \ref{bifurcation_graph_3}, and their data are collected in Table \ref{data_g_double_B0}.\ The orbits are doubly symmetric with respect to the $XOZ$- and $YOZ$-plane, hence they are invariant under $-\sigma$ ($OZ$-symmetry).\ The symmetry $\sigma$ (reflection at the ecliptic) yields the symmetric family.\ Their properties are summarized as follows:\
\begin{center}
    \begin{tabular}{ccc}
         $\Gamma$ & stability behavior & $\mu_{CZ}$ \\
         $(1.383093 , 0.702961)$ & pos.\ hyperbolic ($-$) \&  elliptic ($+$) & $5$ \\
         $(0.702961 , -0.072698)$ & pos.\ hyperbolic ($-$) \&  elliptic ($-$) & $5$ \\
        $(-0.072698, -0.074122)$ & elliptic ($+$) \&  elliptic ($-$) & $4$ \\
        $(-0.074122 , -0.135863)$ & complex instability & $4$ \\
        $(-0.135863 , -0.182001)$ & elliptic ($+$) \&  elliptic ($-$) & $4$
    \end{tabular}
\end{center}
\noindent
Stable regions are provided in the Jacobi constant range $(-0.072698, -0.074122)$ and also in the range $(-0.135863 , -0.182001)$, where orbits in the latter are close to $\mathcal B_0^{(2,2)}$.\

\begin{figure}[t!]
\centering
\includegraphics[width=1\linewidth]{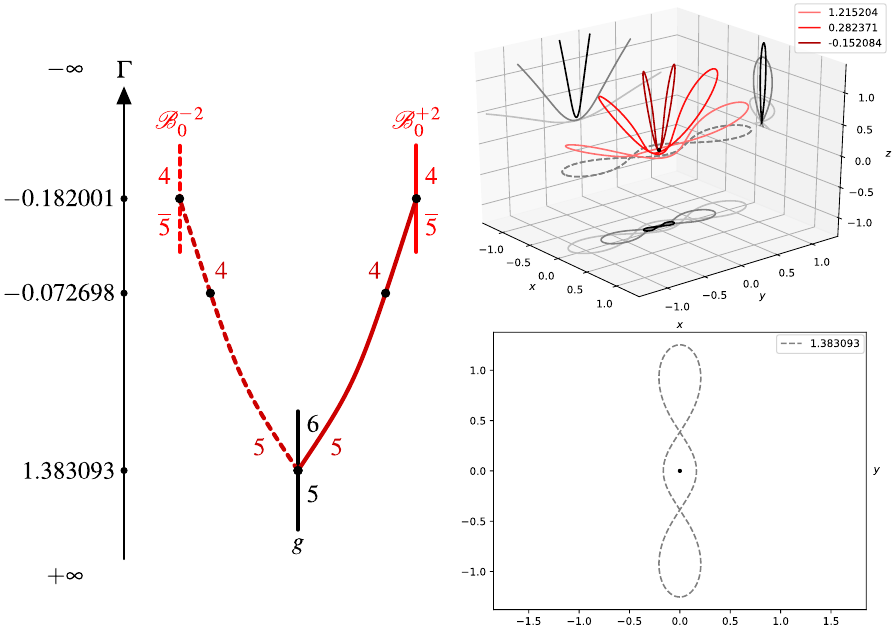}
\caption{Left:\ Bifurcation graph between $g^{(1,1)}$ ($\Gamma = 1.383093$) with $\mathcal B_0^{(2,2)}$ ($\Gamma = -0.182001$); in between is a degenerate orbit at $\Gamma = -0.072698$.\ Right:\ Corresponding orbits, where $g^{(1,1)}$ is gray dashed.}
\label{bifurcation_graph_3}
\end{figure}

\subsection{Bifurcation graph between $g'$ and double cover of halo orbits}
\label{sec:g'_halo}

The critical point $g'^{(1,1)}$ gives rise to a family of spatial orbits at $\Gamma = 3.390159$ where the index jumps from 6 to~7 (planar index is 3 and spatial index jumps from 3 to 4).\ The orbits of this branch of spatial orbits are simple symmetric with respect to the $XOZ$-plane, and they terminate at $\text{halo}^{(2,1)}$ ($\Gamma = 1.095146$).\ At this bifurcation point the double cover of $\text{halo}^{(2,1)}$ are bad orbits before the bifurcation point, because the associated Floquet multipliers become from negative hyperbolic to elliptic (index jump from $\overline{7}$ to 6).\ The corresponding bifurcation graph and some orbits are shown in Figure~\ref{bifurcation_graph_4}.\ Based on their data given in Table \ref{data_g'_double_halo}, their properties are:\
\begin{center}
    \begin{tabular}{ccc}
         $\Gamma$ & stability behavior & $\mu_{CZ}$ \\
         $(3.390159 , 1.126645)$ & neg.\ hyperbolic ($+$) \&  elliptic ($+$) & $6$ \\
         $(1.126645 , 1.115248)$ & neg.\ hyperbolic ($+$) \&  neg. hyperbolic ($-$) & $6$ \\
         $(1.115248 , 1.105626)$ & complex instability & $6$ \\
         $(1.105626 , 1.095146)$ & elliptic ($-$) \&  elliptic ($+$) & $6$
    \end{tabular}
\end{center}
\noindent
Very close to the double cover of $\text{halo}^{(2,1)}$ in the Jacobi constant range $(1.105626 , 1.095146)$ this new family of spatial orbits enters the region of stability.\

\begin{figure}[t!]
\centering
\includegraphics[width=1\linewidth]{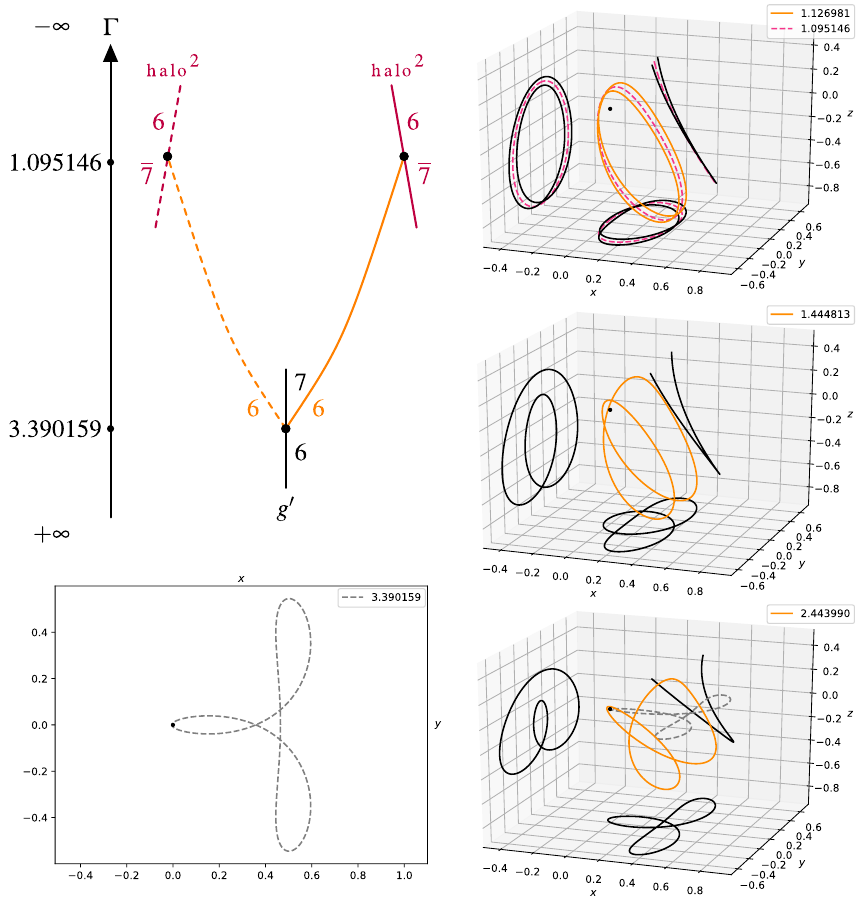}
\caption{Left top:\ Bifurcation graph between $g'^{(1,1)}$ ($\Gamma = 3.390159$) with $\text{halo}^{(2,1)}$ ($\Gamma = 1.095146$).\ Corresponding orbits start bottom left, then right, then up.\ The same graph with the same topological properties holds also for the symmetric families associated to $g'$ and double cover of halo orbits.}
\label{bifurcation_graph_4}
\end{figure}

\subsection{Bifurcation graph between double cover of $g$ and $g'$, and 3rd cover of $\mathcal{B}_0^{\pm}$}
\label{sec5.5}

Our next result is associated to the relations between $g^{(2,1)}$, $g'^{(2,1)}$, $g'^{(2,2)}$ and $\mathcal{B}_0^{(3,1)}$, as illustrated in Figure \ref{bifurcation_graph_7}.\ We use in the following description the same color for the families as in the bifurcation graph.\

At $\Gamma = 3.057470$ the index of the double cover of $g^{(2,1)}$ jumps from 9 to 11.\ At this bifurcation point the spatial Floquet multipliers of $g^{(2,1)}$ are at $-1$ (spatial index is 3) and avoid transition to neg.\ hyperbolicity (point of tangency of vertical stability index), hence the spatial index of its double cover jumps from 5 to 7.\ The planar Floquet multipliers are pos.\ hyperbolic with index 2, so its double cover has a planar index of 4.\ This implies the index jump from 9 to 11.\ For the generated two families, the orbits of one family start with index 9 and are doubly symmetric with respect to the $OX$-axis and $XOZ$-plane (blue family), and the orbits of the other family start with index 10 and are doubly symmetric with respect to the $OY$-axis and $YOZ$-plane (green family).\ Their data are collected in Table~\ref{data_g2_g'2} and some orbits are plotted in Figure~\ref{figure_g2_g'2}.\ These families were described in~\cite{batkhin2009} and denoted as $g_{1v}^{XOZ}$ and $g_{1v}^{YOZ}$ correspondingly. The last family is an open one.

\begin{figure}[t!]
			\centering
			\includegraphics[width=1\linewidth]{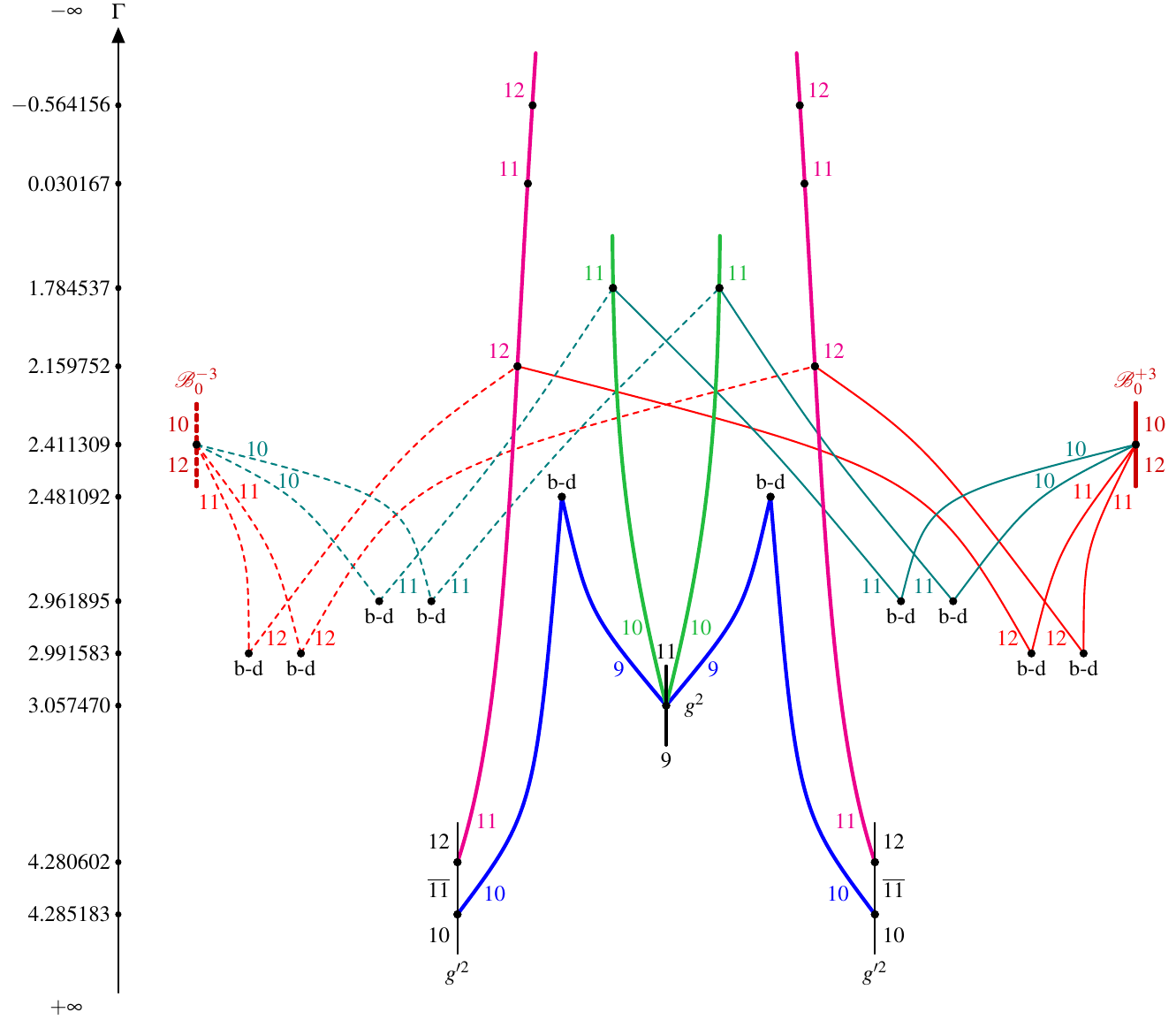}
 			\caption{Bifurcation graph between $g'^{(2,1)}$ ($\Gamma = 4.285183$), $g'^{(2,2)}$ ($\Gamma = 4.280602$), $g^{(2,1)}$ ($\Gamma = 3.057470$) and $\mathcal B_0^{(3,1)}$ ($\Gamma = 2.411309$).\ The corresponding orbits with the same colors are plotted in Figure \ref{figure_g2_g'2}, Figure~\ref{figure_g'2} and Figure \ref{figure_B03}.\ Stable regions can be seen in their data base in Table \ref{data_g2_g'2} and in Table \ref{data_B03}.}
			\label{bifurcation_graph_7}
\end{figure}

\begin{figure}[t!]
	\centering
	\includegraphics[width=0.85\linewidth]{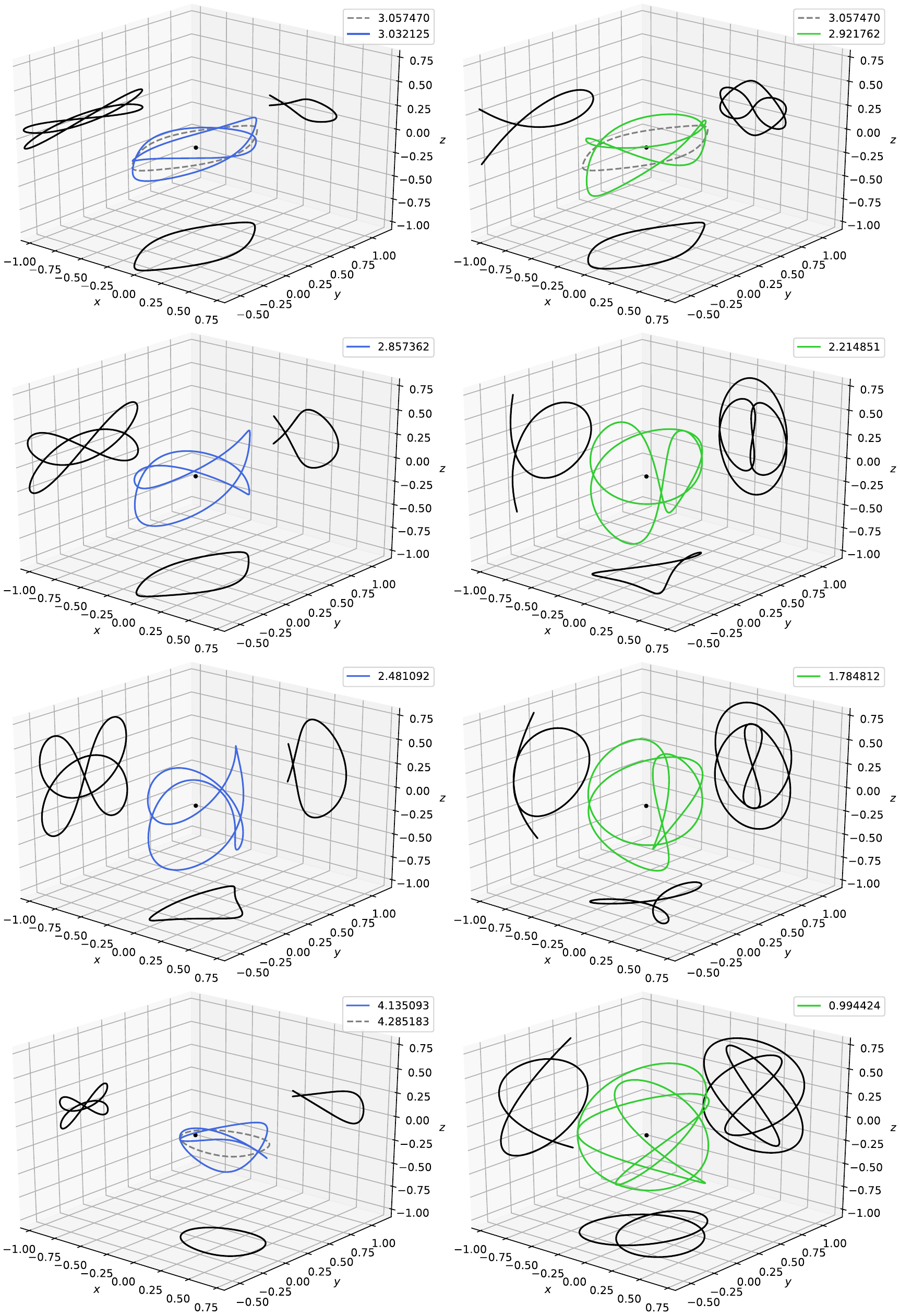}
	\caption{Blue and green orbits (from top to bottom) from Figure \ref{bifurcation_graph_7}, branching out from $g^{(2,1)}$ (gray dashed at the top).\ Blue family terminates at $g'^{(2,1)}$ (gray dashed at the bottom left).\ Third green orbit corresponds to symmetry-breaking bifurcation point at $\Gamma = 1.784812$ related to aqua family branching out from 3rd cover of $\mathcal B_0^{(3,1)}$, see Figure \ref{figure_B03}.}
	\label{figure_g2_g'2}
\end{figure}

The $g'$-orbits in the very small $\Gamma$-range $(4.285183,4.280602)$ between $g'^{(2,1)}$ and $g'^{(2,2)}$ are spatial negative hyperbolic, being spatial elliptic before (neg.\ $B$-sign) and after (pos.\ $B$-sign).\ The spatial index is~3 and the planar Floquet multipliers are elliptic with rotation angle less then~$\pi$, therefore the double cover of $g'^{(2,1)}$ jumps from 10 to~$\overline{11}$, and the double cover of~$g'^{(2,2)}$ jumps from $\overline{11}$ to 12.\ The first family starts with index 10 (blue family) and the second family, which is an open one, starts with index 11 (purple family).\ The orbits of both families are doubly symmetric with respect to the $OX$-axis and $XOZ$-plane.\ Their data can be found in Table \ref{data_g2_g'2} and some purple orbits are plotted in Figure \ref{figure_g'2}.\

\begin{figure}[t!]
	\centering
	\includegraphics[width=0.9\linewidth]{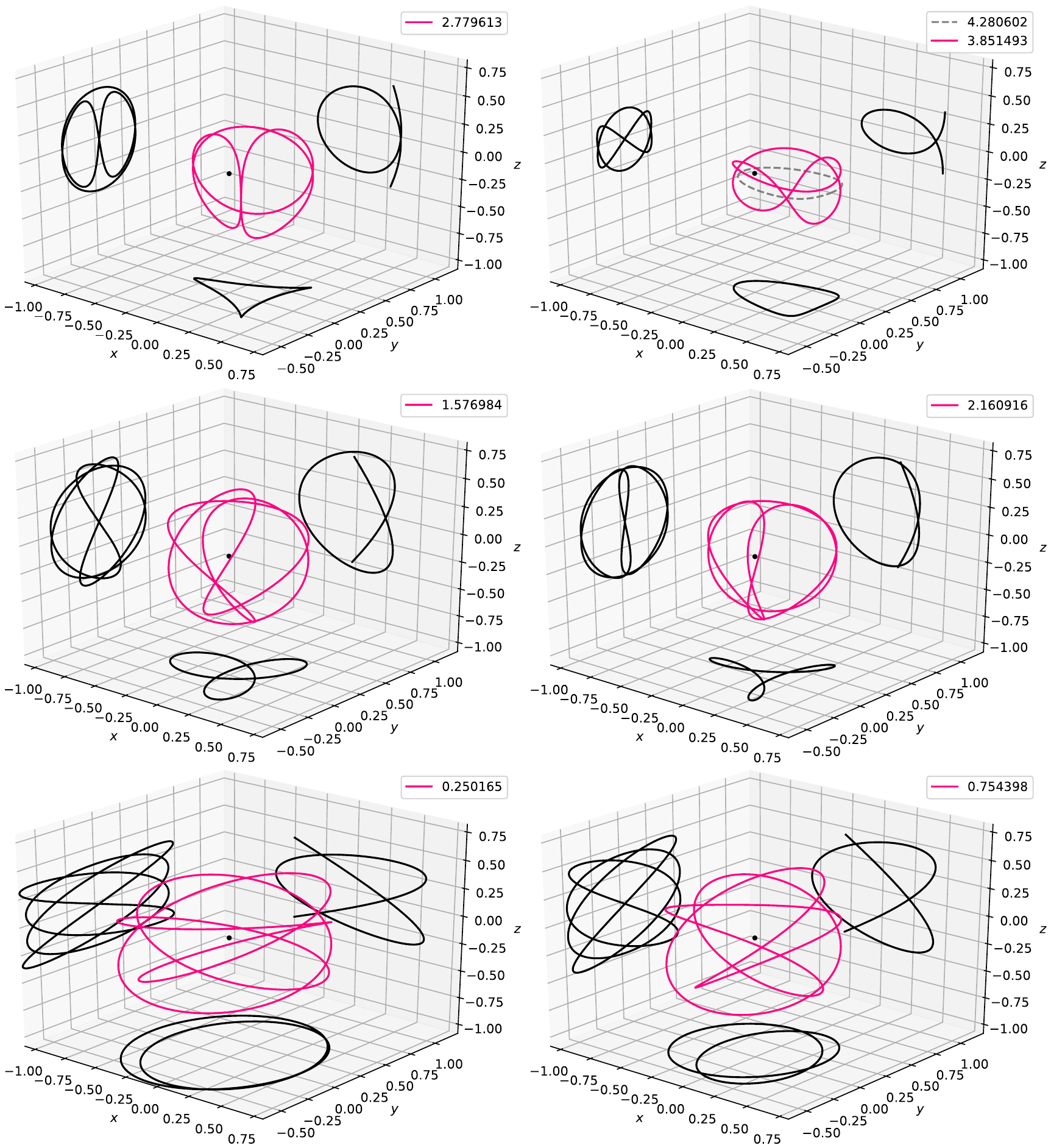}
	\caption{Plot of purple orbits from bifurcation graph \ref{bifurcation_graph_7}, branching out from $g'^2$ (gray dashed on top right).\ Third orbit (middle right) corresponds to symmetry-breaking bifurcation point at $\Gamma = 2.160916$, interacting with red family branching out from third cover of $\mathcal B_0^{+}$, see Figure \ref{figure_B03}.}
	\label{figure_g'2}
\end{figure}

Our computations show that both blue families coincide with each other; if one continues at one of the critical orbits $g^{(2,1)}$ and $g'^{(2,1)}$, after a birth-death bifurcation, one terminates at the other critical orbit.\ This connection has already been computed in \cite{michalodimitrakis} and appears in a bifurcation graph in \cite{aydin_cz}.\ Each second new family (green and purple) is not closed, but owns index jumps.\ The green family makes an index jump from 10 to 11 at $\Gamma = 1.784537$; the index of purple family jumps from 11 to 12 at $\Gamma = 2.159752$.\

\begin{figure}[t!]
	\centering
	\includegraphics[width=0.85\linewidth]{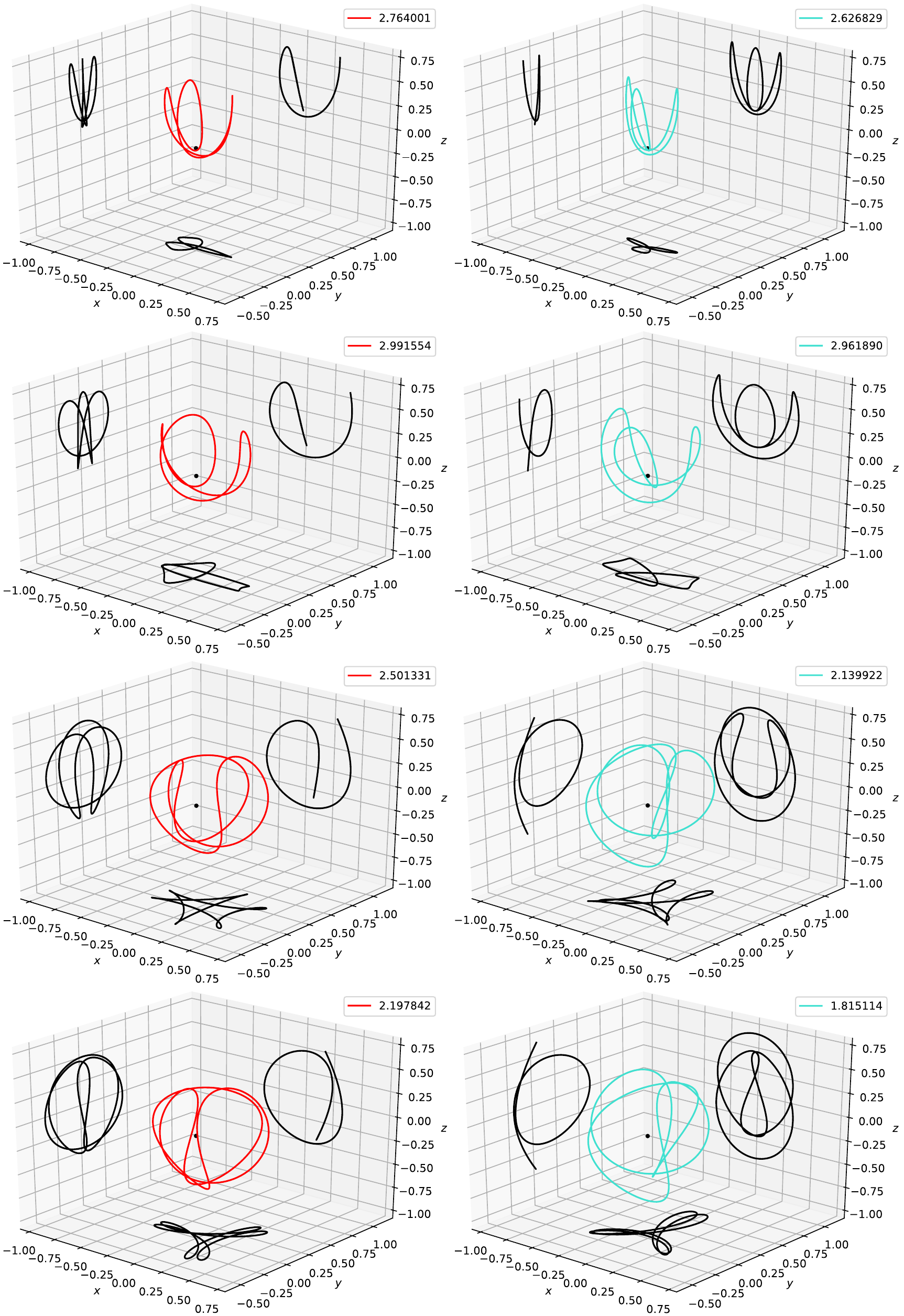}
	\caption{Plot of red and aqua orbits (from top to bottom) from bifurcation graph \ref{bifurcation_graph_7}, branching out from third cover of $\mathcal B_0^{+(3,1)}$.\ Red family terminates at purple family (see Figure \ref{figure_g'2} at $\Gamma = 2.160916$) and aqua family terminates at green family (see Figure \ref{figure_g2_g'2} at $\Gamma = 1.784812$).}
	\label{figure_B03}
\end{figure}

Now it is appropriate to emphasize that Conley--Zehnder indices have been a great support to deduce connections of families branching out from $\mathcal B_0^{(3,1)}$ to green and purple families.\ The period-tripling bifurcation of $\mathcal B_0^{(3,1)}$ generates two families, where the index jumps from 12 to 10.\ The orbits of one family are simple symmetric with respect to the $YOZ$-plane and start with index 10 (aqua family).\ The orbits of the other family are simple symmetric with respect to the $XOZ$-plane and start with index 11 (red family).\ Their data are collected in Table \ref{data_B03} and some orbits are plotted in Figure \ref{figure_B03}.\ Their continuations show that if one starts at $\mathcal B_0^{+(3,1)}$ then the termination is at $\mathcal B_0^{-(3,1)}$, and vice versa.\ This shows especially that in between (at the point going back to the symmetric critical orbit) there is a relation to other families, which is hard to see by pure computation.\ By contemplating the indices and symmetry properties, we grasp that aqua family interacts with green family at $\Gamma = 1.784537$ and red family meets purple family at $\Gamma = 2.159752$, as shown in the bifurcation graph in Figure~\ref{bifurcation_graph_7}.\ Both bifurcation points correspond to symmetry-breaking pitchfork bifurcations.\


\subsection{Family $\mathbf{f}_3$}
\label{sec:f_3}

\begin{figure}[t!]
\centering
\includegraphics[width=0.75\linewidth]{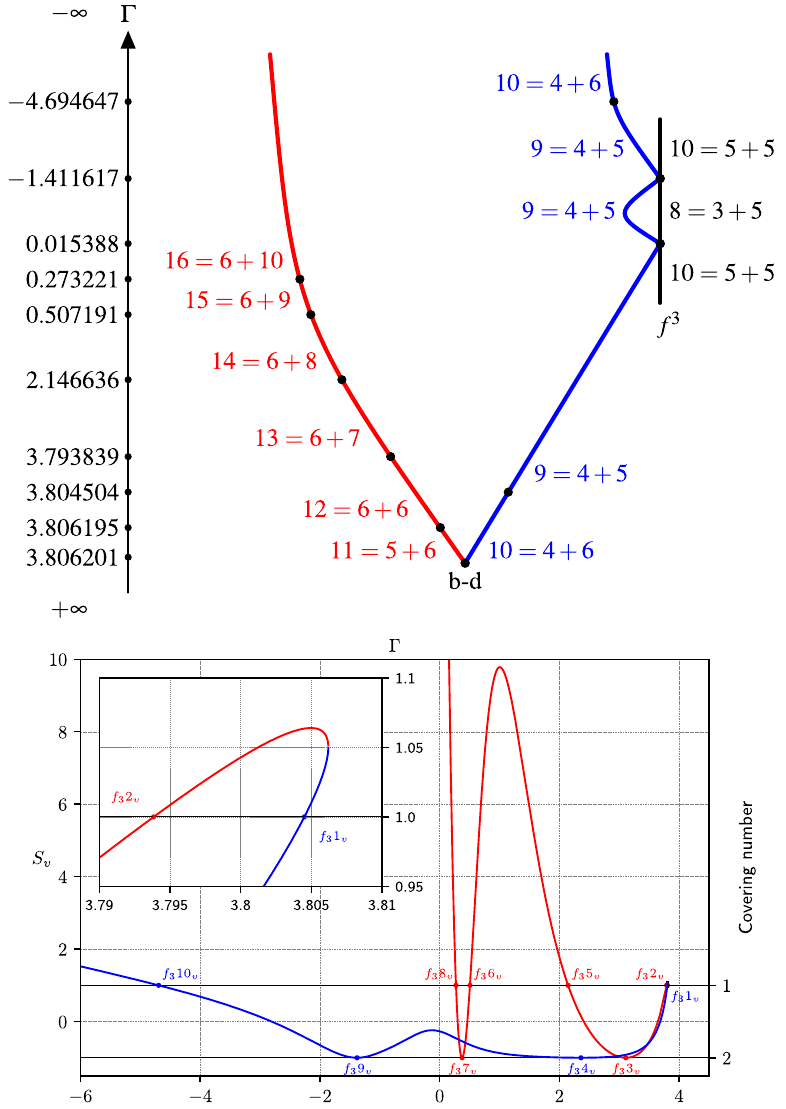}
\caption{Top:\ Bifurcation graph corresponding to family $\mathbf{f}_3$ with indices $\mu_{CZ} = \mu_{CZ}^p + \mu_{CZ}^s$.\ Bottom:\ Plot of vertical stability index of blue and red branches of family $\mathbf{f}_3$.}
\label{bifurcation_graph_9}
\end{figure}


Let us give a brief description of the family $\mathbf{f}_3$ and its key features.\ The orbits of $\mathbf{f}_3$ are doubly symmetric, that is, they are symmetric with respect to both coordinate axes $OX$ and $OY$.\ There are two branches along each of them that $\Gamma$ tend to $-\infty$.\ One branch drawn in blue in Fig.~\ref{bifurcation_graph_9} tends to the limiting arc-solution denoted by $\{+1,-1\}$ (for detail of this notation see~\cite{Henon2003,BatkhinCosRes13I}).\ There are two intersections with the family $f$ for $\Gamma=-1.411618$ and $\Gamma=0.015388$ with 3-covering of the latter.\ For $\Gamma=2.238611$ there is a collision orbit.\ The other branch of the family drawn in red in~\ref{bifurcation_graph_9} tends to the limiting arc-solution denoted $\{+1,i,-1,e\}$.\ The two branches meet with each other at the values $\Gamma=3.806201$.\ The data base of the $f_3$ family is collected in Table \ref{data_f_3} and some orbits for certain $\Gamma$-values are shown in Fig.~\ref{figure_f3}.

The family $\mathbf{f}_3$ is almost planar unstable except for a very narrow interval of values~\cite{BatkhinPCS2020}.\ In contrast, the vertical stability index $S_v$ falls into the interval $[-1;+1]$ on a rather wide range of changing $\Gamma$.\ Ten critical orbits of $\mathbf{f}_3$ were found at which there exist intersections with the families of spatial periodic solutions.

\begin{figure}[t!]
\centering
\includegraphics[width=0.9\linewidth]{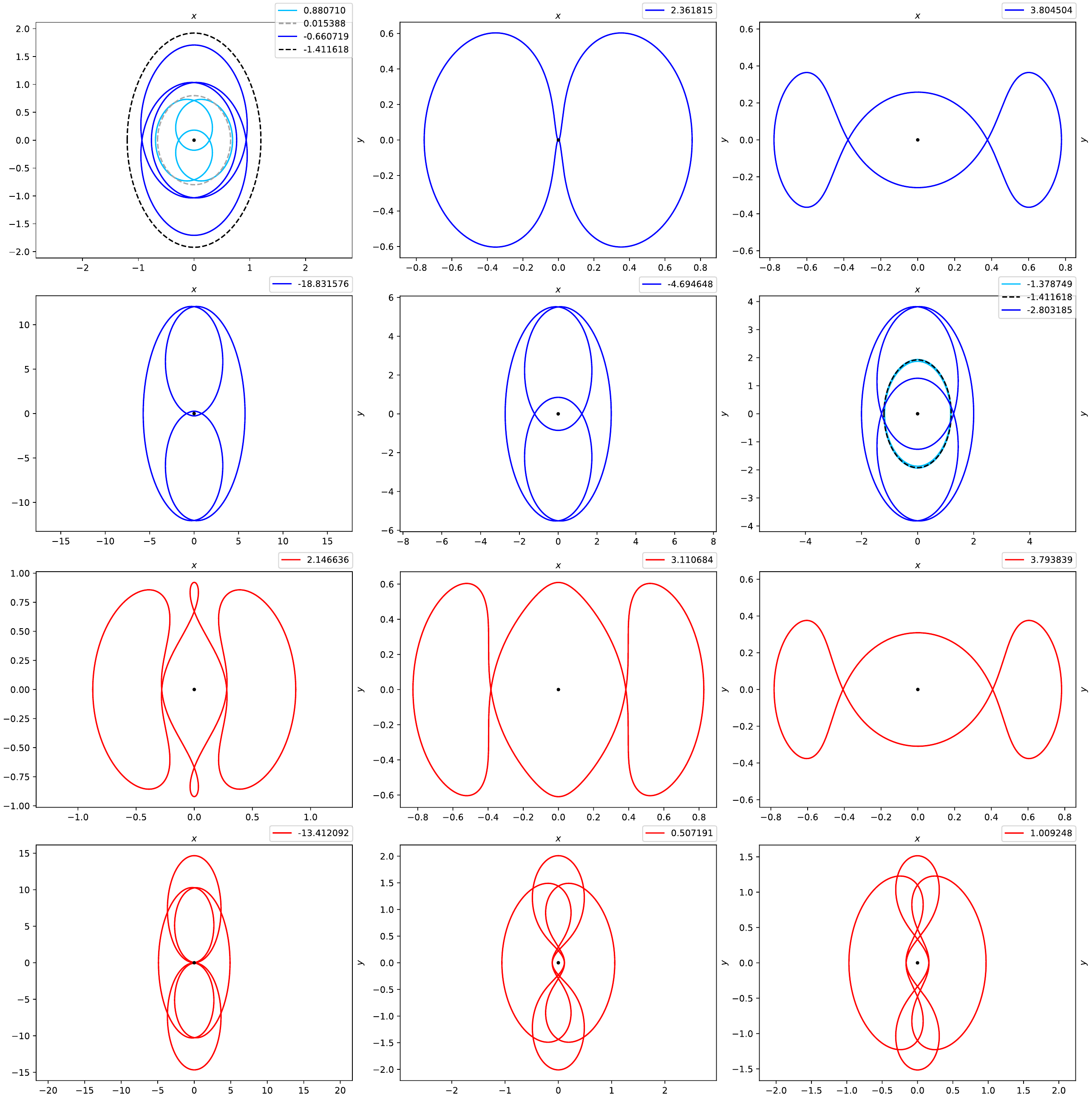}
\caption{Plot of blue and red branches from bifurcation graph \ref{bifurcation_graph_9}, related to family $\mathbf{f}_3$.\ Both branches start close to the common orbit of birth-death type.\ For the blue branch:\ The gray dashed orbit at $\Gamma = 0.015388$ corresponds to $f_p^{(3,1)}$ and the black dashed at $\Gamma = -1.411618$ corresponds to $f_p^{(3,2)}$.}
\label{figure_f3}
\end{figure}

The behavior of the vertical stability index $S_v$ of the $\mathbf{f}_3$ family is rather complicated, and at large intervals of variation of the Jacobi constant $\Gamma$ the orbits of this family interact with a large number of spatial families.\ To denote these families, we will use the same approach as in~\cite{henon74}, where the generating orbits of planar families were denoted by the indices $\mathbf{F}n_v$, where $n$ is the number of the orbit ordered by decreasing parameter $\Gamma$.\ We provide an overview description of the generated spatial families, paying attention only to some peculiarities.\ From Fig.~\ref{figure_f3} we can see that the family $\mathbf{f}_3$ has 6 \textit{vertical critical} (VC) solutions with 1 covering, which generate six spatial families (critical orbits $\mathbf{f}_31_v$, $\mathbf{f}_32_v$, $\mathbf{f}_35_v$, $\mathbf{f}_36_v$, $\mathbf{f}_38_v$, $\mathbf{f}_310_v$) and 4 VC with 2 covering (critical orbits $\mathbf{f}_33_v$, $\mathbf{f}_34_v$, $\mathbf{f}_37_v$, $\mathbf{f}_39_v$).\ In the last cases, each of the VC is the solution of intersection of four families of spatial orbits.\ These four families form two pairs of mutually symmetric orbits with the same dynamical characteristics, so, it is enough to describe only one representative in each pair.


In the VC $\mathbf{f}_31_v$ there is an intersection at $\Gamma\approx 3.80450408$ with the family of spatial periodic solutions of doubly symmetric orbits with respect to the plane $XOZ$ and the axis $OY$.\ This family is open, i.e. it tends to $\Gamma\to-\infty$ and does not have stable orbits.  

In the VC $\mathbf{f}_32_v$ there is an intersection at $\Gamma\approx3.79383981$ with the family of spatial periodic solutions of doubly symmetric orbits with respect to the plane $YOZ$ and the axis $OX$.\ This family is also open and does not have stable orbits.  

As mentioned above, the VC solution $\mathbf{f}_33_v$ is the common solution for two pairs of families of periodic spatial orbits at $\Gamma\approx 3.11068444$. One pair consists of doubly symmetric orbits with respect to the plane $XOZ$ and the axis $OX$ and the other with respect to the plane $YOZ$ and the axis $OY$. Both pairs of families demonstrate the same behavior and are open.\ There are some intervals of vertical stability along the families.

For VC solution $\mathbf{f}_34_v$ at $\Gamma\approx2.36181537$ also exist two pairs of families with symmetric orbits as in the previous case.\ Only one family has been studied. This family demonstrates nonmonotonic changing of the integral $\Gamma$ along the family, but finally it tends to $\Gamma\to-\infty$.

The behavior of family of doubly symmetric periodic solution, which share the common VC $\mathbf{f}_35_v$ at $\Gamma\approx2.14663695$ drastically differs from the behavior of other families.\ Orbits of this family is symmetric with respect to coordinate planes $XOZ$ and $YOZ$.\ The family is closed, i.e. it exists  on a finite interval of $\Gamma\in[0.0792;2.147]$. At $\Gamma\approx 0.75514$ it has the 5th covering of the family $f$ (see Section \ref{sec:5.7}).

In the VC $\mathbf{f}_36_v$ there is an intersection at $\Gamma\approx0.50719121$ with the family of spatial periodic solutions of doubly symmetric orbits with respect to the axis $OX$ and $OY$.\ This family is open and it does not have stable orbits.  

In the VC $\mathbf{f}_37_v$ there is an intersection at $\Gamma\approx0.37546334$ with two families of doubly symmetric spatial periodic solutions.\ The first one consists of orbits symmetrical with respect to the plane $XOZ$ and to the axes $OX$ and the second one consists of orbits symmetrical with respect to the plane $YOZ$ and to the axes $OY$.\ Both families tend to $\Gamma\to-\infty$ and their orbits are very unstable.

In the VC $\mathbf{f}_38_v$ there is an intersection at $\Gamma\approx0.273221$ with the family of spatial periodic solutions of doubly symmetric orbits with respect to the axes $OX$ and the plane $YOZ$.\ This family is open and it does not have stable orbits.

In the VC $\mathbf{f}_39_v$ there is an intersection at $\Gamma\approx-1.37874859$ with two families of spatial periodic solutions.\ The first one contains doubly symmetric orbits with respect to the plane $YOZ$ and axes $OY$, the second one -- with respect to the plane $XOZ$ and axes $OX$.\ The first family is closed and has the 6th covering of the family $f$ (see Section \ref{sec:5.8}), while the second one is open.

Finally, the last VC $\mathbf{f}_310_v$ gives at the value  $\Gamma\approx-4.69464797$ intersection with family of spatial doubly symmetric orbits with respect to the axis $OX$ and the plane $YOZ$.\  All orbits of this family are highly unstable.

\subsection{Bifurcation graph between 3rd cover of $g$ and $g'$, 4th cover of $\mathcal{B}_0^{\pm}$, 5th cover of $f$ and~$\mathbf{f}_3$}
\label{sec:5.7}

The next result is achieved by bifurcations generated by $g^{(3,1)}$, $g'^{(3,1)}$, $\mathcal{B}_0^{(4,1)}$, $f^{(5,1)}$ and $\mathbf{f}_3 5_v$.\ The corresponding bifurcation graph is illustrated in Figure \ref{bifurcation_graph_6}.\ Again as before, we describe these relations with the same color for each family as in the bifurcation graph.\

\begin{figure}[t!]
			\centering
			\includegraphics[width=0.9\linewidth]{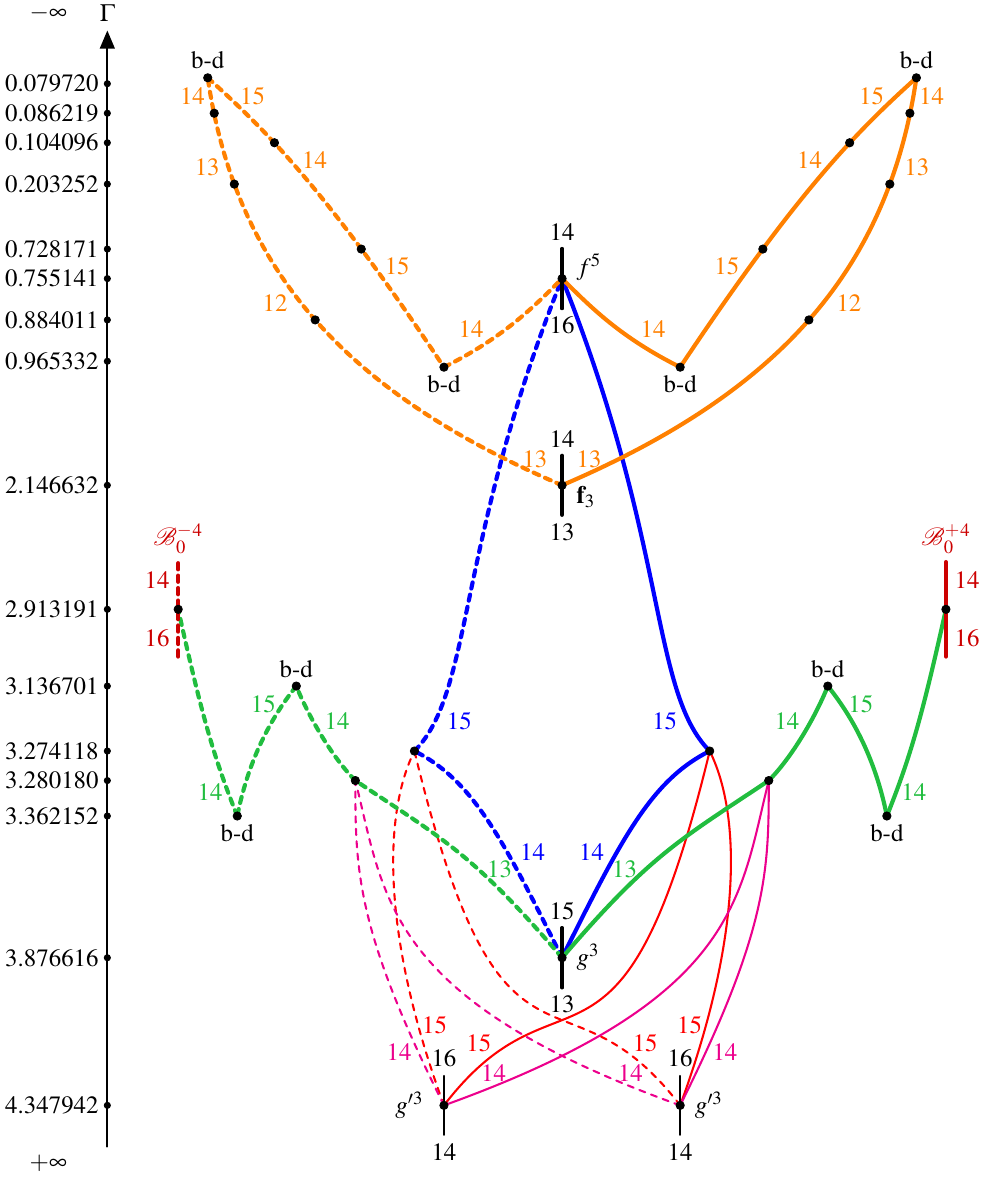}
			\caption{Bifurcation graph showing connections between $g'^{(3,1)}$ $(\Gamma = 4.347942)$, $g^{(3,1)}$ $(\Gamma = 3.876616)$, $\mathcal B_0^{(4,1)}$ $(\Gamma = 2.913191)$, $\mathbf{f}_3 5_v$ $(\Gamma = 2.146632)$ and $f^{(5,1)}$ $(\Gamma = 0.755141)$.\ Orbits with same colors are plotted in Figure \ref{figure_g3_f5_BO4} and Figure \ref{figure_g3_f5_BO4_f3}.\ Stable regions are provided in ranges of:\ both green regions with index 14, orange branch with index 14 between b-d and $f^{(5,1)}$, and some purple regions, see Table \ref{data_g3_B04_f5}.}
			\label{bifurcation_graph_6}
\end{figure}

The Conley--Zehnder index of triple covering of $g^{(3,1)}$ ($\Gamma = 3.876616$) jumps from 13 to 15 (planar index is 6 and spatial index jumps from 7 to 9).\ One family of the two new branches of spatial orbits starts with index 14 (blue family) and the other starts with index 13 (green family).\ Their data are collected in Table \ref{data_g3_B04_f5} and some orbits are plotted in Figure \ref{figure_g3_f5_BO4}.\ The orbits of the blue family are doubly symmetric with respect to the $OX$- and $OY$-axis, and they terminate at $f^{(5,1)}$ ($\Gamma = 0.755141$) whose fifth cover's index jumps from 16 to 14.\ The blue orbits undergo an index jump from 14 to 15 at $\Gamma = 3.274118$.\ The green orbits are doubly symmetric with respect to the $XOZ$- and $YOZ$-plane.\ In the non-regularized system, this family terminates at a collision (close at $\Gamma = 2.913191$), as investigated in \cite{kalantonis} which has also computed the blue family.\ Our studies in the regularized problem show that this collision point corresponds to $\mathcal{B}_0^{(4,1)}$ where the index jumps from 16 to 14.\ Furthermore, the green orbits make an index jump from 13 to 14 at $\Gamma = 3.280180$, and then they undergo two times a birth-death bifurcation, first at $\Gamma = 3.136701$ and second at $\Gamma = 3.362152$.\

\begin{figure}[t!]
	\centering
	\includegraphics[width=0.85\linewidth]{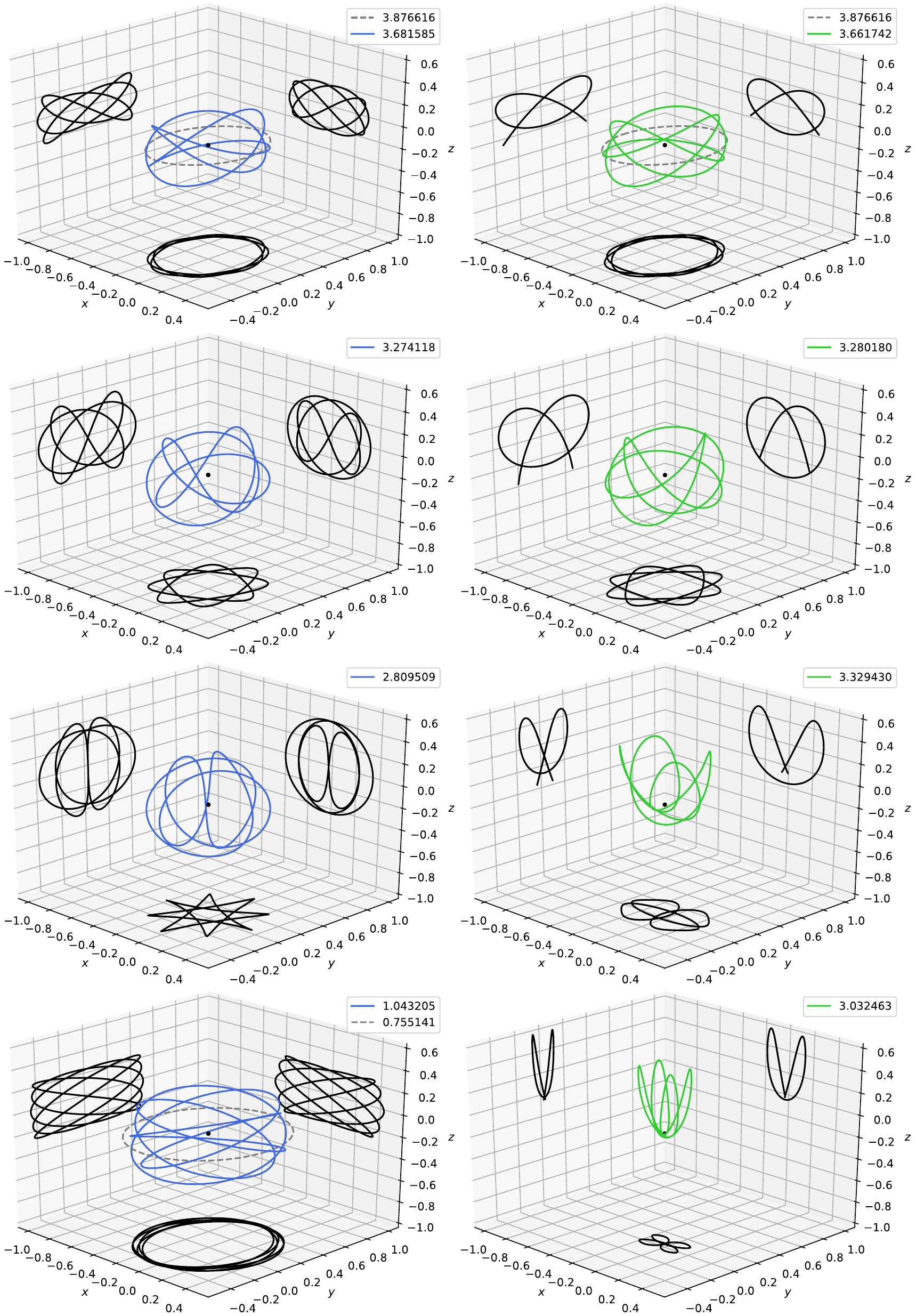}
	\caption{Blue and green orbits from Figure \ref{bifurcation_graph_6}, from top to bottom.\ Gray dashed at top:\ $g^{(3,1)}$.\ Blue family terminates at $f^{(5,1)}$ (gray dashed at bottom left) and green family terminates at $\mathcal{B}_0^{(4,1)}$.}
	\label{figure_g3_f5_BO4}
\end{figure}

\begin{figure}[ht!]
	\centering
	\includegraphics[width=0.85\linewidth]{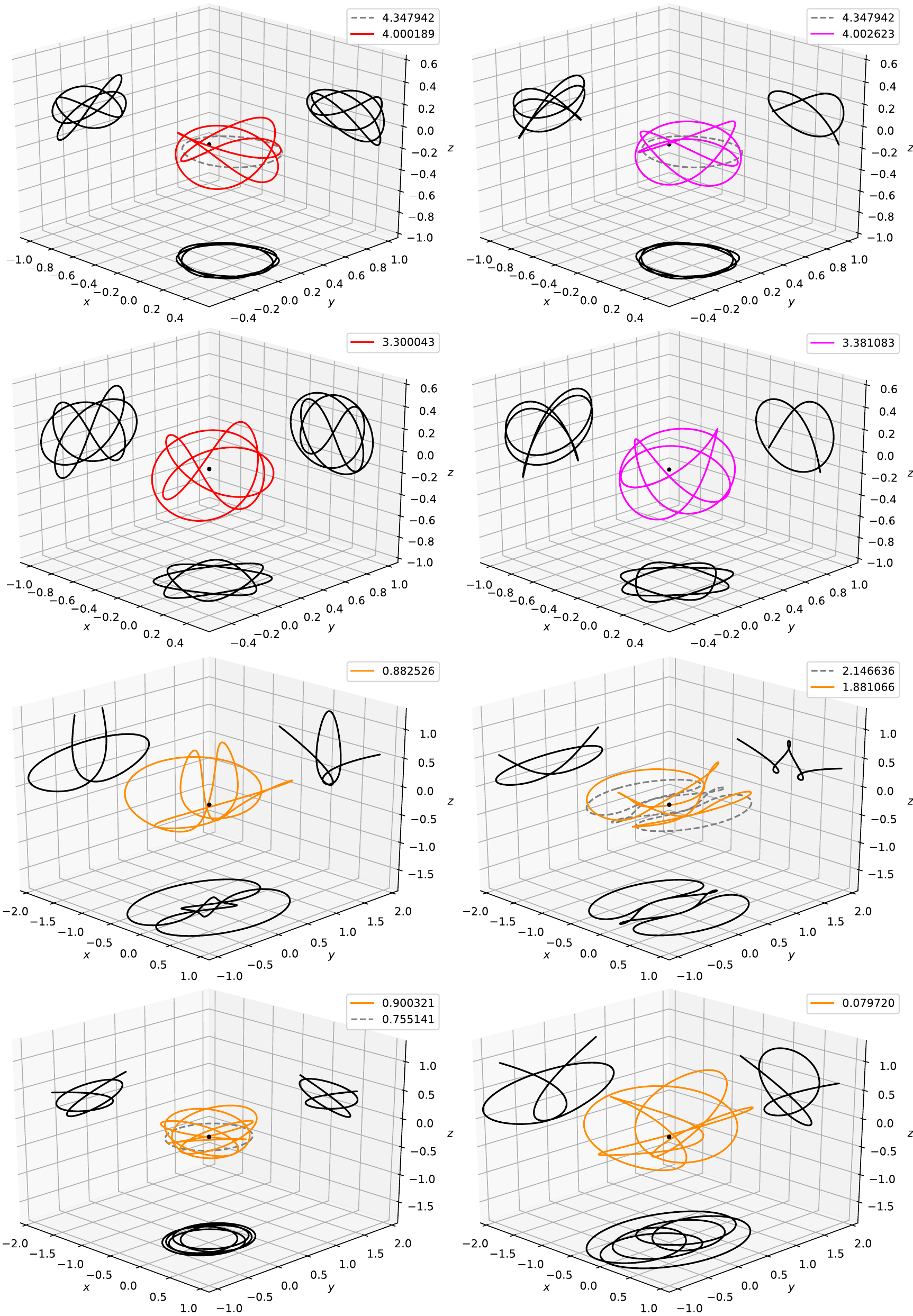}
	\caption{Red, purple and orange orbits from Figure \ref{bifurcation_graph_6}.\ Red and purple orbits branch out from $g'^{(3,1)}$ (gray dashed at the top) and terminate at blue resp.\ green family (see second row in Figure \ref{figure_g3_f5_BO4}).\ Orange orbits bifurcate from $\mathbf{f}_3 5_v$ (gray dashed in first plot) and terminate at $f^{(5,1)}$ (gray dashed in last plot).}
	\label{figure_g3_f5_BO4_f3}
\end{figure}

Triple covering of $g'^{(3,1)}$ ($\Gamma = 4.347942$) generates two families of spatial orbits, the red and purple family.\ At this bifurcation point the index jumps from 14 to 16.\ The red family has index 15 (orbits are simple symmetric with respect to the $OX$-axis) and the purple family has index 14 (orbits are simple symmetric with respect to the $OY$-axis).\ Their data can be found in Table \ref{data_g3_B04_f5} and some orbits are plotted in Figure \ref{figure_g3_f5_BO4_f3}.\ A straightforward computation shows that these branches terminate at corresponding symmetric $g'$-orbit, as shown in \cite{kalantonis}.\ It was investigated in \cite{aydin_cz}, by considering their indices and symmetry properties, that red family interacts with blue family (at $\Gamma = 3.274118$) and purple family interacts with green family (at $\Gamma = 3.280180$), such as shown in the bifurcation graph \ref{bifurcation_graph_6}.\ Each of these two bifurcation points corresponds to symmetry-breaking pitchfork bifurcation.\

The orange branch bifurcates from $\mathbf{f}_35_v$ ($\Gamma = 2.146632$) and, after two birth-death bifurcations and several index jumps, the orbits terminate at $f^{(5,1)}$.\ Their data can be found in Table \ref{data_f3_f5} and some orbits are plotted in Figure \ref{figure_g3_f5_BO4_f3}.\

The part in Figure \ref{bifurcation_graph_6} related to blue, green (without $\mathcal{B}_0^{\pm}$), red and purple family was already constructed in \cite{aydin_cz}, based on computations from \cite{kalantonis}.\ We have completed this network by exploring the families $\mathcal{B}_0^{\pm}$ and $\mathbf{f}_3$ and their relations at bifurcation points.\

\subsection{Bifurcation graph between 4th cover of $g$ and $g'$, 5th cover of $\mathcal{B}_0^{\pm}$, 6th cover of $f$ and double cover of $\mathbf{f}_3$}
\label{sec:5.8}

\begin{figure}[t!]
			\centering
			\includegraphics[width=0.85\linewidth]{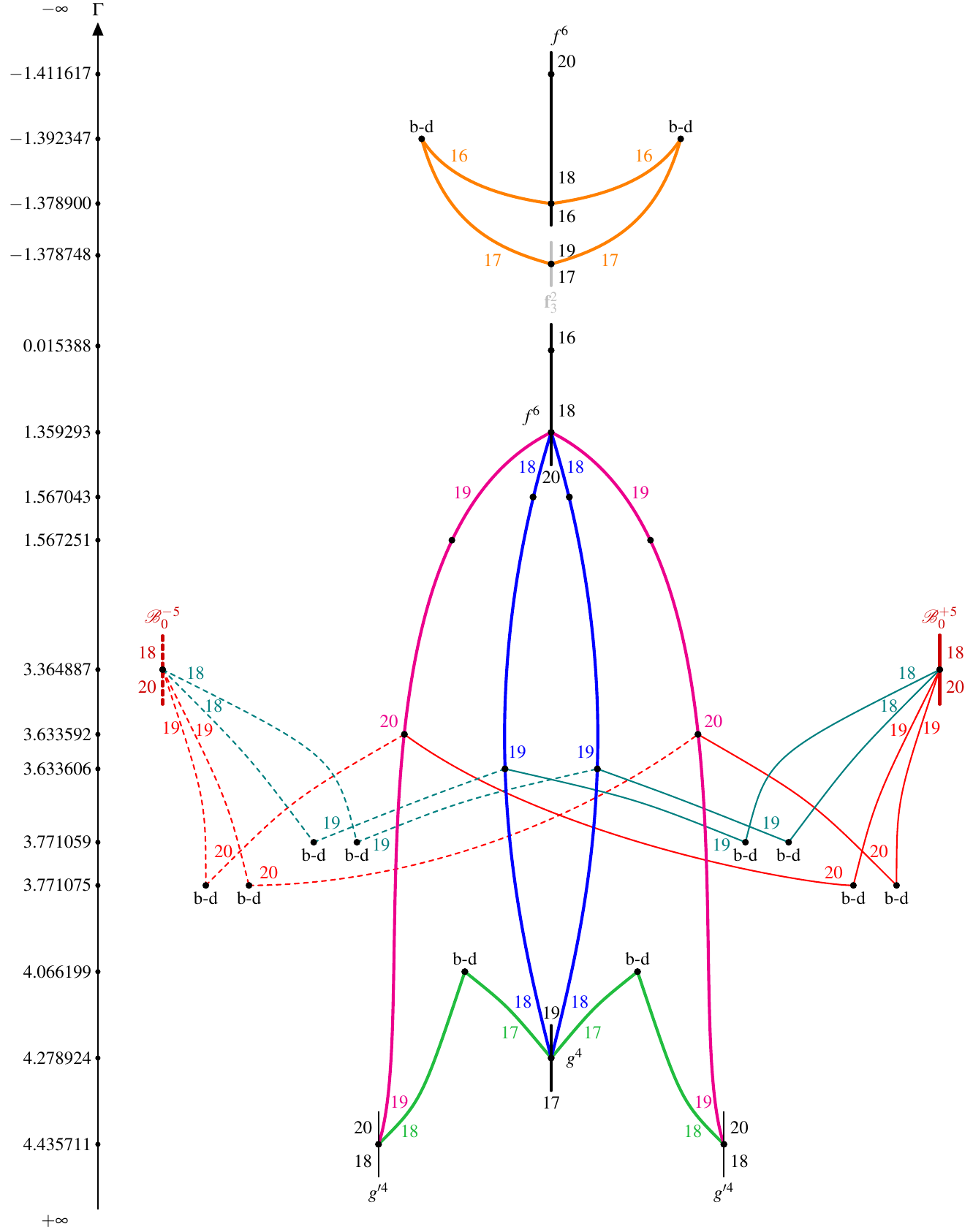}
			\caption{Bifurcation graph between $g'^{(4,1)}$ ($\Gamma = 4.435711$), $g^{(4,1)}$ ($\Gamma = 4.278924$), $\mathcal B_0^{(5,1)}$ ($\Gamma = 3.364887$), $f^{(6,1)}$ ($\Gamma = 1.359293$), $\mathbf{f}_3 9_v$ ($\Gamma = -1.378748$) and $f^{(6,2)}$ ($\Gamma = -1.378900$).\ Degenerate orbits at $\Gamma = 0.015388$ and $\Gamma = -1.411617$ correspond to double covering of $f_p^{(3,1)}$ and $f_p^{(3,2)}$.}
			\label{bifurcation_graph_8}
\end{figure}

\begin{figure}[t!]
	\centering
	\includegraphics[width=0.85\linewidth]{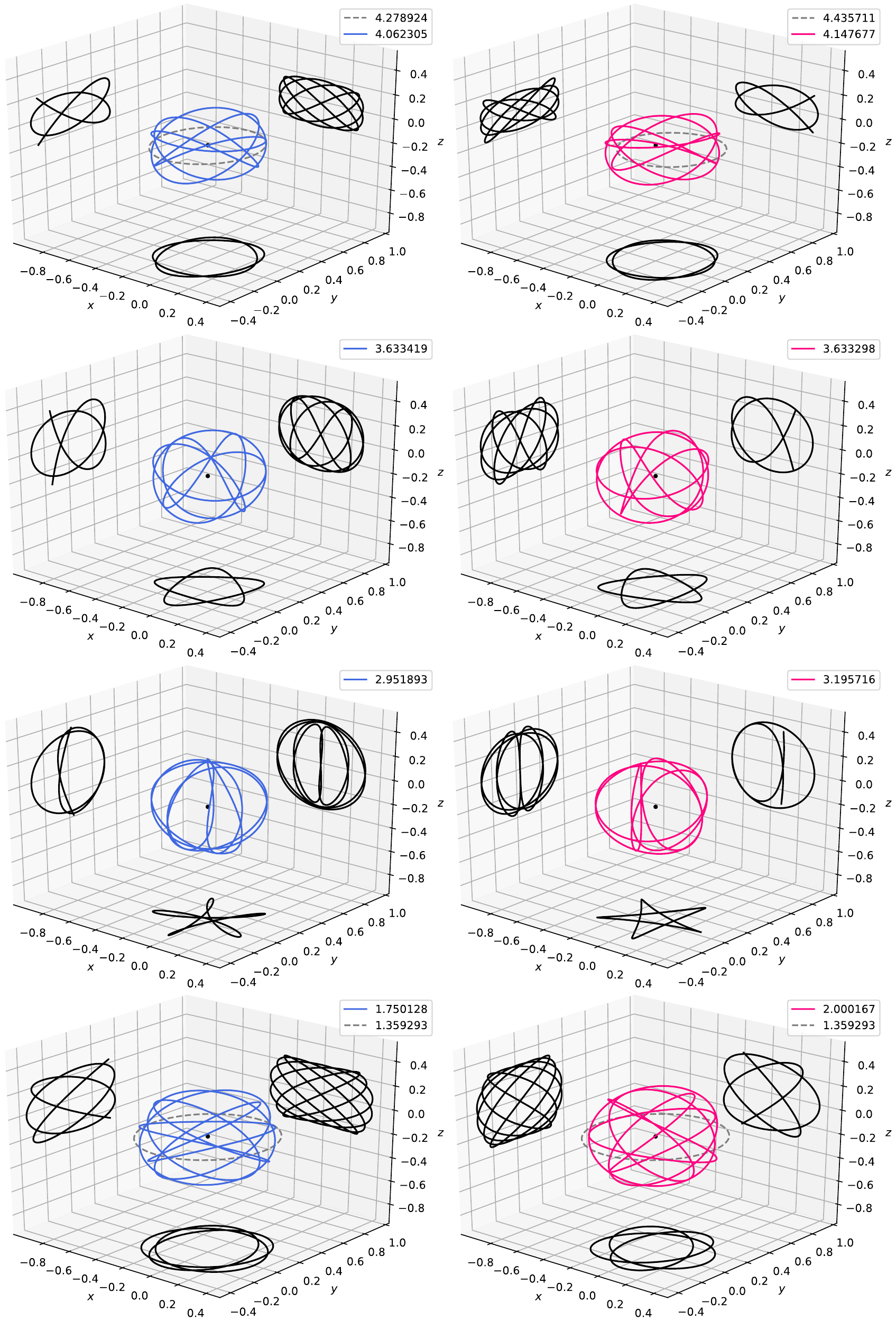}
	\caption{Blue (starting at $g^{(4,1)}$) and purple (starting at $g'^{(4,1)}$) orbits from Figure \ref{bifurcation_graph_8}, from top to bottom.\ Both families terminate at $f^{(6,1)}$.\ Gray dashed orbits (top and bottom) are corresponding planar critical orbits.}
	\label{figure_g4_g'4_f6}
\end{figure}

\begin{figure}[t!]
	\centering
	\includegraphics[width=0.85\linewidth]{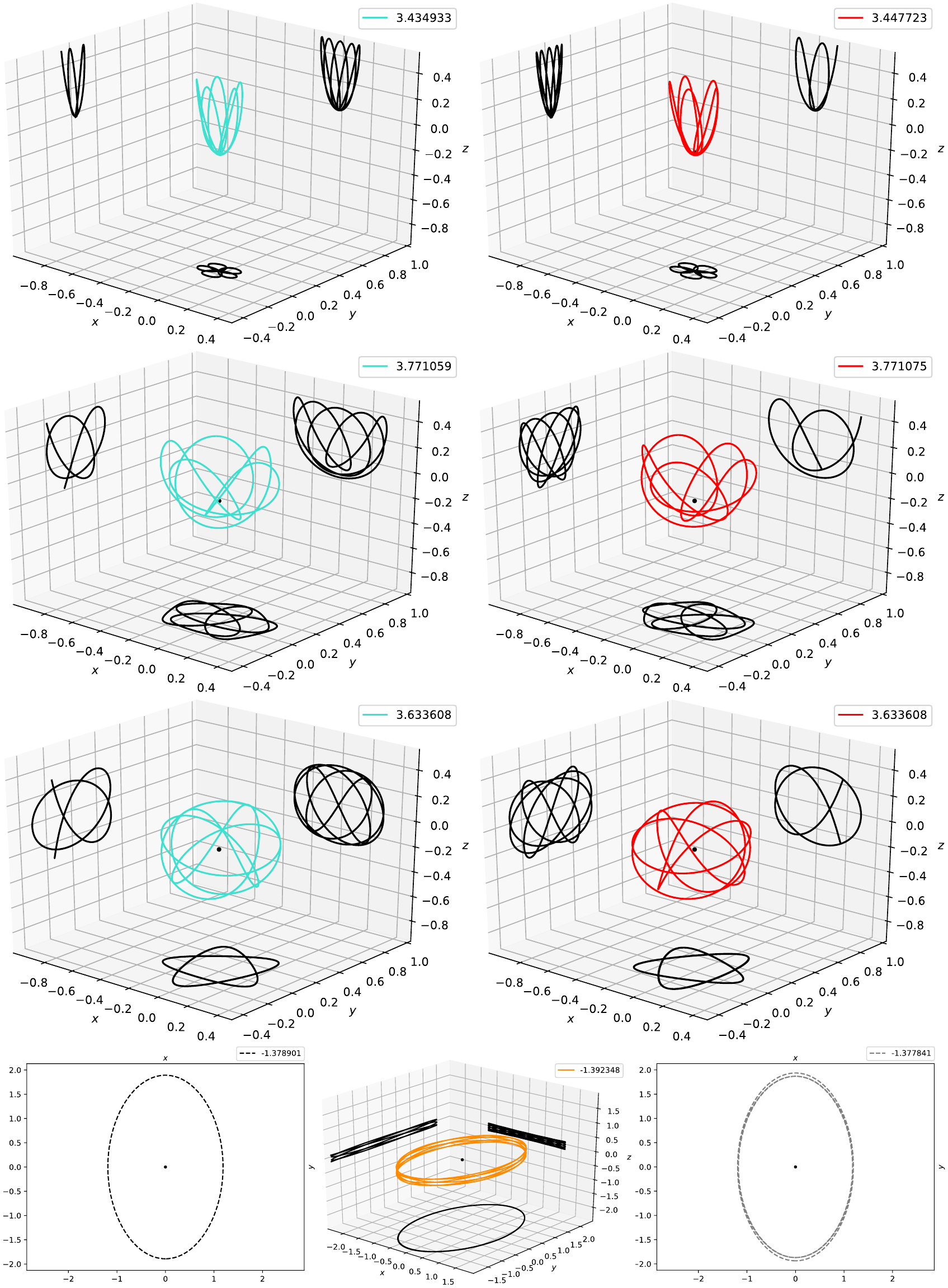}
	\caption{Aqua, red and orange orbits from Figure \ref{bifurcation_graph_8}.\ Aqua and red orbits branch out from $\mathcal B_0^{(5,1)}$ and terminate at blue resp.\ purple family (see second row in Figure \ref{figure_g4_g'4_f6}).\ The orange orbit in the middle of the last row is close to birth-death type related to the connection between $\mathbf{f}_3 9_v$ (right) and $f^{(6,2)}$ (left).}
	\label{figure_B05_f32}
\end{figure}

\begin{figure}[t!]
	\centering
	\includegraphics[width=1\linewidth]{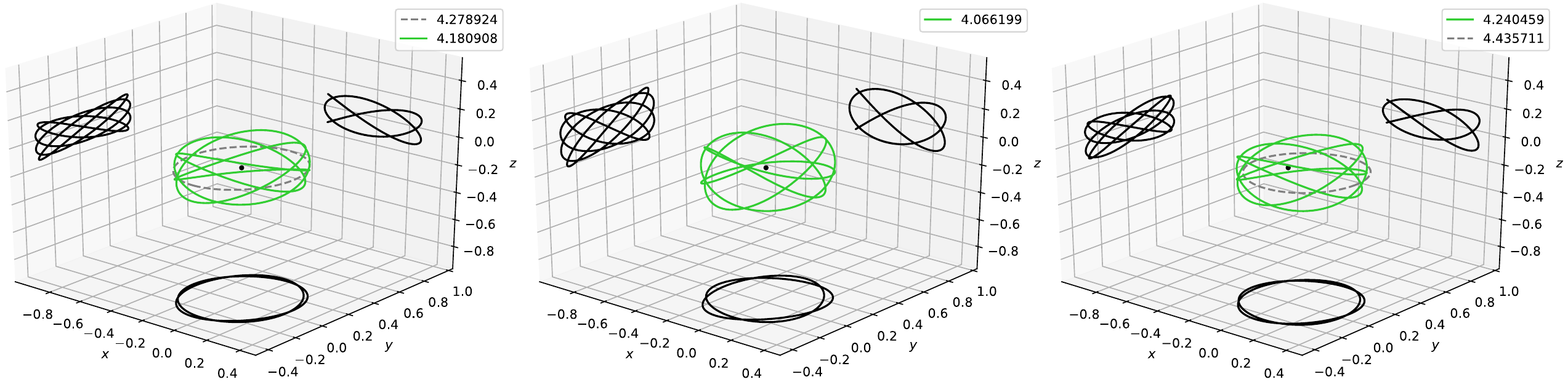}
	\caption{Green family from Figure \ref{bifurcation_graph_8}.\ Left orbit branches out from $g^{(4,1)}$, right orbit branches out from $g'^{(4,1)}$ and middle orbit is of birth-death type in between.}
	\label{figure_g4_g'4}
\end{figure}

Similar to the bifurcation result in Figure \ref{bifurcation_graph_7} associated to $g^{(2,1)}$, $g'^{(2,1)}$, $g'^{(2,2)}$ and $\mathcal B_0^{(3,1)}$, we have derived connections illustrated in Figure \ref{bifurcation_graph_8} related to $g^{(4,1)}$, $g'^{(4,1)}$ and $\mathcal B_0^{(5,1)}$, plus with $f^{(6,1)}$, $f^{(6,2)}$ and $\mathbf{f}_3 9_v$.\

Fourth covering of $g^{(4,1)}$ generates two families of spatial orbits, where the index jumps from 17 to 19.\ Orbits of one family (green) start with index 17 and are doubly symmetric with respect to $OX$-axis and $XOZ$-plane.\ Orbits of the other family (blue) start with index 18 and are doubly symmetric with respect to $OY$-axis and $YOZ$-plane.\ The blue family terminates, after two index jumps, at $f^{(6,1)}$ where the index jumps from 20 to 18.\ The green family interacts, after a birth-death bifurcation, with $g'^{(4,1)}$.\ The latter critical orbit gives rise to a second family of spatial orbits (purple) which are doubly symmetric with respect to the $OX$-axis and $XOZ$-plane, and start with index 19.\ The purple family terminates, after two index jumps, at the same critical orbit as the blue branch, namely at $f^{(6,1)}$.\ The data of blue, green and purple are collected in Table \ref{data_g4_g'4_B05_f6_1}, some blue and purple orbits are plotted in Figure \ref{figure_g4_g'4_f6} and some green orbits are plotted in Figure \ref{figure_g4_g'4}.\ The blue family has a stable region after the second index jump until its termination at $f^{(6,1)}$, corresponding to the $\Gamma$-range $(1.567043,1.359293)$.\ The green family has a stable region after bifurcation from $g'^{(4,1)}$ until before birth-death bifurcation for $\Gamma \in (4.435711,4.068827)$.\ These connections were also calculated in \cite{kalantonis} and were partially described in a bifurcation graph in \cite{aydin_cz}.\

As before in the discussion to Figure \ref{bifurcation_graph_7}, we have discovered that two families branching out from $\mathcal B_0^{(5,1)}$ are related to the families bifurcation from $g^{(4,1)}$ and $g'^{(4,1)}$, namely in the same way provided by indices and symmetry properties.\ The pure computation starts at $\mathcal B_0^{(5,1)}$ and terminates at its symmetric one.\ The index of fifth covering of $\mathcal B_0^{(5,1)}$ jumps from 20 to 18, and hence generates two families (aqua and red family).\ Aqua orbits are simple symmetric with respect to the $YOZ$-plane and start with index 18, and red orbits are simple symmetric with respect to the $XOZ$-plane and start with index 19.\ In particular, aqua/red family interacts with blue/purple family, correspondingly, as shown in Figure \ref{bifurcation_graph_8}.\ Both bifurcation points correspond to symmetry-breaking pitchfork bifurcations.\ The data of aqua and red families can be found in Table \ref{data_g4_g'4_B05_f6_2} and some orbits are plotted in Figure \ref{figure_B05_f32}.\ Aqua orbits are stable until birth-death bifurcation in the $\Gamma$-range $(3.364887,3.771059)$, and red orbits are stable after birth-death bifurcation until before its interaction with purple family for $\Gamma \in (3.771075,3.634014)$.\

Moreover, we have investigated that $\mathbf{f}_39_v$ is related to $f^{(6,2)}$, as illustrated via orange branch in Figure~\ref{bifurcation_graph_8}.\ The data of orange orbits are collected in Table \ref{data_f32_f6} and some orbits are plotted in Figure \ref{figure_B05_f32}.\ The orange branch with index 16 consists of stable periodic orbits.\

\subsection{Bifurcation graph between 5th cover of $g$, 6th cover of $\mathcal{B}_0^{\pm}$ and 7th cover of $f$}
\label{sec:5.9}

Our final result illustrated in Figure \ref{bifurcation_graph_5} indicates similar topological structure as in Figure \ref{bifurcation_graph_6} associated to $g^{(3,1)}$, $\mathcal{B}_0^{(4,1)}$ and $f^{(5,1)}$, namely we have derived connections related to $g^{(5,1)}$, $\mathcal{B}_0^{(6,1)}$ and $f^{(7,1)}$.\

\begin{figure}[t!]
			\centering
			\includegraphics[scale=0.85]{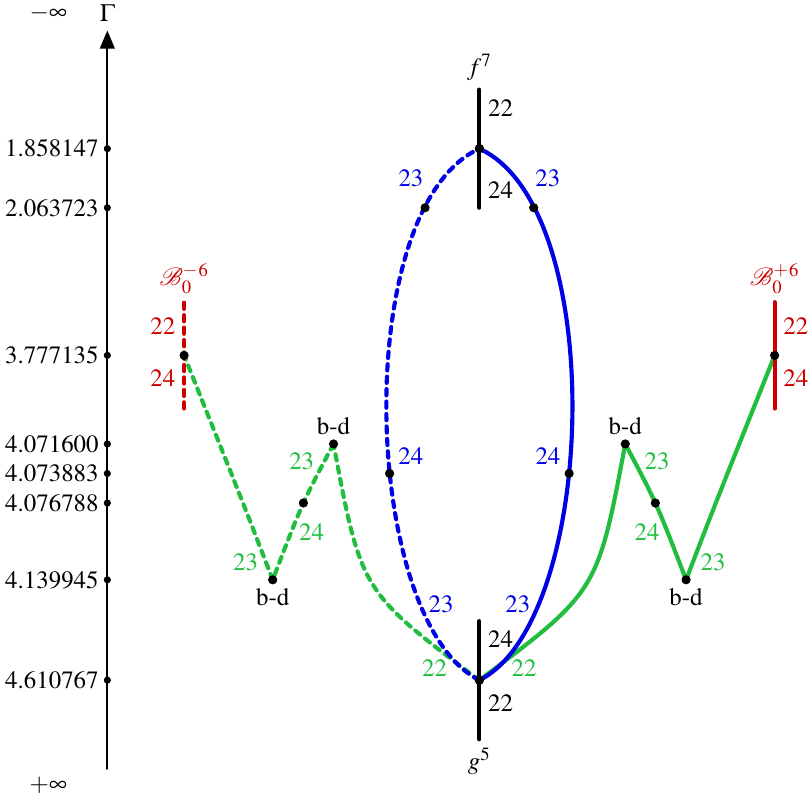}
			\caption{Bifurcation graph showing the connection between $g^{(5,1)}$ ($\Gamma = 4.610767$), $\mathcal{B}_0^{(6,1)}$ ($\Gamma = 3.777135$) and $f^{(7,1)}$ ($\Gamma = 1.858147$).\ Green orbits with index 22 are stable.}
			\label{bifurcation_graph_5}
\end{figure}

At the fifth covering of $g^{(5,1)}$ the index jumps from 22 to 24 generating two families (blue and green).\ Orbits of both branches are doubly symmetric; blue orbits with respect to the $OX$- and $OY$-axis (starting with index 23) and green orbits with respect to the $XOZ$- and $YOZ$-plane (starting with index 22).\ The blue family terminates, after two index jumps, at $f^{(7,1)}$.\ The green family interacts, after undergoing two birth-death bifurcations and one index jump, with $\mathcal{B}_0^{(6,1)}$.\ The data of blue and green orbits can be found in Table \ref{data_g5_B06_f7} and some orbits are plotted in Figure \ref{figure_g_5_f7_BO6}.\

\begin{figure}[t!]
	\centering
	\includegraphics[width=0.85\linewidth]{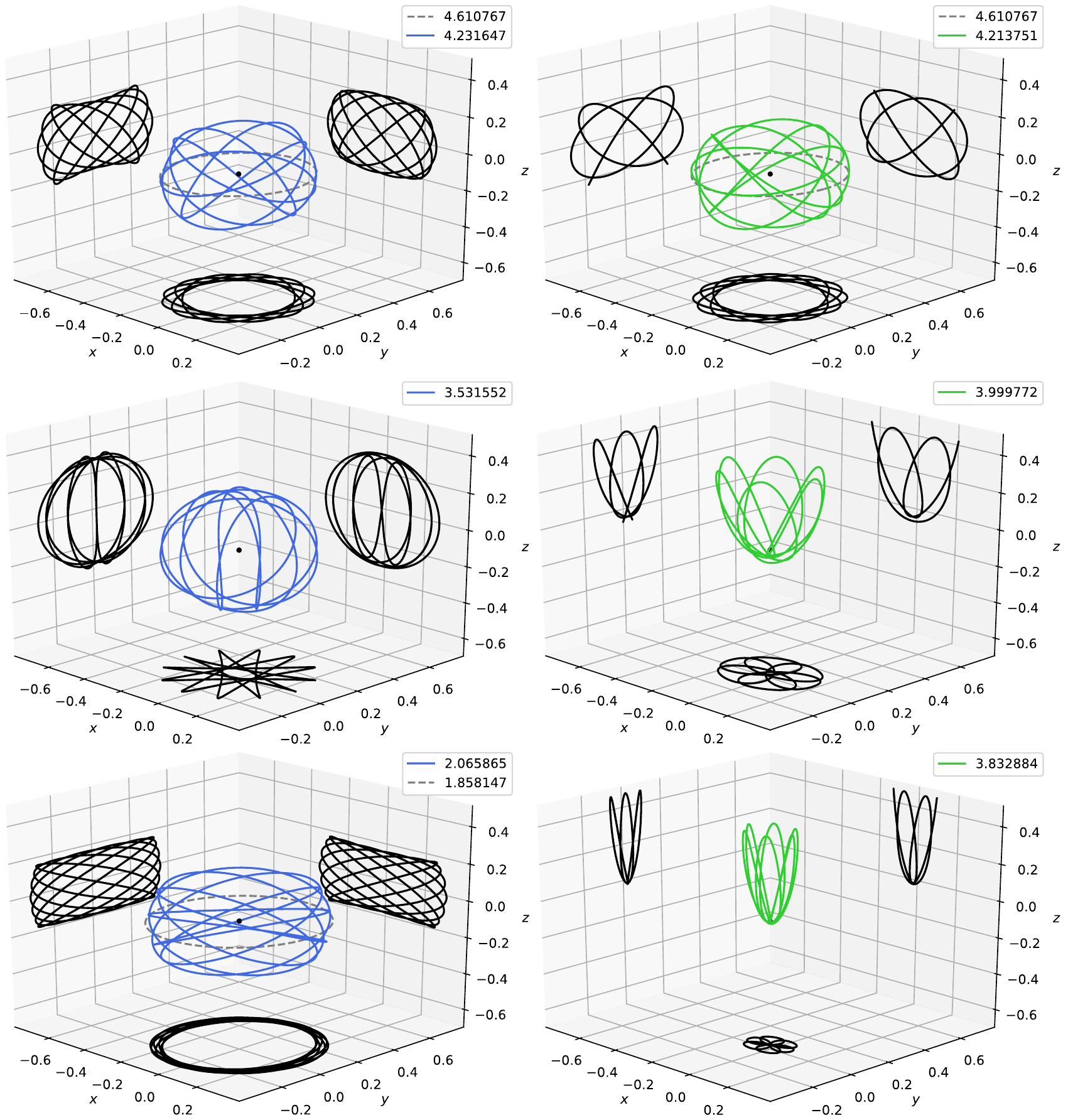}
	\caption{Plot of blue and green orbits from Figure \ref{bifurcation_graph_5}.\ Blue orbits represent connection between $g^{(5,1)}$ and $f^{(7,1)}$, and green orbits show connection among $g^{(5,1)}$ and $\mathcal B_0^{(6,1)}$.}
	\label{figure_g_5_f7_BO6}
\end{figure}

\section{Conclusion}

In this paper, we systematically explored the bifurcations of symmetric periodic orbit families in the spatial Hill three-body problem using symplectic invariants:\ Conley--Zehnder index, Krein signature, local Floer homology and its Euler characteristics.\ This analytical framework enabled us to construct bifurcation graphs that demonstrating a common network between natural symmetric periodic orbit families and reveal the intricate interconnections between these families.\ To understand such network structure, consider four factors:\vspace{-0.5em}

\begin{itemize}[noitemsep]
	\item[1)] Working in a \textbf{spatial system} ensures connections through vertical bifurcations, such as halo orbits, beyond isolated planar studies.
	\item[2)] \textbf{Regularized coordinates} are crucial, as the family $\mathcal{B}_0^{\pm}$, featuring vertical rectilinear collision orbits, cannot be fully comprehended in physical coordinates.
	\item[3)] \textbf{Symmetries of periodic solutions} are essential; they impose constraints on the spatial orbits' symmetries in interacting families.
	\item[4)] The \textbf{Conley--Zehnder index} is used to measure a winding number of the linearized flow along an orbit.\ It plays a crucial role as it provides a grading on the \textbf{local Floer homology} and its \textbf{Euler characteristics}, which serve as bifurcation invariants.\ This information is essential for understanding how different families of periodic orbits are related to each other at bifurcation points, which is not easily observable through direct computations.
\end{itemize}

The study of spatial orbit families demonstrates two main ways in which the natural families $g$, $f$, $a,c$ and $\mathcal B_0$ interact.\ The first way is an explicit interaction by means of bridge.\ Such interaction is possible if the symmetry conditions and relations between the coverings of the families involved are satisfied.\ The another way of interaction is an implicit one, when the families of spatial solutions branching from the natural families share the same periodic orbits, undergoing symmetry-breaking pitchfork bifurcation.\ Table \ref{table_1} in Introduction provides an overview of our main bifurcation results.\ As a final conlusion, in Table \ref{table_2} we outline our results in the context of explicit and implicit interactions.\

\begin{table}[bh]\centering
\begin{tcolorbox}[width=11cm]
\begin{tabular}
{P{2cm}P{2cm}P{1.5cm}P{1.5cm}P{1.2cm}}
     symmetries & explicit & implicit & symmetries & Ref.
\end{tabular}
\tcbline
\begin{tabular}
{P{2cm}P{2cm}P{1.5cm}P{1.5cm}P{1.2cm}}
     $XOZ-YOZ$ & $g \leftrightarrow \mathcal{B}_0^2$ & & & \ref{sec5.3} \\
     $XOZ$ & $g' \leftrightarrow \text{halo}^2$ & & & \ref{sec:g'_halo}
\end{tabular}
\tcbline
\begin{tabular}
{P{2cm}P{2cm}P{1.5cm}P{1.5cm}P{1.2cm}}
     $OX-XOZ$ & $g'^2 \leftrightarrow g^2$ & & & \ref{sec5.5} \\
     $OX-XOZ$ & $g'^2$ (open) & $\mathcal{B}_0^3$ & $XOZ$ & \ref{sec5.5} \\
     $OY-YOZ$ & $g^2$ (open) & $\mathcal{B}_0^3$ & $YOZ$ & \ref{sec5.5}
\end{tabular}
\tcbline
\begin{tabular}
{P{2cm}P{2cm}P{1.5cm}P{1.5cm}P{1.2cm}}
     $OX-OY$ & $g^3 \leftrightarrow f^5$ & $g'^3$ & $OX$ & \ref{sec:5.7} \\
     $XOZ-YOZ$ & $g^3 \leftrightarrow \mathcal{B}_0^4$ & $g'^3$ & $XOZ$ & \ref{sec:5.7} \\
     $XOZ-YOZ$ & $f_3 \leftrightarrow f^5$ & & & \ref{sec:5.7}
\end{tabular}
\tcbline
\begin{tabular}
{P{2cm}P{2cm}P{1.5cm}P{1.5cm}P{1.2cm}}
     $OY-YOZ$ & $g^4 \leftrightarrow f^6$ & $\mathcal{B}_0^5$ & $YOZ$ & \ref{sec:5.8} \\
     $OX-XOZ$ & $g'^4 \leftrightarrow f^6$ & $\mathcal{B}_0^5$ & $XOZ$ & \ref{sec:5.8}\\
     $OX-XOZ$ & $g'^4 \leftrightarrow g^4$ & & & \ref{sec:5.8} \\
     $OY-YOZ$ & $f_3^2 \leftrightarrow f^6$ & & & \ref{sec:5.8}
\end{tabular}
\tcbline
\begin{tabular}
{P{2cm}P{2cm}P{1.5cm}P{1.5cm}P{1.2cm}}
     $OX-OY$ & $g^5 \leftrightarrow f^7$ & & & \ref{sec:5.9} \\
     $XOZ-YOZ$ & $g^5 \leftrightarrow \mathcal{B}_0^6$ & & & \ref{sec:5.9}
\end{tabular}
\end{tcolorbox}
\caption{Bifurcation results showing explicit and implicit interactions with corresponding symmetries.\ This table is split according to our bifurcation graphs and coverings of the corresponding families.\ Families in the ''implicit`` column interact with the corresponding branch in the same row left to it.\ On the right column the corresponding subsections are referred to, where the details can be read.\ In particular, to infer the implicit connections, the Conley--Zehnder indices have played a significant role.}
\label{table_2}
\end{table}

\clearpage
\begin{appendix}
\setcounter{secnumdepth}{0}
\renewcommand{\thetable}{A.\arabic{table}}
\section{Appendix Tables of data}
	
In this appendix we collect the tables of symplectic data base (initial data, ($C/B$)-sign with corresponding Floquet multipliers, and Conley--Zehnder indices) associated to the families of periodic orbits we have studied.\ If $\lambda, 1/\lambda$ is a hyperbolic pair of Floquet multipliers lying on the real axis then we denote by $\lambda$ the one with $|\lambda| > 1$.\ For Floquet multipliers lying on the unit circle $S^1 \setminus \{ \pm 1 \}$, we denote by $\varphi$ the rotation angle modulo $2 \pi$, and in addition in case of a $k$-th root of unity we write $(\frac{2 \pi m}{k})$ next to the rotation angle, where $m \in \{ 1,...,k-1 \}$.

 \begin{table}[ht]\fontsize{10}{10}\selectfont \centering
		\caption{Data for family $g$ (1st block), $g'$ (2nd block) and $f$ (3rd block).}
\begin{tabular}{c|c|c|c|c|c}
			$\Gamma$ & $x(0)$ & $\dot{y}(0)$ & $T$ & ($C/B$)-sign \& Floquet multipliers & $\mu_{CZ}$ / $\mu_{CZ}^p$ / $\mu_{CZ}^s$\\
			\hline 6.50888000 & 0.17610000 & 2.22291184 & 0.50799 & $(+/-)$ $\varphi_p = 0.449$, $(+/-)$ $\varphi_s = 0.535$ & 6 / 3 / 3\\
            4.61076708 & 0.27482506 & 1.70093381 & 1.14495 & $(+/-)$ $\varphi_p = 0.403$, $(+/-)$ $\varphi_s = 1.256 (\frac{2\pi}{5})$ & 6 / 3 / 3\\
            4.49999001 & 0.28350000 & 1.67206473 & 1.22588 & $\varphi_p = 0$, $(+/-)$ $\varphi_s = 1.349$ & $6 \shortto 5$ / $3 \shortto 2$ / 3\\
            4.27892454 & 0.30115821 & 1.62301941 & 1.41824 & $(-/-)$ $\lambda_p = 2.265$, $(+/-)$ $\varphi_s = 1.570 (\frac{2\pi}{4})$ & 5 / 2 / 3\\
            3.87661638 & 0.32764501 & 1.59674819 & 1.88441 & $(-/-)$ $\lambda_p = 8.531$, $(+/-)$ $\varphi_s = 2.094 (\frac{2\pi}{3})$ & 5 / 2 / 3\\
            3.05747078 & 0.31084341 & 1.91481277 & 1.47270 & $(-/-)$ $\lambda_p = 154$, $(+/-)$ $\varphi_s = 3.141(\pi)$ & 5 / 2 / 3\\
            1.38309359 & 0.16471506 & 3.29248839 & 5.13725 & $(-/-)$ $\lambda_p = 1916$, $\varphi_s = 0$ & $5 \shortto 6$ / 2 / $3 \shortto 4$ \\
\end{tabular}
		\label{data_g_g'_f}
\end{table}

 \begin{table}[ht]\fontsize{10}{10}\selectfont \centering
		\caption{Data for planar Lyapunov orbits, associated to $L_2$.}
		\begin{tabular}{c|c|c|c|c|c}
			$\Gamma$ & $x(0)$ & $\dot{y}(0)$ & $T$ & ($C/B$)-sign \& Floquet multipliers & $\mu_{CZ}$ / $\mu_{CZ}^p$ / $\mu_{CZ}^s$\\
            \hline 4.32674461 &	0.69301979 & 0.00226843 & 1.51650 & $(-/-)$ $\lambda_p = 2013$, $(-/+)$ $\varphi_s = 6.066$ & 3 / 2 / 1\\
			4.29958936 & 0.66424043 & 0.18712196 & 1.51842 & $(-/-)$ $\lambda_p = 1988$, $(-/+)$ $\varphi_s = 6.075$ & 3 / 2 / 1\\
			4.19914612 & 0.62673486 & 0.41277806 & 1.52572 & $(-/-)$ $\lambda_p = 1896$, $(-/+)$ $\varphi_s = 6.113$ & 3 / 2 / 1\\
			4.00531266 & 0.58126467 & 0.67012429 & 1.54072 & $(-/-)$ $\lambda_p = 1729$, $\varphi_s = 0$ & $3 \shortto 4$ / 2 / $1 \shortto 2$ \\
			3.53719088 & 0.50138826 & 1.09813675 & 1.58282 & $(-/-)$ $\lambda_p = 1368$, $(+/+)$ $\lambda_s = 1.310$ & 4 / 2 / 2 \\
			2.50801682 & 0.36280204 & 1.84377527 & 1.71899 & $(-/-)$ $\lambda_p = 771$, $(+/+)$ $\lambda_s = 1.610$ & 4 / 2 / 2 \\
			1.27025055 & 0.21731494 & 2.84159454 & 2.04661 & $(-/-)$ $\lambda_p = 344$, $(+/+)$ $\lambda_s = 1.203$ & 4 / 2 / 2 \\
            1.22806326 & 0.21266090 & 2.88309804 & 2.06338 & $(-/-)$ $\lambda_p = 334$, $\varphi_s = 0$ & $4 \shortto 5$ / 2 / $2 \shortto 3$ \\
			0.35543459 & 0.12430197 & 3.97250147 & 2.54074 & $(-/-)$ $\lambda_p = 183$, $(+/-)$ $\varphi_s = 1.570 (\frac{2\pi}{4})$ & 5 / 2 / 3 \\
			0.14728017 & 0.10657722 & 4.31885774 & 2.69141 & $(-/-)$ $\lambda_p = 136$, $(+/-)$ $\varphi_s = 2.094 (\frac{2\pi}{3})$ & 5 / 2 / 3 \\
			$-$0.0293890 & 0.09298784 & 4.64365350 & 2.82554 & $(-/-)$ $\lambda_p = 151$, $(+/+)$ $\lambda_s = -1.005$ & 5 / 2 / 3 \\
			$-$1.0101392 & 0.04351254 & 6.85416494 & 3.47956 & $(-/-)$ $\lambda_p = 161$, $(+/+)$ $\lambda_s = -7.287$ & 5 / 2 / 3 \\
			$-$5.0382072 & 0.00668157 & 17.4461722 & 4.16967 & $(-/-)$ $\lambda_p = 660$, $(+/+)$ $\lambda_s = -19.93$ & 5 / 2 / 3 \\
			$-$20.072597 & 0.00084643 & 48.8151941 & 4.36630 & $(-/-)$ $\lambda_p = 4470$, $(+/+)$ $\lambda_s = -104.7$ & 5 / 2 / 3 
		\end{tabular}
		\label{data_a}
\end{table}

\begin{table}[ht]\fontsize{10}{10}\selectfont \centering
		\caption{Data for vertical Lyapunov orbits, associated to $L_2$.}
		\begin{tabular}{c|c|c|c|c|c|c}
			$\Gamma$ & $x(0)$ & $\dot{y}(0)$ & $\dot{z}(0)$ & $T/4$ & ($C/B$)-sign \& Floquet multipliers & $\mu_{CZ}$\\
			\hline 4.32674501 & 0.69336109 & $-$0.0000003 & 0.00192343 & 0.78539 & $(-/-)$ $\lambda = 2643$, $(+/-)$ $\varphi = 0.224$ & 5 \\
			4.00068753 & 0.67716025 & $-$0.0337281 & 0.57212194 & 0.79637 & $(-/-)$ $\lambda = 2277$, $(+/-)$ $\varphi = 0.312$ & 5 \\
			3.44534553 & 0.64980144 & $-$0.0977783 & 0.94323002 & 0.81850 & $(-/-)$ $\lambda = 1747$, $(+/-)$ $\varphi = 0.426$ & 5 \\
			1.76907143 & 0.57526073 & $-$0.3534217 & 1.60483043 & 0.92540 & $(-/-)$ $\lambda = 741$, $(+/-)$ $\varphi = 0.610$ & 5 \\
			0.99685263 & 0.55075032 & $-$0.4979227 & 1.81565624 & 1.00336 & $(-/-)$ $\lambda = 517$, $(+/-)$ $\varphi = 0.505$ & 5 \\
			0.51340845 & 0.54017719 & $-$0.5858314 & 1.92905572 & 1.06094 & $(-/-)$ $\lambda = 434$, $(+/-)$ $\varphi = 0.026$ & 5 \\
			0.51243056 & 0.54010280 & $-$0.5857445 & 1.92940524 & 2.12213 & $(-/-)$ $\lambda = 434$, $\varphi = 0$ & $5 \shortto 4$ \\
\end{tabular}
		\label{data_vertical_lyapunov}
	\end{table}

\begin{table}[ht]\fontsize{10}{10}\selectfont \centering
		\caption{Data for spatial orbits connecting families of planar Lyapunov (second out-of-plane bifurcation) and vertical Lyapunov orbits.}
		\begin{tabular}{c|c|c|c|c|c|c}
			$\Gamma$ & $x(0)$ & $\dot{y}(0)$ & $\dot{z}(0)$ & $T/2$ & ($C/B$)-sign \& Floquet multipliers & $\mu_{CZ}$\\
			\hline 1.22806326 & 0.21266090 & 2.88309804 & 0 & 2.06338 & $(-/-)$ $\lambda_p = 334$, $\varphi_s = 0$ & $4 \shortto 5$ \\
            1.22802128 & 0.21266912 & 2.88287276 & 0.03142903 & 2.06338 & $(-/-)$ $\lambda_1 = 334$, $(-/-)$ $\lambda_2 = 1.080$ & 4 \\
			1.22108891 & 0.21403149 & 2.84579684 & 0.40273942 & 2.06390 & $(-/-)$ $\lambda_1 = 335$, $(-/-)$ $\lambda_2 = 1.133$ & 4 \\
			1.17618650 & 0.22307958 & 2.61172890 & 1.05706649 & 2.06728 & $(-/-)$ $\lambda_1 = 341$, $(-/-)$ $\lambda_2 = 1.337$ & 4 \\
			0.88594065 & 0.29420463 & 1.29568353 & 2.11965140 & 2.09010 & $(-/-)$ $\lambda_1 = 381$, $(-/-)$ $\lambda_2 = 1.696$ & 4 \\
			0.71338374 & 0.35510469 & 0.59978837 & 2.22200585 & 2.10452 & $(-/-)$ $\lambda_1 = 405$, $(-/-)$ $\lambda_2 = 1.592$ & 4 \\
			0.53270943 & 0.47870335 & $-$0.2813710 & 2.06241226 & 2.12031 & $(-/-)$ $\lambda_1 = 431$, $(-/-)$ $\lambda_2 = 1.183$ & 4 \\
			0.51343480 & 0.52627680 & $-$0.5220859 & 1.96091173 & 2.12204 & $(-/-)$ $\lambda_1 = 433$, $(-/-)$ $\lambda_2 = 1.039$ & 4 \\
			0.51243056 & 0.54010280 & $-$0.5857445 & 1.92940524 & 2.12213 & $(-/-)$ $\lambda = 434$, $\varphi = 0$ & $5 \shortto 4$
		\end{tabular}
		\label{data_a_2v}
	\end{table}

    \begin{table}[ht]\fontsize{10}{10}\selectfont \centering
		\caption{Data for rectilinear vertical consecutive collision oribits (family $\mathcal{B}_0^+$).}
		\begin{tabular}{c|c|c|c|c|c}
			$h = \frac{1}{2}\Gamma$ & $a$ & $T_{\text{phys}}$ & $T_{\text{reg}}$ & ($C/B$)-sign \& Floquet multipliers & $\mu_{CZ}$\\
			\hline 99.99994999 & 0.01000000 & 0.00222144 & 0.44428819 & $(+/-)$ $\varphi_1 = 0.002$, $(-/+)$ $\varphi_2 = 6.280$ & $4$ \\
            4.97999999 & 0.20000000 & 0.19814798 & 1.98345564 & $(+/-)$ $\varphi_1 = 0.216$, $(-/+)$ $\varphi_2 = 6.103$ & 4 \\
            1.88856793 & 0.49700000 & 0.74773069 & 3.05324244 & $(+/-)$ $\varphi_1 = 0.485$, $(-/+)$ $\varphi_2 = 5.234 (\frac{2\pi 5}{6})$ & 4 \\
            1.68244383 & 0.54600000 & 0.85050489 & 3.17538556 & $(+/-)$ $\varphi_1 = 0.505$, $(-/+)$ $\varphi_2 = 5.026 (\frac{2\pi 4}{5})$ & 4 \\
            1.45659562 & 0.60900000 & 0.98340707 & 3.31417136 & $(+/-)$ $\varphi_1 = 0.509$, $(-/+)$ $\varphi_2 = 4.715 (\frac{2\pi 3}{4})$ & 4 \\
            1.20565470 & 0.69200000 & 1.15671788 & 3.46760468 & $(+/-)$ $\varphi_1 = 0.464$, $(-/+)$ $\varphi_2 = 4.183 (\frac{2\pi 2}{3})$ & 4 \\
            1.02698947 & 0.76000000 & 1.29470522 & 3.57019296 & $(+/-)$ $\varphi_1 = 0.365$, $(-/+)$ $\varphi_2 = 3.244$ & 4 \\
\end{tabular}
		\label{data_B0}
	\end{table}

	\begin{table}[ht]\fontsize{10}{10}\selectfont \centering
		\caption{Data for halo orbits; connection between planar Lyapunov orbits and $\mathcal{B}_0^{\pm}.$}
		\begin{tabular}{c|c|c|c|c|c|c}
			$\Gamma$ & $x(0)$ & $z(0)$ & $\dot{y}(0)$ & $T/2$ & ($C/B$)-sign \& Floquet multipliers & $\mu_{CZ}$\\
			\hline 4.00531266 & 0.58126467 & 0 & 0.67012429 & 1.54072 & $(-/-)$ $\lambda_p = 1729$, $\varphi_s = 0$ & $3 \shortto 4$ \\
            4.00515953 & 0.58123897 & 0.00401081 & 0.67021204 & 1.54072 & $(-/-)$ $\lambda = 1728$, $(-/+)$ $\varphi = 6.277$ & 3 \\
			3.37457564 & 0.47879236 & 0.24044127 & 0.99410946 & 1.52904 & $(-/-)$ $\lambda = 936$, $(-/+)$ $\varphi = 5.636$ & 3 \\
			2.39415881 & 0.31610954 & 0.33704920 & 1.45608154 & 1.49250 & $(-/-)$ $\lambda = 288$, $(-/+)$ $\varphi = 4.718$ & 3 \\
			1.56135153 & 0.15584832 & 0.33687541 & 1.94583333 & 1.40672 & $(-/-)$ $\lambda = 65.81$, $(-/+)$ $\varphi = 3.596$ & 3 \\
            1.35081531 & 0.10748085 & 0.31650269 & 2.13708475 & 1.35823 & $(-/-)$ $\lambda = 36.28$, $(-/+)$ $\varphi = 3.177$ & 3 \\
            1.32815770 & 0.10188012 & 0.31327120 & 2.16242790 & 1.35134 & $(-/-)$ $\lambda_1 = 33.52$, $(+/+)$ $\lambda_2 = -1.006$ & 3 \\
\end{tabular}
		\label{data_halo}
	\end{table}

    \begin{table}[ht]\fontsize{10}{10}\selectfont \centering
			\caption{Data for spatial orbits relating $g^{(1,1)}$ to $\mathcal{B}_0^{(2,2)}$.}
			\begin{tabular}{c|c|c|c|c|c|c}
				$\Gamma$ & $y(0)$ & $z(0)$ & $\dot{x}(0)$ & $T/2$ & ($C/B$)-sign \& Floquet multipliers & $\mu_{CZ}$\\
				\hline 1.38309359 &	$-$1.2529662 & 0 &	$-$0.461647 & 2.56862 & $(-/-)$ $\lambda_p = 1916$, $\varphi_s = 0$ & $5 \shortto 6 $ \\
				1.38309305 & $-$1.2529661 & 0.00071779 & $-$0.461647 & 2.56862 & 
				$(-/-)$ $\lambda = 1915$, $(-/+)$ $\varphi = 6.284$ & 5 \\
				1.28937954 & $-$1.2249340 & 0.30063425 & $-$0.453785 & 2.55455 &
				$(-/-)$ $\lambda = 1729$, $(-/+)$ $\varphi = 5.108$ & 5 \\
				0.70296101 & $-$1.0086006 & 0.82343385 & $-$0.393753 & 2.43480 &
				$(-/-)$ $\lambda = 750$, $(+/-)$ $\varphi = 3.140$ & 5 \\
				0.11441123 & $-$0.6259848 & 1.15611372 & $-$0.265029 & 2.25850 & $(-/-)$ $\lambda = 91.07$, $(+/-)$ $\varphi = 2.093$ & 5 \\
				$-$0.0724601 & $-$0.3894917 &	1.25238451 & $-$0.170009 &	2.19137 &
				$(-/-)$ $\lambda = 1.508$, $(+/-)$ $\varphi = 1.238$ & 5 \\
				$-$0.0729364 & $-$0.3886668 &	1.25262752 & $-$0.169662 &	2.19119 &
				$(-/+)$ $\varphi_1 = 6.164$, $(+/-)$ $\varphi_2 = 1.191$ & 4 \\
\end{tabular}
			\label{data_g_double_B0}
		\end{table}

\begin{table}[ht]\fontsize{10}{10}\selectfont \centering
		\caption{Data for spatial orbits connecting $g'^{(1,1)}$ with $\text{halo}^{(2,1)}$.}
		\begin{tabular}{c|c|c|c|c|c|c}
			$\Gamma$ & $x(0)$ & $z(0)$ & $\dot{y}(0)$ & $T/2$ & ($C/B$)-sign \& Floquet multipliers & $\mu_{CZ}$\\
			\hline 3.39015957 & 0.46469701 & 0 & 1.24961994 & 1.81528 & $(+/+)$ $\lambda_p = -313$, $\varphi_s = 0$ & $6 \shortto 7 $ \\
			2.44399049 & 0.30902151 & 0.25094832 & 1.67439737 & 1.92782 & $(+/+)$ $\lambda = -235$, $(-/+)$ $\varphi = 4.974$ & 6 \\
			2.00005693 & 0.22711359 & 0.27473750 & 1.92093361 & 1.99454 & $(+/+)$ $\lambda = -174$, $(-/+)$ $\varphi = 4.512$ & 6 \\
			1.44481307 & 0.11024344 & 0.26163775 & 2.35956487 & 2.08186 & $(+/+)$ $\lambda = -79.24$, $(-/+)$ $\varphi = 3.854$ & 6 \\
			1.12698114 & 0.02171818 & 0.19957730 & 2.96596840 & 2.12652 & $(+/+)$ $\lambda = -5.956$, $(-/+)$ $\varphi = 3.235$ & 6 \\
			1.12675701 & 0.02162326 & 0.19946262 & 2.96705177 & 2.12655 & $(+/+)$ $\lambda = -5.891$, $(-/+)$ $\varphi = 3.197$ & 6 \\
			1.12653383 & 0.02152852 & 0.19934800 & 2.96813498 & 2.12657 & $(+/+)$ $\lambda_1 = -5.82$, $(-/-)$ $\lambda_2 = -1.04$ & 6 \\
\end{tabular}
		\label{data_g'_double_halo}
	\end{table}

\begin{table}[ht]\fontsize{10}{10}\selectfont \centering
		\caption{First two blocks:\ Data for two families of spatial orbits branching out from $g^{(2,1)}$.\ Orbits in the first block terminate at $g'^{(2,1)}$, and orbits in the second block interact at index jump ($10 \shortto 11$) with orbits branching out from $\mathcal{B}_0^{(3,1)}$ (see first block in Table \ref{data_B03}).\ Third block:\ Data for family branching out from $g'^{(2,2)}$, whose orbits interact at first index jump ($11 \shortto 12$) with orbits branching out from $\mathcal{B}_0^{(3,1)}$ (see second block in Table \ref{data_B03}).}
\begin{tabular}{c|c|c|c|c|c|c}
			$\Gamma$ & $x(0)$ & $\dot{y}(0)$ & $\dot{z}(0)$ & $T/4$ & ($C/B$)-sign \& Floquet multipliers & $\mu_{CZ}$\\
			\hline 3.05747078 & 0.31084341 & 1.91481277 & 0 & 		1.47270 & $(-/-)$ $\lambda_p = 23815$, $\varphi_s = 0$ & $9 \shortto 11 $ \\
            3.05716289 & 0.31087152 & 1.91410150 & 0.05001676 & 1.47272 & $(-/-)$ $\lambda = 23812$, $(-/+)$ $\varphi = 6.283$ & 9 \\
            3.03212558 & 0.31338563 & 1.85273226 & 0.46022029 & 1.47408 & $(-/-)$ $\lambda = 21671$, $(-/+)$ $\varphi = 6.253$ & 9 \\
            2.96788006 & 0.31997032 & 1.69864936 & 0.83930502 & 1.47753 & $(-/-)$ $\lambda = 16844$, $(-/+)$ $\varphi = 6.175$ & 9 \\
            2.51039425 & 0.38439023 & 0.66230773 & 1.64233599 & 1.48011 & $(-/-)$ $\lambda = 717$, $(-/+)$ $\varphi = 5.519$ & 9 \\
            2.48109234 & 0.40172637 & 0.52122149 & 1.64617771 & 1.45914 & $(-/-)$ $\lambda = 178$, $(-/+)$ $\varphi = 6.191$ & 9 \\
            & & & & & & b-d \\
\end{tabular}
		\label{data_g2_g'2}
	\end{table}

\begin{table}[ht]\fontsize{10}{10}\selectfont \centering
		\caption{Data for two families branching out from $\mathcal{B}_0^{(3,1)}$ (first row).\ Orbits in the first block terminate at first index jump in third block in Table \ref{data_g2_g'2} and orbits in the second block terminate at index jump in second block in Table \ref{data_g2_g'2}.}
		\begin{tabular}{c|c|c|c|c|c|c}
			$\Gamma$ & $x(0)$ & $z(0)$ & $\dot{y}(0)$ & $T/2$ & ($C/B$)-sign \& Floquet multipliers & $\mu_{CZ}$\\
			\hline 2.41130941 & 0 & 0.69200000 & 0 & 1.73507 & $(+/-)$ $\varphi_1 = 1.394$, $\varphi_2 = 0$ & $12 \shortto 10 $ \\
            \hline 2.41752499 & 0.02588751 & 0.69087817 & $-$0.003022 & 1.73353 & $(+/-)$ $\varphi = 1.391$, $(-/-)$ $\lambda = 1.006$ & 11 \\
            2.42547679 & 0.05559667 & 0.68943016 & $-$0.006474 & 1.73603 & $(+/-)$ $\varphi = 1.412$, $(-/-)$ $\lambda = 1.003$ & 11\\
            2.55016039 & 0.20244802 & 0.66638662 & $-$0.019746 & 1.77947 & $(+/-)$ $\varphi = 1.618$, $(-/-)$ $\lambda = 1.078$ & 11 \\
            2.76400139 & 0.32243638 & 0.62523161 & $-$0.001058 & 1.88212 & $(+/-)$ $\varphi = 1.982$, $(-/-)$ $\lambda = 1.387$ & 11 \\
            2.77077169 & 0.32542301 & 0.62388324 & 0.00041078 & 1.88629 & $(+/-)$ $\varphi = 1.994$, $(-/-)$ $\lambda = 1.403$ & 11 \\
            2.97595347 & 0.40524984 & 0.57871202 & 0.11268542 & 2.10610 & $(+/-)$ $\varphi = 2.979$, $(-/-)$ $\lambda = 2.179$ & 11 \\
\end{tabular}
		\label{data_B03}
	\end{table}

\begin{table}[ht]\fontsize{10}{10}\selectfont \centering
		\caption{Data for family $\mathbf{f}_3$.}
		\begin{tabular}{c|c|c|c|c|c}
			$\Gamma$ & $x(0)$ & $\dot{y}(0)$ & $T/2$ & ($C/B$)-sign \& Floquet multipliers & $\mu_{CZ}$ / $\mu_{CZ}^p$ / $\mu_{CZ}^s$\\
			\hline $-$13.41209 & $-$4.8957388 & 9.25880043 & 14.8143 & $(-/-)$ $\lambda_p = 2.4 \cdot 10^{13}$, $(+/+)$ $\lambda_s = 16683$ & 16 / 6 / 10\\
            $-$8.009081 & $-$3.8026329 & 7.20521228 & 14.5818 & $(-/-)$ $\lambda_p = 1.3 \cdot 10^{12}$, $(+/+)$ $\lambda_s = 5121$ & 16 / 6 / 10\\
            $-$2.026588 & $-$2.0411151 & 3.93762588 & 13.1814 & $(-/-)$ $\lambda_p = 2.7 \cdot 10^{9}$, $(+/+)$ $\lambda_s = 1090$ & 16 / 6 / 10\\
            0.26914524 & $-$1.1181292 & 2.29569070 & 9.05739 & $(-/-)$ $\lambda_p = 1.3 \cdot 10^{8}$, $(+/+)$ $\lambda_s = 1.802$ & 16 / 6 / 10\\
            0.27322105 & $-$1.1170102 &	2.29355807 & 9.04601 & $(-/-)$ $\lambda_p = 1.3 \cdot 10^{8}$, $\varphi_s = 0$ & $ 15 \shortto 16$ / 6 / $9 \shortto 10$\\
            0.28326989 & $-$1.1142676 & 2.28831996 & 9.01799 & $(-/-)$ $\lambda_p = 1.3 \cdot 10^{8}$, $(-/+)$ $\varphi_s = 5.351$ & 15 / 6 / 9\\
            0.36886730 & $-$1.0918350 & 2.24482130 & 8.78148 & $(-/-)$ $\lambda_p = 1.4 \cdot 10^{8}$, $(-/+)$ $\varphi_s = 3.259$ & 15 / 6 / 9\\
\end{tabular}
		\label{data_f_3}
\end{table}

\begin{table}[ht]\fontsize{10}{10}\selectfont \centering
		\caption{First and second block:\ Data for two families of spatial orbits branching out from $g^{(3,1)}$ (1st row).\ Orbits in the first block terminate at $f^{(5,1)}$, and orbits in the second block terminate at $\mathcal{B}_0^{(4,1)}$.\ Third and fourth block:\ Data for two families of spatial orbits bifurcation from $g'^{(3,1)}$.\ Orbits in the third/fourth block interact with family from first/second block at resp.\ first index jumps.}
		\begin{tabular}{c|c|c|c|c|c|c}
			$\Gamma$ & $x(0)$ & $\dot{y}(0)$ & $\dot{z}(0)$ & $T/2$ & ($C/B$)-sign \& Floquet multipliers & $\mu_{CZ}$\\
			\hline 3.87661638 & 0.32764501 & 1.59674819 & 0 & 2.82662 & $(-/-)$ $\lambda_p = 621$, $\varphi_s = 0$ & $13 \shortto 15 $ \\
            \hline 3.87640420 & 0.32765275 & 1.59639061 & 0.03500000 & 	2.82665 & $(-/-)$ $\lambda_1 = 620$, $(-/-)$ $\lambda_2 = 1.005$ & 14 \\
            3.86641314 & 0.32801741 & 1.57957157 & 0.24200000 & 	2.82814 & $(-/-)$ $\lambda_1 = 591$, $(-/-)$ $\lambda_2 = 1.007$ & 14 \\
            3.69242900 & 0.33442700 & 1.29263199 & 0.97599999 & 2.85507 & $(-/-)$ $\lambda_1 = 238$, $(-/-)$ $\lambda_2 = 1.024$ & 14 \\
            3.27414894 & 0.35033597 & 0.64724536 & 1.54399999 & 2.92956 & $(-/-)$ $\lambda_1 = 1.228$, $(-/-)$ $\lambda_2 = 1.052$ & 14 \\
            3.27411847 & 0.35034414 & 0.64717939 & 1.54400000 & 2.92957 & $(-/-)$ $\lambda_1 = 1.229$, $\lambda_2 = 1$ & $14 \shortto 15$ \\
            3.27283294 & 0.35038725 & 0.64531070 & 1.54499999 & 2.92982 & $(-/-)$ $\lambda = 1.230$, $(+/-)$ $\varphi = 0.334$ & 15 \\
\end{tabular}
		\label{data_g3_B04_f5}
\end{table}

\begin{table}[ht]\fontsize{10}{10}\selectfont \centering
		\caption{Data for spatial orbits connecting $\mathbf{f}_3 5_v$ with $f^{(5,1)}$.}
		\begin{tabular}{c|c|c|c|c|c|c}
			$\Gamma$ & $x(0)$ & $z(0)$ & $\dot{y}(0)$ & $T/4$ & ($C/B$)-sign \& Floquet multipliers & $\mu_{CZ}$\\
			\hline 2.14663695 & $-$0.8726486 & 0 & 1.55877618 & 2.78707 & $(-/-)$ $\lambda_p = 3.5 \times 10^8$, $\varphi_s = 0$ & $13 \shortto 14$\\
            2.14662066 & $-$0.8726474 & 0.00199999 & 1.55877716 & 2.78707 & $(-/-)$ $\lambda = 3.5 \times 10^8$, $(-/+)$ $\varphi = 6.265$ & 13 \\
            2.10513484 & $-$0.8695323 & 0.10069161 & 1.56134706 & 2.78918 & $(-/-)$ $\lambda = 3.0 \times 10^8$, $(-/+)$ $\varphi = 5.393$ & 13 \\
            1.88106637 & $-$0.8539132 & 0.25060633 & 1.57829249 & 2.80094 & $(-/-)$ $\lambda = 1.2 \times 10^8$, $(-/+)$ $\varphi = 3.997$ & 13 \\
            1.40446782 & $-$0.8292488 & 0.40020666 & 1.63414222 & 2.83390 & $(-/-)$ $\lambda = 1.2 \times 10^7$, $(+/-)$ $\varphi = 2.339$ & 13 \\
            0.88549561 & $-$0.8205325 & 0.48672003 & 1.73026092 & 2.91164 & $(-/-)$ $\lambda = 3.6 \times 10^5$, $(+/-)$ $\varphi = 0.115$ & 13 \\
            0.88252640 & $-$0.8205461 & 0.48707376 & 1.73091656 & 2.91232 & $(-/-)$ $\lambda_1 = 3.5 \times 10^5$, $(+/+)$ $\lambda_2 = 1.129$ & 12 \\
\end{tabular}
		\label{data_f3_f5}
\end{table}

\begin{table}[ht]\fontsize{10}{10}\selectfont \centering
		\caption{Data for spatial orbits connecting $\mathbf{f}_3 9_v$ with $f^{(6,2)}$.}
		\begin{tabular}{c|c|c|c|c|c|c}
			$\Gamma$ & $y(0)$ & $z(0)$ & $\dot{x}(0)$ & $T/4$ & ($C/B$)-sign \& Floquet multipliers & $\mu_{CZ}$\\
			\hline $-$1.3778408 & $-$1.9362212 & 0 & $-$1.552668 & 6.84227 & $(-/-)$ $\lambda_p = 1.188$, $\varphi_s = 0$ & $17 \shortto 19$\\
            $-$1.3787982 & $-$1.9357491 & 0.00699775 & $-$1.553040 & 6.84383 & $(-/-)$ $\lambda = 1.182$, $(-/+)$ $\varphi_ = 6.283$ & 17\\
            $-$1.3813013 & $-$1.9306848 & 0.05015015 & $-$1.553814 & 6.84540 & $(-/-)$ $\lambda = 1.156$, $(-/+)$ $\varphi_ = 6.282$ & 17\\
            $-$1.3923324 & $-$1.9028470 & 0.11526201 & $-$1.558262 & 6.85217 & $(-/-)$ $\lambda = 1.011$, $(-/+)$ $\varphi_ = 6.280$ & 17\\
            & & & & & & b-d \\
            $-$1.3923474 & $-$1.9021982 & 0.11530898 & $-$1.558377 & 6.85218 & $(-/+)$ $\varphi_1 = 6.278$, $(-/+)$ $\varphi_2 = 6.281$ & 16 \\
            $-$1.3815975 & $-$1.8931081 & 0.05155514 & $-$1.560452 & 6.84586 & $(-/+)$ $\varphi_1 = 6.200$, $(-/+)$ $\varphi_2 = 6.283$ & 16 \\
            $-$1.3789032 & $-$1.8922720 & 0.00161344 & $-$1.560714 & 6.84426 & $(-/+)$ $\varphi_1 = 6.186$, $(-/+)$ $\varphi_2 = 6.283$ & 16 \\
            $-$1.3789005 & $-$1.8922713 & 0 & $-$1.560715 & 6.84426 & $(-/+)$ $\varphi_p = 6.186$, $\varphi_s = 0$ & $16 \shortto 18$
		\end{tabular}
		\label{data_f32_f6}
\end{table}

\begin{table}[ht]\fontsize{10}{10}\selectfont \centering
		\caption{Data for families of spatial orbits connecting $g^{(4,1)}$ with $f^{(6,1)}$ (first block), $g^{(4,1)}$ with $g'^{(4,1)}$ (second block) and $g'^{(4,1)}$ with $f^{(6,1)}$ (third block).\ The bifurcation point at the first index jump in the first and third block interacts with $\mathcal{B}_0^{(5,1)}$, see Table \ref{data_g4_g'4_B05_f6_2}.}
		\begin{tabular}{c|c|c|c|c|c|c}
			$\Gamma$ & $y(0)$ & $z(0)$ & $\dot{x}(0)$ & $T/4$ & ($C/B$)-sign \& Floquet multipliers & $\mu_{CZ}$\\
			\hline 4.27892454 &	$-$0.343725 & 0 & 1.24083658 & 1.41824 & $(-/-)$ $\lambda_p = 26.3$, $\varphi_s = 0$ & $17 \shortto 19 $ \\
            4.27828277 & $-$0.343558 & $-$0.010998 & 1.24098553 & 1.41829 & $(-/-)$ $\lambda_1 = 26.1$, $(+/+)$ $\lambda_2 = 1.002$ & 18 \\
            4.15265303 & $-$0.310870 & $-$0.150810 & 1.27003203 & 1.42798 & $(-/-)$ $\lambda_1 = 8.429$, $(+/+)$ $\lambda_2 = 1.001$ & 18 \\
            4.06432421 & $-$0.287904 & $-$0.193377 & 1.29033505 & 1.43510 & $(-/-)$ $\lambda_1 = 1.223$, $(+/+)$ $\lambda_2 = 1.022$ & 18 \\
            4.06230567 & $-$0.287379 & $-$0.194210 & 1.29079786 & 1.43527 & $(+/-)$ $\varphi_1 = 0.253$, $(-/+)$ $\varphi_2 = 6.266$ & 18 \\
            3.88885478 & $-$0.242317 & $-$0.251881 & 1.33037023 & 1.45005 & $(+/-)$ $\varphi_1 = 3.132$, $(-/+)$ $\varphi_2 = 6.279$ & 18 \\
            3.63647673 & $-$0.176812 & $-$0.306318 & 1.38724038 & 1.47363 & $(-/+)$ $\varphi_1 = 6.018$, $(-/+)$ $\varphi_2 = 6.277$ & 18 \\
\end{tabular}
		\label{data_g4_g'4_B05_f6_1}
	\end{table}

\begin{table}[ht]\fontsize{10}{10}\selectfont \centering
		\caption{Data for two families branching out from $\mathcal{B}_0^{(5,1)}$.\ Orbits in the first block terminate at first index jump in the first block in Table \ref{data_g4_g'4_B05_f6_1}, and orbits in the second block terminate at first index jump in the third block in Table \ref{data_g4_g'4_B05_f6_1}.}
		\begin{tabular}{c|c|c|c|c|c|c}
            $\Gamma$ & $y(0)$ & $z(0)$ & $\dot{x}(0)$ & $T/2$ & ($C/B$)-sign \& Floquet multipliers & $\mu_{CZ}$\\
            \hline 3.36488766 & 0 & 0.54600000 & 0 & 2.12626 & $(+/-)$ $\varphi_1 = 2.525$, $\varphi_2 = 0$ & $20 \shortto 18 $ \\
            3.37479676 & 0.03562396 & 0.54319462 & $-$0.064589 & 2.13067 & $(+/-)$ $\varphi_1 = 2.545$, $(-/+)$ $\varphi_2 = 6.282$ & 18 \\
            3.43493391 & 0.09481356 & 0.52561548 & $-$0.182816 & 2.16090 & $(-/+)$ $\varphi_1 = 3.611$, $(-/+)$ $\varphi_2 = 6.163$ & 18 \\
            3.73607430 & 0.20550126 & 0.41629195 & $-$0.631374 & 2.42835 & $(+/-)$ $\varphi_1 = 2.685$, $(-/+)$ $\varphi_2 = 6.271$ & 18 \\
            3.77105931 & 0.21254635 & 0.38130612 & $-$0.815469 & 2.57879 & $(+/-)$ $\varphi_1 = 2.237$, $(-/+)$ $\varphi_2 = 6.278$ & 18 \\
            & & & & & & b-d \\
            3.77105921 & 0.21254546 & 0.38118325 & $-$0.816220 & 2.57944 & $(+/-)$ $\varphi = 2.235$, $(-/-)$ $\lambda = 1.002$ & 19 \\
\end{tabular}
		\label{data_g4_g'4_B05_f6_2}
	\end{table}

\begin{table}[ht]\fontsize{10}{10}\selectfont \centering
		\caption{Data for two families of spatial orbits branching out from $g^{(5,1)}$ (1st row).\ Orbits in the first block terminate at $f^{(7,1)}$, and orbits in the second block terminate at $\mathcal{B}_0^{(6,1)}$.}
		\begin{tabular}{c|c|c|c|c|c|c}
			$\Gamma$ & $x(0)$ & $\dot{y}(0)$ & $\dot{z}(0)$ & $T/2$ & ($C/B$)-sign \& Floquet multipliers & $\mu_{CZ}$\\
			\hline 4.61076708 & 0.27482506 & 1.70093381 & 0 & 2.86239 & $(+/-)$ $\varphi_p = 2.016$, $\varphi_s = 0$ & $22 \shortto 24 $ \\
            \hline 4.61035355 & 0.27483162 & 1.70024296 & 0.05099409 & 2.86245 & $(+/-)$ $\varphi = 2.020$, $(-/-)$ $\lambda = 1.010$ & 23 \\
            4.60442515 & 0.27492355 & 1.69021380 & 0.20062545 & 2.86335 & $(+/-)$ $\varphi = 2.078$, $(-/-)$ $\lambda = 1.012$ & 23 \\
            4.46762212 & 0.27707582 & 1.46073223 & 0.92043086 & 2.88444 & $(-/+)$ $\varphi = 3.148$, $(-/-)$ $\lambda = 1.011$ & 23 \\
            4.23164740 & 0.28093445 & 1.07379907 & 1.40398678 & 2.92351 & $(-/+)$ $\varphi = 4.618$, $(-/-)$ $\lambda = 1.010$ & 23 \\
            4.07657010 & 0.28357750 & 0.82577549 & 1.59233125 & 2.95123 & $(-/+)$ $\varphi = 6.010$, $(-/-)$ $\lambda = 1.011$ & 23 \\
            4.07119716 & 0.28367068 & 0.81727260 & 1.59771853 & 2.95223 & $(-/-)$ $\lambda_1 = 1.004$, $(-/-)$ $\lambda_2 = 1.041$ & 24 \\
\end{tabular}
		\label{data_g5_B06_f7}
	\end{table}

\end{appendix}

\clearpage

\noindent
\textbf{Acknowledgements.}\ CA acknowledges support by the Deutsche Forschungsgemeinschaft (DFG, German Research Foundation) – Project-ID 541062288.\ Furthermore, CA wishes to thank Urs Frauenfelder for valuable discussions on Conley--Zehnder indices of planar and vertical Lyapunov orbits.\ AB wants to express gratitude to Prof. Alexander Bruno for useful discussion on bifurcations of periodic solutions in Hamiltonian case.

\addcontentsline{toc}{section}{References}

\hypertarget{bibliolist}{}\printbibliography
		
\end{document}